\documentclass[11pt]{amsart}
\usepackage[T1]{fontenc}
\usepackage[utf8]{inputenc}
\usepackage[english]{babel}
\usepackage{amsmath}
\usepackage{amssymb}
\usepackage{amscd}
\begin{document}
	
\title{The Dimension Conjecture for 2-Nondegenerate Hypersurfaces}
	
\author{M. A. Stepanova}
\address{Steklov Mathematical Institute of the Russian Academy of Sciences, Moscow}
\email{step\_masha@mail.ru}
	
\thanks{This research was supported by Russian Science Foundation grant No. 24-11-00196, https://rscf.ru/project/24-11-00196/.}

	 \maketitle
	 
	\begin{abstract}
	
The dimension conjecture in CR geometry is proved for 2-nondegenerate real-analytic hypersurfaces with sign-definite Levi form. Namely, it is proved that, in the indicated class of hypersurfaces in $\mathbb{C}^{N}$, the maximum dimension of finite-dimensional Lie algebras of infinitesimal holomorphic automorphisms is strictly smaller than the corresponding dimension for nondegenerate hyperquadrics in $\mathbb{C}^{N}$.

Bibliography: 16 titles.
	
	%UDC: 517.55.
	
	
	\end{abstract}

\textbf{1. Introduction}

The dimension conjecture in CR geometry is a question with a long and rich history. In its most general formulation for real-analytic manifolds of arbitrary codimension, the conjecture turned out to be false for codimension $K>1$ (see \cite{37}). However, for hypersurfaces, i.e. for $K=1$, the conjecture remains open. We state it for this case, which historically appeared first and, in essence, served as a basis for more general conjectures (all necessary definitions are given in Section 2).
	
	\vspace{3ex}
	
\textbf{The dimension conjecture ($K=1$):} 
	
1) The dimension of the Lie algebra of infinitesimal holomorphic automorphisms of the germ of a holomorphically nondegenerate real-analytic hypersurface in $\mathbb{C}^{N}$ does not exceed $N^{2}+2N$, i.e. the dimension of the automorphism algebra of any nondegenerate hyperquadric in $\mathbb{C}^{N}$ (the defining equation for a hyperquadric is written in Section 5, formula \eqref{eq43}).
	
2) Moreover, the maximum dimension of the automorphism algebra is attained only on hyperquadrics.

	\vspace{3ex}
This conjecture has been proved for the cases $N=2$ (see \cite{1}), $N=3$ (see \cite{5} or \cite{40}), and $N=4$ (see \cite{6}).
	
	\vspace{3ex}

We note that the condition of holomorphic nondegeneracy of a hypersurface (see Definition 1 in Section 2) is a criterion for finite-dimensionality of the automorphism algebra and is therefore natural. Indeed, if one abandons finite-dimensionality, then the maximum dimension (equal to infinity) is attained on any holomorphically degenerate hypersurface -- for example, on a real hyperplane. Moreover, for holomorphically nondegenerate hypersurfaces in $\mathbb{C}^{N}$ there is an upper bound for the dimension of the automorphism algebra that depends only on $N$. Namely, it is known (see Theorem 12.3.1, \cite{41}) that automorphisms of the germ of a real-analytic holomorphically nondegenerate hypersurface at a generic point $\xi$ (i.e. at a point $\xi$ outside a certain proper analytic subset of the hypersurface; below, for brevity, we shall use the term ``generic point'') are uniquely determined by the jet of order $\delta=2(N-1)$. It is enough to obtain an estimate for the dimension of the automorphism algebra of the germ $M_{\xi}$ at a generic point $\xi$ on $M$. Hence the dimension of the automorphism algebra does not exceed the dimension of the space of jets of order $2(N-1)$. Therefore the maximum dimension is certainly attained on some nondegenerate hypersurfaces.

In addition to searching for the maximum dimension for fixed $N$, one can also ask questions about the so-called submaximal dimensions (i.e. the second largest values of the dimension), as well as about which hypersurfaces realize the submaximal value. A number of works are devoted to this circle of problems (see, for example, \cite{32}). Here the following phenomenon occurs: the difference between the maximal value of the dimension and the next value is greater than one. The first result of this type apparently appears in \cite{31}. At the same time, it is clear that it makes sense to begin the search for submaximal hypersurfaces only after the hypersurfaces with maximal dimension of the automorphism algebra have been found. In all cases studied, these are nondegenerate hyperquadrics, with arbitrary signature of the Hermitian form.
	
In the present paper we prove the conjecture formulated above in a special case -- for 2-nondegenerate hypersurfaces (see the definition in Section 2) with sign-definite Levi form. At the same time, we restrict ourselves to proving that the dimension of the automorphism algebra of such a hypersurface is smaller than the dimension of the automorphism algebra of the corresponding quadric by at least one. Most likely, this estimate can be improved (i.e. one can prove that the difference of the dimensions exceeds one), but this question lies beyond the scope of our problem. Sign-definiteness of the Levi form is important for item 3.b of the proof of Lemma 13.
	
We shall use the fact that one can associate with a 2-nondegenerate hypersurface a set of biholomorphic invariants -- the 2-type (see the definitions in Section 2).
At the same time, from the $2$-type one can choose a convenient coordinate system in which the defining equation takes a sufficiently simple form. We rely on the results of \cite{44}, where a normal form for 2-nondegenerate hypersurfaces convenient for our purposes is described. Then, using the normal form of the equation, we estimate the dimension of the automorphism algebra. 
	
	\vspace{3ex}
	
We give a brief outline of the proof of the conjecture (for details and the necessary definitions see Sections 2--5).

1) We consider a generic point on the hypersurface and carry out the estimate at this point. We note that this assumption does not restrict generality.
	
2) At a generic point a holomorphically nondegenerate manifold is $l$-nondegenerate; we consider the case $l=2$. Under these conditions the rank of the Levi form is maximal and is equal to a certain quantity, which we denote by $m$. We also restrict ourselves to the case of a sign-definite Levi form. In the space $\mathbb{C}^{N}=\mathbb{C}^{n+1}$ we introduce three groups of variables: $z=(z_{1},..,z_{m}), \, y=(y_{1},...,y_{n-m}), \, w=u+iv$. Using the result of \cite{44}, we write the normal form for a 2-nondegenerate hypersurface for this value of $m$. In order to estimate the dimension of the stabilizer of the automorphism algebra, it is enough to consider the model surface (in the terminology of 
\cite{44}), which is what we do. Next we write out the $2$-model surface $\mathcal{Q}$ for the model surface $M$ (see Section 3).
	
3) We introduce suitable weights (then the $2$-model surface becomes the sum of two weighted homogeneous polynomials of weights two and three, respectively).

4) We specify the coefficients of vector fields $X$ belonging to the stabilizer of the automorphism algebra as power series with variable coefficients. We assume that each coefficient is the sum of terms whose weight is not less than a prescribed one (for each coefficient this weight is its own). We consider the lowest weighted component of the tangency condition. More precisely, from the $2$-model surface we write down a linear operator (the homological operator) and then look for a certain subspace of its kernel that corresponds to tangent fields $X$ from the stabilizer of the automorphism algebra. Then we estimate the dimension of this subspace.
	
5) We consider those weighted components of the field $X$ (the supporting component $\hat{X}$ of the field $X$) that occur in the lowest weighted component of the tangency condition.
	
6) We divide the fields into two types according to the form of the field $\hat{X}$: we consider separately fields of the first type (the coefficient of the field $\hat{X}$ at $\frac{\partial}{\partial w}$ is zero) and of the second type (the coefficient of the field $\hat{X}$ at $\frac{\partial}{\partial w}$ is nonzero).
	
7) We prove that the coefficients of the field $\hat{X}$ do not depend on $y$.
	
8) We choose a basis of the automorphism algebra of a special form, distribute the basis elements by types, and then estimate the number of basis elements.
	
9) We obtain the estimate $n^{2}+2n+1$ for the dimension of the stabilizer (a brief outline of the argument is given at the beginning of the proof of Lemma 13).

10) Taking into account the estimate $2n+1$ for the number of basis fields not lying in the stabilizer of the automorphism algebra, we obtain the total estimate $n^{2}+4n+2$, which is strictly less than the dimension of the automorphism algebra of a nondegenerate quadric (it is equal to $n^{2}+4n+3$).

	\vspace{3ex}
	
We note that the automorphism algebra for the germ of a 2-nondegenerate hypersurface with sign-definite Levi form at a generic point (on this hypersurface) is parametrized by the automorphism algebra of any holomorphically nondegenerate hyperquadric (see Remark 19). This means that the space of parameters defining the automorphism algebra of such hypersurface is a subspace of the parameter space defining the automorphism algebra of a hyperquadric.
	
The author thanks the participants of the Vitushkin seminar, the conference ``International Conference on Subtropical Geometry of k-Nondegenerate Submanifolds in Complex Space'' in Shenzhen (November 2025), as well as I. V. Beloshapka and Yu. A. Neretin, for valuable comments and discussions.

	\vspace{3ex}
	
\textbf{2. $l$-type}
	
	\vspace{3ex}
Let $T=(T_{1},...,T_{n};T_{n+1})=(T_{1},...,T_{n};w=u+iv)=(t;w=u+iv)$ be local coordinates in some neighborhood $U$ of the origin in $\mathbb{C}^{n+1}$. (We single out the coordinate $w$ for convenience of notation: below the coordinates $t=(T_{1},...,T_{n})$ will be coordinates in the complex tangent space of the hypersurface at some point of it, and the coordinate $w$ will be transversal to the complex tangent space.)
Let $M$ be a real-analytic hypersurface in $\mathbb{C}^{n+1}$ given in the neighborhood $U$ of the origin by the equation
	
	\begin{equation}\label{eq20}
	\rho(T,\overline{T})=0.
	\end{equation}

We give the necessary definitions (see also \cite{41}).
	
	\vspace{3ex}
	
\textbf{Definition 1.} $M$ is called \textit{holomorphically nondegenerate} if there is no germ of a nonzero holomorphic vector field tangent to $M$. 
	
	\vspace{3ex}
	
Here, by a germ at the point $\xi$ of a holomorphic vector field we mean a field of type $(1,0)$ of the form
	
	$$\mathcal{X}=f_{1}(T)\frac{\partial}{\partial T_{1}}+...+f_{n+1}(T)\frac{\partial}{\partial T_{n+1}},$$ 
where $f_{j}$ are germs at the point $\xi$ of holomorphic functions. The tangency condition has the form $\mathcal{X}\rho=0$ when equality \eqref{eq20} holds. We note that for a connected real-analytic hypersurface, holomorphic nondegeneracy at one point implies holomorphic nondegeneracy at an arbitrary point on $M$ (see \cite{41}).

	\vspace{3ex}
	
\textbf{Definition 2.} 
A \textit{CR-field} on $M$ is a vector field tangent to $M$ of the form
		
		$$f_{1}(T,\overline{T})\frac{\partial}{\partial \overline{T}_{1}}+...+f_{n+1}(T,\overline{T})\frac{\partial}{\partial \overline{T}_{n+1}}$$
		(a field of type $(0,1)$),
		where $f_{j}$ are smooth complex-valued functions.
	
	\vspace{3ex}
	
Let $L_{1},...,L_{n}$ be a basis of the space of CR-fields in a neighborhood of the point $\xi$ on $M$ (i.e. a set of $n$ CR-fields that are linearly independent everywhere in some neighborhood of the point $\xi$), $\frac{\partial \rho}{\partial T}=(\frac{\partial \rho}{\partial T_{1}},...,\frac{\partial \rho}{\partial T_{n+1}})$, $L^{\alpha}=L_{1}^{\alpha_{1}}\cdot...\cdot L_{n}^{\alpha_{n}},$ where $\alpha_{j}$ are nonnegative integers, and $|\alpha|=\alpha_{1}+...+\alpha_{n}$. 
	
	\vspace{3ex}
	
\textbf{Definition 3.} 
The hypersurface $M$ is called \textit{$l$-nondegenerate} at the point $\xi$ if the complex linear span of the system of vectors $\{L^{\alpha}\Big(\frac{\partial \rho}{\partial T}\Big)(\xi,\bar{\xi}), \, |\alpha|\leq l\}$ coincides with $\mathbb{C}^{n+1}$.
	\vspace{3ex}
	
At the same time, for an $l$-nondegenerate hypersurface $M$ one can introduce a finer characteristic as follows. Consider the sequence of spaces $E(\nu)={\rm span_{\mathbb{C}}}\{L^{\alpha}\Big(\frac{\partial \rho}{\partial T}\Big)(\xi,\bar{\xi}), \, |\alpha|\leq \nu\}$, where ${\rm span_{\mathbb{C}}}\{...\}$ denotes the complex linear span. It is clear that the inclusion $E(0)\subset E(1)\subset E(2) \subset ... \subset E(l) = \mathbb{C}^{n+1}$ holds.

Denote by $\sigma_{j}, \, 1\leq j \leq l,$ those numbers $\nu\geq 1$ for which a jump in the dimension of the spaces $E(\nu)$ occurs, i.e. ${\rm dim} \, E(\sigma_{j})-{\rm dim} \, E(\sigma_{j-1})>0$, where $\sigma_{0}=0$. At the same time, $E(0)$ is generated by the vector $\frac{\partial \rho}{\partial T}(\xi,\bar{\xi})$, and therefore ${\rm dim} \, E(0)=1$. We denote the sizes of the jumps ${\rm dim} \, E(\sigma_{j})-{\rm dim} \, E(\sigma_{j-1})$ by $k_{j}$. Suppose that in total $s$ jumps occurred at the values of $\nu$ equal to $\sigma_{1},...,\sigma_{s}$. We call the set of numbers $((\sigma_{1},k_{1}),...,(\sigma_{s},k_{s}))$ the $l$-\textit{type} of the hypersurface $M$ at the point $\xi$.
Both the number $l$ and the $l$-type are biholomorphic invariants (see \cite{41}). We note that the $l$-type can be introduced in an analogous way also for hypersurfaces that are not $l$-nondegenerate (the only difference is that for all $\nu$ the dimension of the spaces $E(\nu)$ in this case will be less than $(n+1)$).
	\vspace{3ex}
	
Holomorphic nondegeneracy is equivalent to each of the following two conditions (see \cite{41}): 
	
1) the Lie algebra of infinitesimal holomorphic automorphisms ${\rm aut} \, M_{\xi}$ of the germ $M_{\xi}$ is finite-dimensional for any point $\xi \in M\bigcap U$ (recall that $U$ is some neighborhood of the origin, and ${\rm aut} \, M_{\xi}$ consists of fields tangent to $M_{\xi}$ that generate local one-parameter biholomorphic mappings of the germ into itself),
	
2) there exists an $l$ such that $M$ is $l$-nondegenerate at a generic point, i.e. outside a certain proper analytic subset $\Sigma$ of the set $M\bigcap U$.

In what follows we shall assume that $M$ is holomorphically nondegenerate, with $l=2$, and that the Levi form is sign-definite; we shall estimate the dimension of the automorphism algebra of the germ $M_{\xi}$. We note that for this it is enough to obtain the corresponding estimate for the germ $M_{\xi}$ at a generic point $\xi$. Therefore, without loss of generality, one may consider the germ of the hypersurface $M$ centered at $\xi$ at a generic point and assume that $\xi$ is the origin, i.e. $\xi=0$.

For a 2-nondegenerate hypersurface with nonzero Levi form, the 2-type has the form $((1,k_{1}),(2,k_{2}))$, where $k_{1}$ is the rank of the Levi form and $k_{1}+k_{2}=n$.

At a generic point, the defining equation of a 2-nondegenerate hypersurface $M$ can be brought to a certain convenient form, which we shall describe in the next section.
	
First we note that at a generic point the $2$-type has a rather special form, namely:
	
	\vspace{3ex}
	
\textbf{Lemma 4.} At a generic point the $2$-type satisfies the following condition:
the quantity $k_{1}$ takes the maximal possible value (among the values computed at all points $\xi \in M$).
	
\textbf{Proof.} The quantity $k_{1}$ is equal to the dimension of the space $E(1)$. At a generic point this dimension is maximal (since it is equal to the rank of the Levi form).

Lemma 4 is proved.

	\vspace{3ex}
If the $2$-type satisfies the condition of Lemma 4, then we shall call it \textit{minimal}. In what follows we shall assume that the $2$-type is minimal (this is a condition on the choice of a generic point).

	\vspace{3ex}
	
\textbf{3. Normal Form} 
	 
	\vspace{3ex}

Thus, suppose that at a generic point the $2$-type is equal to $((1,k_{1}),(2,k_{2}))=((1,m),(2,k))$. According to the $2$-type, one can choose a suitable coordinate system $(z_{1},...,z_{m};y_{1},...,y_{k};w=u+iv)=(t;w=u+iv)=T$ in $\mathbb{C}^{n+1}$ such that in it the defining equation is written in the form \eqref{eq39} (see below).
Here $z,y$ are vector coordinates of dimensions $m,k$, respectively, with $k=n-m$. 
	
In accordance with the results of \cite{44}, which we shall formulate in a form convenient for us, the hypersurface $M$ in some neighborhood $U$ of the origin can be given by the equation 
	
	\begin{equation}\label{eq39}
	v=P(z,y,\bar{z},\bar{y})+R(z,y,\bar{z},\bar{y},u),
	\end{equation}
where $R(z,y,\bar{z},\bar{y},u)$ is a series all of whose terms have weight greater than two for the choice of weights $[w]=2, \, [z]=1, \, [y]=0$, and $P(z,y,\bar{z},\bar{y})$ has weight two (i.e. depends quadratically on $z,\bar{z}$). At the same time, the hypersurface $\tilde{M}$ given by the equation
	$$v=P(z,y,\bar{z},\bar{y}),$$
has the following properties (we list only those of them that we shall need):
	
1) it has the same rank of the Levi form (equal to $m$) and is 2-nondegenerate in some neighborhood $U$ of the origin,
	
2) $P(z,y,\bar{z},\bar{y})=\langle z,\bar{z}\rangle+2 \, {\rm Re} \, \Big(\mathcal{P}(z)\bar{y}\Big)+...$, where $\langle z,\bar{z}\rangle$ is a nondegenerate Hermitian form of rank $m$, $\mathcal{P}(z)$ is a vector polynomial of dimension $k$ whose coordinates are linearly independent homogeneous polynomials of degree two, and dots denote terms of degree greater than three,
	
3) all terms of bidegrees $(1,j)$ and $(j,1)$ entering $P(z,y,\bar{z},\bar{y})$ are exactly those contained in the expression $\langle z,\bar{z}\rangle+2 \, {\rm Re} \, \Big(\mathcal{P}(z)\bar{y}\Big)$ (i.e. $j$ can be equal only to one or two),

4) $P(z,y,\bar{z},\bar{y})$ contains no pluriharmonic terms,
	
5) for terms of the form $c \, z^{a}\bar{z}^{b}\bar{y}_{s}y_{j}$ in $P$ one has $|a|=|b|=1$ ($a,b$ are multi-degrees, $c\neq 0$ is a constant),
	
6) $P(z,y,\bar{z},\bar{y})$ contains no terms of the form $c \, z^{a} \bar{z}^{b} \bar{y}_{j}\bar{y}_{s}$ ($a,b$ are multi-degrees, $c\neq 0$ is a constant),
	
7) the inequality ${\rm dim \, aut} \, M_{0} \leq {\rm dim \, aut} \, \tilde{M}_{0}$ holds, where $M_{0}$ and $\tilde{M}_{0}$ are the germs at zero of the hypersurfaces $M,\tilde{M}$.

	\vspace{3ex}
If conditions 2)--6) are satisfied for the defining equation of the germ $M_{0}$, then we shall say that $M_{0}$ is given in \textit{normal form}.
	\vspace{3ex}
	
From property 7) we obtain that it is enough to get an estimate for the dimension of the automorphism algebra only for those $M$ for which $R(z,y,\bar{z},\bar{y},u)$ in their defining equation is identically zero. These are the $M$ that we shall consider below, i.e.
	
	\begin{equation}\label{eq40}
	M=\{v=\langle z,\bar{z}\rangle+2 \, {\rm Re} \, \Big(\mathcal{P}(z)\bar{y}\Big)+...\},
	\end{equation}
where dots denote terms of degree greater than three, and the germ $M_{0}$ is given in normal form. Since, by our assumption, the Levi form is sign-definite, we shall also assume that $\langle z,\bar{z}\rangle=|z_{1}|^{2}+|z_{2}|^{2}+...+|z_{m}|^{2}$ (this is easily achieved by a linear change of coordinates that does not change the general form of the defining equation).

We introduce one more convention. When analyzing the defining relation or tangency conditions, we shall need to single out in them individual monomials or groups of monomials of the form $c \, T^{\alpha} \bar{T}^{\beta}$, where $c$ is a constant and $\alpha,\beta$ are multi-degrees. We shall assume that all expressions we consider are represented as Taylor series in the variables $T, \, \overline{T}$. At the same time, we shall either assume that all such terms have been collected, or else consider an expression in which they have not been collected (when considering such monomials we specify, when necessary, which of the variants is meant). And if a certain expression contains a monomial $c \, T^{\alpha} \bar{T}^{\beta}$ as a term, then we shall say that the expression contains the monomial $c \, T^{\alpha} \bar{T}^{\beta}$ (or that the monomial $c \, T^{\alpha} \bar{T}^{\beta}$ enters the expression, or that the monomial occurs in the expression).

Next,
we denote the coordinates of the vector $\mathcal{P}(z)$ by $p_{j}(z)$.

We shall need the following lemma, which follows easily from the linear independence of the polynomials $p_{j}$.
	
	\vspace{3ex}
	
\textbf{Lemma 5.} 
There exists a linear change of coordinates such that in the new coordinates $p_{j}=\tilde{p}_{j}+...$, where $\tilde{p}_{j}, \, 1\leq j\leq k,$ are linearly independent monomials, and dots denote a sum of monomials different from $\tilde{p}_{j}$ for all $j$. Moreover, in the new coordinates the equation of the germ will also be written in normal form.

\textbf{Proof.}

Consider the polynomial $p_{1}$ and choose an arbitrary nonzero monomial $\tilde{p}_{1}$ in it. If $c \, \tilde{p}_{1}$ with $c\in \mathbb{C}, \, c\neq 0,$ enters $p_{j}, \, j\neq 1,$ then 
after the linear change $\{y_{1}=\check{y}_{1}-\bar{c}\check{y}_{j}; \,  y_{j}=\check{y}_{j}, j\neq 1\}$, instead of the term $2 \, {\rm Re} \, (p_{j}\overline{y}_{j})$ in \eqref{eq40} we obtain the term $2 \, {\rm Re} \, (\check{p}_{j}\overline{\check{y}}_{j})$, such that the monomial $\tilde{p}_{1}$ enters the polynomial $\check{p}_{j}$ with zero coefficient. Therefore one may assume that $p_{2},...,p_{k}$ do not contain the monomial $c\tilde{p}_{1}$ for $c\neq 0$.
	
Now we do the same for the polynomials $p_{2},p_{3},...,p_{k}$, choosing $\tilde{p}_{2}$. And so on: we choose $\tilde{p}_{3},...,\tilde{p}_{k}$.

Lemma 5 is proved.
	
	\vspace{3ex}

In what follows we shall assume that the coordinate system is such that the polynomials $p_{j}$ satisfy the condition formulated in Lemma 5.

Note that the monomials $p_{j}$ may fail to depend on all variables $z_{j}$. Without loss of generality one may assume that $p_{j}$ do not depend on $z_{1},...,z_{r}$
(if $p_{j}$ depend on all variables $z_{j}$, we put $r=0$). Moreover $0\leq r\leq m-1$. 
We shall call the quantity $r$ the \textit{defect} of the germ $M_{0}$.
	
	\vspace{3ex}
	
Now one can associate with the germ $M_{0}$ of the hypersurface $M$ the polynomial hypersurface
	
	$$\mathcal{Q}=\{v=|z_{1}|^{2}+|z_{2}|^{2}+...+|z_{m}|^{2}+2 \, {\rm Re} \,\Big(\mathcal{P}(z)\bar{y}\Big) \}=$$
	
	$$=\{v=\langle z,\bar{z}\rangle+ 2 \, {\rm Re} \, \mathcal{P}(z)\bar{y}\},$$	
which we shall call the \textit{$2$-model surface} of the germ $M_{0}$.

	\vspace{3ex}
	
\textbf{4. Estimate of the Dimension of the Automorphism Algebra of the Light Cone in $\mathbb{C}^{3}$}

	\vspace{3ex}

We illustrate our method with the example of the light cone $C$ in $\mathbb{C}^{3}$. The cone $C$ is the simplest (from the point of view of dimension) example of a Levi-degenerate CR-manifold with finite-dimensional automorphism algebra. (We note that the automorphism algebra of the cone is known (see, for example, \cite{34} or \cite{35}), and therefore considering this example gives no new results; however, the example of the cone will serve as a model for our considerations in the general case.) We shall discuss in detail the part of the argument for the general case that simplifies when the cone is considered, and for the remaining parts we shall restrict ourselves to a reference to the general case.
	
Let $(\zeta_{1},\zeta_{2},\zeta_{3})$ be coordinates in $\mathbb{C}^{3}$. The light cone is the tube hypersurface of the following form:
	
	$$C=\{({\rm Im} \, \zeta_{1})^{2}+({\rm Im} \, \zeta_{2})^{2}=({\rm Im} \, \zeta_{3})^{2}, \, {\rm Im} \, \zeta_{3}>0\}.$$
	
The cone $C$ is 2-nondegenerate everywhere at nonsingular points and in suitable coordinates $(z,y,w=u+iv)$ in $\mathbb{C}^{3}$ is given by the equation
	
	\begin{equation}\label{eq2}
	C=\Big\{v=\frac{|z|^{2}+{\rm Re} (z^{2}\bar{y})}{1-|y|^{2}}\Big\}
	\end{equation}
(see \cite{34}). The defining equation has the form
	
	\begin{equation}\label{eq23}
	v=|z|^{2}+{\rm Re} (z^{2}\bar{y})+o(3),
	\end{equation}
where $o(3)$ denotes terms of degree greater than three.
The cone given by equation \eqref{eq2} is holomorphically homogeneous, and therefore as a generic point at which we shall estimate the dimension of the automorphism algebra one may consider any of its points -- for example, the origin.
	
Denote by $C_{0}$ the germ at the origin of the cone $C$. Let $X=2 \, {\rm Re} \, \Big(f(z,y,w)\frac{\partial}{\partial z}+g(z,y,w)\frac{\partial}{\partial y}+h(z,y,w)\frac{\partial}{\partial w}\Big)$ be a vector field tangent to $C_{0}$. We assume that $f,g,h$ are power series whose coefficients are unknown parameters to be found. Choose the weights of the variables as follows: $[z]=[\bar{z}]=1, \, [w]=[\bar{w}]=2, \, [y]=[\bar{y}]=1$. Thus we obtain a grading on the space of power series and vector fields (if we also put $[\frac{\partial}{\partial z}]=[\frac{\partial}{\partial \bar{z}}]=-1, \, [\frac{\partial}{\partial w}]=[\frac{\partial}{\partial \bar{w}}]=-2, \, [\frac{\partial}{\partial y}]=[\frac{\partial}{\partial \bar{y}}]=-1$).  
	
Let 
	
	$$f=\sum_{j=\alpha}^{\infty}f_{j}, \ \ \ g=\sum_{j=\beta}^{\infty}g_{j}, \ \ \ h=\sum_{j=\gamma}^{\infty}h_{j},$$
where $f_{j},g_{j},h_{j}$ are components of weight $j$. That is, we assume that all components of the functions $f,g,h$ whose weight is less than $\alpha,\beta,\gamma$, respectively, are equal to zero, and the coefficients of the remaining weighted components are unknown parameters. Then we have the corresponding representations for the fields:
	
	$$\Xi_{1}=2 \, {\rm Re} \, \Big(f\frac{\partial}{\partial z}\Big)=\sum_{j=\alpha-1}^{\infty}\Xi_{1}^{(j)}, \ \ \ \Xi_{2}=2 \, {\rm Re} \, \Big(g\frac{\partial}{\partial y}\Big)=\sum_{j=\beta-1}^{\infty}\Xi_{2}^{(j)},$$
	
	$$\Xi_{3}=2 \, {\rm Re} \, \Big(h\frac{\partial}{\partial w}\Big)=\sum_{j=\gamma-2}^{\infty}\Xi_{3}^{(j)},$$ 
where $\Xi_{1}^{(j)},\Xi_{2}^{(j)},\Xi_{3}^{(j)}$ are components of weight $j$.
	
Write the tangency condition for the field $X$ to the cone $C$:
	
	\begin{equation}\label{eq6}
X(v-|z|^{2}-{\rm Re} (z^{2}\bar{y})+o(3))=0 \ \ \ \ \mbox{under the condition} \ \eqref{eq2}.
	\end{equation}
	
Denote by $\lfloor F\rfloor$ the lowest nonzero weighted component of the function $F$, and by $[F]_{s}$ the component of weight $s$ of the function $F$. 
For the tangency condition to hold, it is necessary that it hold for its lowest weighted component, namely 
	
	\begin{equation}\label{eq3}
	\lfloor X(v-|z|^{2}-{\rm Re} (z^{2}\bar{y})+o(3))\rfloor=\lfloor X(v-|z|^{2}-{\rm Re} (z^{2}\bar{y}))\rfloor=
	\end{equation}
	
	$$=\lfloor (\Xi_{1}^{\alpha-1}+\Xi_{2}^{\beta-1}+\Xi_{3}^{\gamma-2})(v-|z|^{2}-{\rm Re} (z^{2}\bar{y}))\rfloor=\lfloor{\rm Im} \, h_{\gamma}-2 \, {\rm Re} \, (f_{\alpha}\bar{z})-{\rm Re} (z^{2}\bar{g}_{\beta})\rfloor=0$$
	
$$\mbox{under the condition} \ \Big\lfloor v-\frac{|z|^{2}+{\rm Re} (z^{2}\bar{y})}{1-|y|^{2}}\Big\rfloor=v-|z|^{2}=0.$$
	
	\vspace{3ex}
	
Let $q$ be the weight of the polynomial $\lfloor{\rm Im} \, h_{\gamma}-2 \, {\rm Re} \, (f_{\alpha}\bar{z})-{\rm Re} (z^{2}\bar{g}_{\beta})\rfloor$ after substituting $w=u+i |z|^{2}$. We shall call the quantity $q$ the \textit{order} of the tangent vector field $X\in {\rm aut} \, C_{0}$. We shall call the field
	
	$$\hat{X}=2 \, {\rm Re} \, \Big([f_{\alpha}]_{q-1}\frac{\partial}{\partial z}+[g_{\beta}]_{q-2}\frac{\partial}{\partial y}+[h_{\gamma}]_{q}\frac{\partial}{\partial w}\Big)$$
	the \textit{supporting component} of $X$. 
(Here, if $\alpha=q-1$, then $[f_{\alpha}]_{q-1}=f_{\alpha}$, otherwise $[f_{\alpha}]_{q-1}=0$; similarly for the polynomials $g_{\beta},h_{\gamma}$: $[g_{\beta}]_{q-2}\neq 0$ only when $\beta=q-2$, and $[h_{\gamma}]_{q}\neq 0$ only when $\gamma=q$.) We note that 
	
	\begin{equation}\label{eq21}
	\lfloor{\rm Im} \, h_{\gamma}-2 \, {\rm Re} \, (f_{\alpha}\bar{z})-{\rm Re} (z^{2}\bar{g}_{\beta})\rfloor=\hat{X}(v-|z|^{2}-{\rm Re} (z^{2}\bar{y}))
	\end{equation}
$$ \ \ \ \ \mbox{for} \ \ \ w=u+i |z|^{2}.$$
That is, the homogeneous weighted components of the functions $f,g$ and $h$ that enter the lowest component of the tangency condition are the coefficients of the field $\hat{X}$. We also note that, in general, the supporting component of a field does not coincide with the lowest weighted component of the field. 
	
Introduce the homological operator $\Lambda$:
$$\Lambda(z,y,w,f_{\alpha},g_{\beta},h_{\gamma})=\lfloor{\rm Im} \, h_{\gamma}-2 \, {\rm Re} \, (f_{\alpha}\bar{z})-{\rm Re} (z^{2}\bar{g}_{\beta})\rfloor \ \ \ \ \mbox{for} \ \ \ w=u+i |z|^{2}.$$

Denote by ${\rm st} \, C_{0}$ the stabilizer of the automorphism algebra of the cone. It consists of fields vanishing at the origin, i.e. fields $X$ such that $f(0,0,0)=g(0,0,0)=h(0,0,0)=0$. Below we shall estimate the dimension of the stabilizer ${\rm st} \, C_{0}$, and we shall consider the operator $\Lambda$ for fields $X\in {\rm st} \, C_{0}$.
	
Denote by ${\rm \widehat{st}} \, C_{0}$ the vector space generated by all possible supporting components $\widehat{X}$ of fields $X\in{\rm st} \, C_{0}$. By Lemma 9 (see below), to estimate the dimension of the space ${\rm st} \, C_{0}$ it is enough to estimate the dimension of the space ${\rm \widehat{st}} \, C_{0}$, which lies in the kernel ${\rm ker} \,\Lambda$ of the operator $\Lambda$. For this we estimate the dimension of the space ${\rm ker} \,\Lambda$ and discard those of its elements that certainly are not supporting components of any field from ${\rm st} \, C_{0}$.

Thus we shall estimate the dimension of a certain subspace of the kernel of the homological operator $\Lambda$. This subspace corresponds to elements of the stabilizer of the automorphism algebra of the germ $C_{0}$. Recall that in the classical scheme the dimension of the family of solutions of the original equation is estimated by the dimension of the space of solutions $(f,g,h)$ of the linear equation $\Lambda(z,y,w,f,g,h)=0$ (the kernel of the homological operator $\Lambda$). However, unlike in the classical scheme, in our case the weights of the polynomials $f_{\alpha},g_{\beta},h_{\gamma}$ are different. A similar construction was discussed in \cite{6} in connection with the proof of the dimension conjecture for hypersurfaces in $\mathbb{C}^{4}$.

The fields tangent to $C$ are divided into the following types (according to the form of the supporting component of the field):
	
	1.1) $\gamma>q; \, \alpha=q-1, \, \beta>q-2$,
	
	1.2) $\gamma>q; \, \alpha>q-1, \, \beta=q-2$,
	
	1.3) $\gamma>q; \, \alpha=q-1, \, \beta=q-2$,
	
	2.1) $\gamma=q; \, \alpha>q-1, \, \beta>q-2$,
	
	2.2) $\gamma=q; \, \alpha=q-1, \, \beta>q-2$,
	
	2.3) $\gamma=q; \, \alpha>q-1, \, \beta=q-2$,
	
	2.4) $\gamma=q; \, \alpha=q-1, \, \beta=q-2$.
	
	\vspace{3ex}
	
That is, the type of a field is determined by the quantities $\gamma, \, \alpha, \, \beta$ as follows: if $\gamma=q$, then the quantity $h_{\gamma}$ enters (nontrivially) the lowest weighted component of the tangency condition (namely, the component of weight $q$); if $\alpha=q-1$, then the quantity $f_{\alpha}$ enters (nontrivially) the lowest weighted component of the tangency condition; and if $\beta=q-2$, then the quantity $g_{\beta}$ enters (nontrivially) the lowest weighted component of the tangency condition. Thus, for example, for a field of type 1.1 the component of weight $q$ is a nontrivial relation on the quantity $f_{\alpha}$.
	
From the formal point of view, type 2.4 is the most general among the types (it includes all the other types): for such fields the quantities $h_{\gamma}, \, f_{\alpha}$, and $g_{\beta}$ enter the lowest weighted component of the tangency condition. At the same time, type 1.3 includes types 1.1 and 1.2, type 2.2 includes types 1.1 and 2.1, and type 2.3 includes types 1.2 and 2.1. However, below, when estimating the number of linearly independent fields of different types, we consider the spaces generated by fields of different types and quotient spaces by the spaces of fields of the types entering them. Thus, for fields of type 1.3 we consider the quotient space of the space of fields of type 1.3 by the space of fields of types 1.1 and 1.2; for type 2.2, the quotient space by the space of fields of types 1.1 and 2.1; for type 2.3, the quotient space by the space of fields of types 1.2 and 2.1; and for type 2.4, the quotient space by the space of fields of all the remaining types.

	\vspace{3ex}
	
We shall examine in detail the case of fields of the first type whose supporting component does not depend on $w$, because for the cone the argument is simpler. The other types of fields are considered in the same way as in the general case.
	
	\vspace{3ex}

Before passing to the proofs, let us give, as an example, the decomposition by types of fields from the automorphism algebra of the cone $C$. The stabilizer of this algebra is generated as a linear space by the following fields:
	
	$$e_{1}=2 \, {\rm Re}\Big(z\frac{\partial}{\partial z}+2w\frac{\partial}{\partial w}\Big), \ \ \ e_{2}=2 \, {\rm Re}\Big(iz\frac{\partial}{\partial z}+2iy\frac{\partial}{\partial y}\Big),$$
	
	$$e_{3}=2 \, {\rm Re}\Big((z^{2}+iw+iyw)\frac{\partial}{\partial z}+2(z+zy)\frac{\partial}{\partial y}+2zw\frac{\partial}{\partial w}\Big),$$
	
	$$e_{4}=2 \, {\rm Re}\Big((iz^{2}+w+yw)\frac{\partial}{\partial z}+2i(-z+zy)\frac{\partial}{\partial y}+2izw\frac{\partial}{\partial w}\Big),$$
	
	$$e_{5}=2 \, {\rm Re}\Big(zw\frac{\partial}{\partial z}-iz^{2}\frac{\partial}{\partial y}+w^{2}\frac{\partial}{\partial w}\Big).$$

The field $e_{1}$ has type 2.2, the field $e_{2}$ has type 1.1, and the remaining fields have type 2.4.

	\vspace{3ex}

\textbf{Lemma 6.} The coefficients of the field $\hat{X}$ for $X\in {\rm st} \, C_{0}$ do not depend on $y$.
	
\textbf{Proof.} See Lemma 10 (the proof is simple and does not use the specific features of the cone case).

	\vspace{3ex}
	
We now pass to the estimate of the dimension of the automorphism algebra.

	\vspace{3ex}

\textbf{Proposition 7.} 
The number of linearly independent fields $X\in {\rm st} \, C_{0}$ of the first type whose supporting component $\hat{X}$ does not depend on $w$ does not exceed two.
	
\textbf{Proof.} 
The proof consists of considering a tree of cases of depth three. The subcases are numbered by means of sequences of letters and numbers. 
	
Next we estimate the dimension of the space of fields of types 1.1 and 1.2, and the dimension of the quotient space of fields of type 1.3 by the space of fields of types 1.1 and 1.2.
To estimate the dimensions of the quotient spaces, we shall consider a special basis of the space ${\rm aut} \, C_{0}$, which we choose in the course of the proof. In the proof of this proposition $a,b,c,e$ are nonnegative integer degrees, and $d,d_{j},d_{a}\in \mathbb{C}$. These notations are valid only within each of the numbered subcases; that is, the same parameters may take different values when different subcases are considered.

	\vspace{3ex}
	
a) Let $X$ be a field of type 1.1. By Lemma 6, $f_{\alpha}$ is a monomial of the form $d_{a}\cdot z^{a}$. Separating in condition \eqref{eq3} monomials of the form $d \cdot \bar{z}z^{a}$, we obtain that $d_{a}=0$ for $a>1$. For $a=1$ we obtain that $f_{\alpha}$ may be equal to a monomial of the form $d\cdot z$, where $d \in i\mathbb{R}$. In total there is at most one such linearly independent field.
	
	\vspace{3ex}
	
b) Next, let $X$ be a field of type 1.2. By Lemma 6, $g_{\beta}$ is a monomial of the form $d_{a}\cdot z^{a}$. Separating in condition \eqref{eq3} monomials of the form $d \cdot \bar{z}^{2}z^{a}$, we obtain that $g_{\beta}$ may be equal to a monomial of the form $d\cdot z^{2}$, where $d \in i\mathbb{R}$. In total there is at most one such linearly independent field.
	
	\vspace{3ex}
	
c) Now let $X$ be a field of type 1.3. By Lemma 6, $f_{\alpha}$ is a monomial of the form $d_{a}\cdot z^{a}$, and $g_{\beta}$ is a monomial of the form $d_{b}\cdot z^{b}$. The cases $a=1$ and $b=2$ have already been considered above (in parts a and b). Therefore, separating in condition \eqref{eq3} monomials of the form $d \cdot \bar{z}^{b}z^{a}$, we obtain that $f_{\alpha}$ may be equal to a monomial of the form $d_{1}\cdot z^{2}$ with complex coefficient $d_{1}$, and $g_{\beta}$ may be equal to a monomial of the form $d_{2}\cdot z$ with complex coefficient $d_{2}$ (i.e. $a=2$ and $b=1$). This means that there are no more than two fields of type 1.3 (since the coefficient $d_{2}$ is uniquely recovered from $d_{1}$). We note that $q=3$ in the case (c) under consideration.
	
	\vspace{3ex}
	
Thus altogether in cases a, b, and c there are no more than $1+1+2=2^{2}$ fields, i.e. no more than first-type fields for a nondegenerate quadric of CR-dimension two. However, as we shall now show, the fields described in case c do not exist.

If we consider the tangency condition for a field $X$ of type 1.3, then in the left-hand side of equality \eqref{eq6} there also occurs the term $f_{\alpha} \cdot \frac{\partial}{\partial z}(\frac{1}{2}z^{2}) \cdot \bar{y}=d_{1}\cdot z^{3} \bar{y}$ (before collecting like terms in \eqref{eq6}). 
We show that this is impossible.

	\vspace{3ex}
	
We have that the term $\sigma_{1}=-d_{1}\cdot z^{3} \bar{y}$ must also occur in the left-hand side of equality \eqref{eq6}. 
	
At the same time, $\sigma_{1}$ may occur as one of the terms in the expression $\Theta(\chi)$, where $\Theta$ is some term of the field $X$, and $\chi$ is some term in the defining relation. The terms $\Theta$ and $\chi$ must satisfy one of several conditions, which we list below.

Before passing to these conditions, let us make three remarks (we comment on these remarks in detail in the proof of Lemma 13; see items 1.c.2.1--1.c.2.3 of that proof). We number the remarks in the same way as the logical subcases.
	
	\vspace{3ex}
	
\textit{Remarks.}
	
c.1) First, since $\sigma_{1}$ is a monomial, it is enough to consider only such $\chi$ and $\Theta$ that $\chi$ is a monomial and $\Theta$ is equal to one of the following expressions: $\theta(z,y,w)\frac{\partial}{\partial z}, \, \theta(\bar{z},\bar{y},\bar{w})\frac{\partial}{\partial \bar{z}}, \,\theta(z,y,w)\frac{\partial}{\partial y}, \,\theta(\bar{z},\bar{y},\bar{w})\frac{\partial}{\partial \overline{y}}, \, \theta(z,y,w)\frac{\partial}{\partial w}$ or $\theta(\bar{z},\bar{y},\bar{w})\frac{\partial}{\partial \bar{w}}$, where $\theta$ is a monomial.
	
	\vspace{3ex}
	
c.2) Second, $\Theta$ cannot be equal to $\theta(z,y,w)\frac{\partial}{\partial w}$ or $\theta(\bar{z},\bar{y},\bar{w})\frac{\partial}{\partial \bar{w}}$.

	\vspace{3ex}

c.3) Third, we note that the coefficient of the field $\Theta$ cannot depend on the variables $w$ and $\bar{w}$.  
	
	\vspace{3ex}
	
Now we write down the conditions on $\Theta$ and $\chi$:
	
	\vspace{3ex}

c.a) The field $X$ contains a term of the form $\Theta=d_{3}z^{c}\frac{\partial}{\partial y}$, and the defining function contains a term of the form $\chi=d_{4}z^{e} y \cdot \bar{y}$. This condition cannot hold, since the defining function contains no terms of the indicated form.
	
	\vspace{3ex}
	
c.b) The field $X$ contains a term of the form $\Theta=d_{3}z^{c}\frac{\partial}{\partial z}$, and the defining function contains a term of the form $\chi=d_{4}z^{e} \cdot \bar{y}$. Then $e=2$, which is impossible, since the term $\frac{1}{2}z^{2}\bar{y}$ was already used earlier when obtaining the monomial $\sigma_{1}$ (if $\chi$ is equal to $\frac{1}{2}z^{2}\bar{y}$, then $\Theta$ enters the supporting component of the field $X$ and was also already used earlier). Therefore this case is impossible.

	\vspace{3ex}

c.c) The field $X$ contains a term of the form $\Theta=d_{3} \cdot \bar{y} \frac{\partial}{\partial \bar{z}}$, and the defining function contains a term of the form $\chi=d_{4}  \cdot z^{3} \cdot \bar{z}$. The defining function contains no terms of the indicated form, and therefore this case is impossible.
	
	\vspace{3ex}
	
c.d) The field $X$ contains a term of the form $\Theta=d_{3} \cdot \bar{y} \frac{\partial}{\partial \bar{y}}$, and the defining function contains a term of the form $\chi=d_{4}  \cdot z^{3} \cdot \bar{y}$. This case is impossible because the defining function contains no terms of the indicated form.

	\vspace{3ex}
	
Thus altogether we have no more than two linearly independent fields.
	
Proposition 7 is proved.

	\vspace{3ex}

The proof of the following lemma, which gives an estimate for fields of the second type, follows from Lemma 13 (see below).

	\vspace{3ex}

\textbf{Lemma 8.} ${\rm dim \, st} \, C_{0}\leq 9$.

	\vspace{3ex}

Since the cone is five-dimensional and homogeneous, the total estimate for the dimension of the full automorphism algebra is $5+9=14$, which is less than the corresponding dimension, equal to 15, for Levi-nondegenerate quadrics in $\mathbb{C}^{3}$.

	\vspace{3ex}
	
\textbf{5. Estimate of the Dimension of the Automorphism Algebra in the General Case}

We now pass to the consideration of the general case, in which we shall use the notation of Section 3. Below, the variables $z$ and $y$ will be vector variables, with $z\in \mathbb{C}^{m}, \, y\in \mathbb{C}^{k}$. Suppose that the generic point at which we carry out the estimate is the origin, which we denote by $\mathbf{0}$. 
	
Let 
	
	$$X=2 \, {\rm Re} \, \Big(\sum_{\nu=1}^{m}f_{\nu}(z,y,w)\frac{\partial}{\partial z_{\nu}}+\sum_{\mu=1}^{k}g_{\mu}(z,y,w)\frac{\partial}{\partial y_{\mu}}+h(z,y,w)\frac{\partial}{\partial w}\Big)$$ 
be a vector field tangent to $M_{0}$. We assume that $f_{\nu},g_{\mu},h$ are power series whose coefficients are unknown parameters to be found. In vector notation we have $X=2 \, {\rm Re} \, \Big(f(z,y,w)\frac{\partial}{\partial z}+g(z,y,w)\frac{\partial}{\partial y}+h(z,y,w)\frac{\partial}{\partial w}\Big)$, where $f=(f_{1},...,f_{m}), \, g=(g_{1},...,g_{k})$. Choose the weights of the variables as follows: $[z]=[\bar{z}]=[z_{j}]=[\bar{z}_{j}]=1, \, [w]=[\bar{w}]=2, \, [y]=[\bar{y}]=[y_{s}]=[\overline{y}_{s}]=1$ for all $j,s$. Thus we obtain a grading on the space of power series and vector fields (if we also put $[\frac{\partial}{\partial z}]=[\frac{\partial}{\partial \bar{z}}]=[\frac{\partial}{\partial z_{j}}]=[\frac{\partial}{\partial \bar{z}_{j}}]=-1, \, [\frac{\partial}{\partial w}]=[\frac{\partial}{\partial \bar{w}}]=-2, \, [\frac{\partial}{\partial y}]=[\frac{\partial}{\partial \bar{y}}]=[\frac{\partial}{\partial y_{s}}]=[\frac{\partial}{\partial \overline{y}_{s}}]=-1$).  
	
With this choice of weights the defining relation for $M$ takes the form
	
	\begin{equation}\label{eq28}
	v=|z_{1}|^{2}+|z_{2}|^{2}+...+|z_{m}|^{2}+2 \, {\rm Re} \,\Big(\mathcal{P}(z)\bar{y}\Big)+o(3),
	\end{equation}
where $o(3)$ denotes terms of degree greater than three, and these terms have the special form described in Section 3 (normal form).

Let 
	
	$$f_{\nu}=\sum_{j=\alpha_{\nu}}^{\infty}f_{\nu, j}, \ \ \ g_{\mu}=\sum_{j=\beta_{\mu}}^{\infty}g_{\mu, j}, \ \ \ h=\sum_{j=\gamma}^{\infty}h_{j},$$
where $f_{\nu, j},g_{\mu, j},h_{j}$ are components of weight $j$. 
That is, we assume that all components of the functions $f_{\nu},g_{\mu},h$ whose weight is less than $\alpha_{\nu},\beta_{\mu},\gamma$, respectively, are equal to zero, and the coefficients of the remaining weighted components are unknown parameters. Then we have the following representations for the fields:
	
	$$\Xi_{1,\nu}=2 \, {\rm Re} \, \Big(f_{\nu}\frac{\partial}{\partial z_{\nu}}\Big)=\sum_{j=\alpha_{\nu}-1}^{\infty}\Xi_{1,\nu}^{(j)}, \ \ \ \Xi_{2,\mu}=2 \, {\rm Re} \, \Big(g_{\mu}\frac{\partial}{\partial y_{\mu}}\Big)=\sum_{j=\beta_{\mu}-1}^{\infty}\Xi_{2,\mu}^{(j)},$$
	
	$$\Xi_{3}=2 \, {\rm Re} \, \Big(h\frac{\partial}{\partial w}\Big)=\sum_{j=\gamma-2}^{\infty}\Xi_{3}^{(j)},$$ 
where $\Xi_{1,\nu}^{(j)},\Xi_{2,\mu}^{(j)},\Xi_{3}^{(j)}$ are components of weight $j$.
	
Write the tangency condition for the field $X$ to the hypersurface $M$:

	\begin{equation}\label{eq36}
	X\Big(v-\langle z,\bar{z}\rangle-2 \, {\rm Re} \,
	\mathcal{P}(z)\bar{y}+o(3)\Big)=0
	\end{equation}	
$$  \ \ \ \ \mbox{under the condition} \ \eqref{eq28}.$$
	
Denote by $\lfloor F\rfloor$ the lowest nonzero weighted component of the function $F$, and by $[F]_{s}$ the component of weight $s$ of the function $F$.

Let $f_{\alpha}=(f_{1, \alpha_{1}},...,f_{m, \alpha_{m}}), \, g_{\beta}=(g_{1, \beta_{1}},...,g_{k, \beta_{k}})$ be the vectors formed from the lowest weighted components of the vector-functions $f$ and $g$. For the tangency condition to hold, it is necessary that it hold for the lowest weighted component of the tangency condition, namely 
	
	\begin{equation}\label{eq13}
	\lfloor X\Big(v-\langle z,\bar{z}\rangle-2 \, {\rm Re} \,\mathcal{P}(z)\bar{y}+o(3)\Big)\rfloor=
	\end{equation}
	$$=\lfloor X\Big(v-\langle z,\bar{z}\rangle-2 \, {\rm Re} \,\mathcal{P}(z)\bar{y}\Big)\rfloor=$$
	
	$$=\lfloor (\sum_{\nu}\Xi_{1,\nu}^{(\alpha_{\nu}-1)}+\sum_{\mu}\Xi_{2,\mu}^{(\beta_{\mu}-1)}+\Xi_{3}^{(\gamma-2)})\Big(v-\langle z,\bar{z}\rangle-2 \, {\rm Re} \,\mathcal{P}(z)\bar{y}\Big)\rfloor=$$
	
	$$=\lfloor{\rm Im} \, h_{\gamma}-2 \, {\rm Re} \, \langle f_{\alpha},\bar{z}\rangle-2 \, {\rm Re} (\mathcal{P}(z)\bar{g}_{\beta})\rfloor=0$$
	
$$\mbox{under the condition} \ \Big\lfloor v-\langle z,\bar{z}\rangle-2 \, {\rm Re} \,\mathcal{P}(z)\bar{y}+o(3)\Big\rfloor=v-\langle z,\bar{z}\rangle=0.$$
That is, the lowest weighted component of the tangency condition is formed by those terms of minimal weight in the tangency condition that give nontrivial relations on $f_{\alpha}, \, g_{\beta}, \, h$.
	\vspace{3ex}
	
Let $q$ be the weight of the lowest component of the tangency condition, i.e. the weight of the polynomial $\lfloor{\rm Im} \, h_{\gamma}-2 \, {\rm Re} \, \langle f_{\alpha},\bar{z}\rangle-2 \, {\rm Re} (\mathcal{P}(z)\bar{g}_{\beta})\rfloor$, in which $w=u+i \langle z,\bar{z}\rangle$. We shall call the quantity $q$ the \textit{order} of the tangent vector field $X\in {\rm aut} \, M_{0}$. We shall call the field
	
	$$\hat{X}=2 \, {\rm Re} \, \Big(\sum_{\nu}[f_{\nu,\alpha_{\nu}}]_{q-1}\frac{\partial}{\partial z_{\nu}}+\sum_{\mu}[g_{\mu,\beta_{\mu}}]_{q-2}\frac{\partial}{\partial y_{\mu}}+[h_{\gamma}]_{q}\frac{\partial}{\partial w}\Big)$$
	the \textit{supporting component} of $X$. 
(Here, if $\alpha_{\nu}=q-1$, then $[f_{\nu,\alpha_{\nu}}]_{q-1}=f_{\nu,\alpha_{\nu}}$, otherwise $[f_{\nu,\alpha_{\nu}}]_{q-1}=0$; similarly for the polynomials $g_{\mu,\beta_{\mu}},h_{\gamma}$: $[g_{\mu,\beta_{\mu}}]_{q-2}\neq 0$ only when $\beta_{\mu}=q-2$, and $[h_{\gamma}]_{q}\neq 0$ only when $\gamma=q$.) We note that 
	
	\begin{equation}\label{eq29}
	\lfloor{\rm Im} \, h_{\gamma}-2 \, {\rm Re} \, \langle f_{\alpha},\bar{z}\rangle-2 \, {\rm Re} (\mathcal{P}(z)\bar{g}_{\beta})\rfloor=
	\end{equation}
$$=\hat{X}\Big(v-\langle z,\bar{z}\rangle-2 \, {\rm Re} \,\mathcal{P}(z)\bar{y}\Big) \ \ \ \ \mbox{for} \ \ \ w=u+i \langle z,\bar{z}\rangle.$$
That is, the homogeneous weighted components of the functions $f_{\nu},g_{\mu}$ and $h$ that enter the lowest weighted component of the tangency condition are the coefficients of the field $\hat{X}$. We also note that, in general, the supporting component of a field does not coincide with the lowest weighted component of the field.
	
Introduce the homological operator $\Lambda$:
$$\Lambda(z,y,w,f_{\alpha},g_{\beta},h_{\gamma})=\lfloor{\rm Im} \, h_{\gamma}-2 \, {\rm Re} \, \langle f_{\alpha},\bar{z}\rangle-2 \, {\rm Re} (\mathcal{P}(z)\bar{g}_{\beta})\rfloor \ \ \ \ \mbox{for} \ \ \ w=u+i \langle z,\bar{z}\rangle.$$

Denote by ${\rm st} \, M_{0}$ the stabilizer of the automorphism algebra of the hypersurface. It consists of fields vanishing at the origin, i.e. fields $X$ such that $f_{\nu}(\mathbf{0})=g_{\mu}(\mathbf{0})=h(\mathbf{0})=0$. Below we shall estimate the dimension of the stabilizer ${\rm st} \, M_{0}$, and we shall consider the operator $\Lambda$ for fields $X\in {\rm st} \, M_{0}$.
	
As the following lemma shows, it is enough to estimate only the number of all possible linearly independent supporting components of the field (the proof follows easily from the definition).

Denote by ${\rm \widehat{st}} \, M_{0}$ the vector space generated by all possible supporting components $\widehat{X}$ of fields $X\in{\rm st} \, M_{0}$, and by ${\rm ker} \, \Lambda$ the kernel of the operator $\Lambda$.
	
	\vspace{3ex}

\textbf{Lemma 9.} 1) ${\rm dim \, st} \, M_{0}\leq{\rm dim \, \widehat{st}} \, M_{0} \leq {\rm dim \, ker} \, \Lambda$;
	
2) If the supporting components $\hat{X}_{1}$ and $\hat{X}_{2}$ are different, then the fields $X_{1}$ and $X_{2}$ themselves are different.
	
\textbf{Proof.} Item 2 is obvious; let us prove item 1. The second inequality in item 1 is also obvious; we shall prove the first inequality. By the finite-dimensionality condition of the algebra ${\rm aut} \, M_{0}$, there exists a number $q_{max}<\infty$ such that $q\leq q_{max}$ for all possible values of $q$.

Consider $\mathcal{V}_{q}\subset {\rm aut} \, M_{0}$, the space generated by fields of order at least $q$. Moreover, if $q>q_{max}$, then the space $\mathcal{V}_{q}$ is zero.

Define a homomorphism $\Psi_{q}: \mathcal{V}_{q}\longrightarrow{\rm \widehat{st}} \, M_{0}$ by the formula $\Psi_{q}(X)=[X]_{q}$, where $[X]_{q}=\hat{X}$, if $X\in \mathcal{V}_{q} \setminus \mathcal{V}_{q+1}$, and $[X]_{q}=0$, if $X\in \mathcal{V}_{q+1}$. The kernel of the homomorphism $\Psi_{q}$ is exactly the space $\mathcal{V}_{q+1}$, i.e. the image of the homomorphism $\Psi_{q}$ is isomorphic to the quotient space $\mathcal{V}_{q} \diagup\mathcal{V}_{q+1}$.
	
At the same time it is clear that the algebra ${\rm st} \, M_{0}$ is isomorphic to the direct sum 
	
	$$\bigoplus_{q=0}^{q_{max}} \mathcal{V}_{q} \diagup\mathcal{V}_{q+1}.$$
	
It follows that the space ${\rm st} \, M_{0}$ can be mapped linearly and injectively into ${\rm \widehat{st}} \, M_{0}$, which implies the required inequality for the dimension.
	
Lemma 9 is proved.

	\vspace{3ex}
Thus, below we shall estimate the dimension of the space ${\rm \widehat{st}} \, M_{0} \subseteq {\rm ker} \,\Lambda$. To do this, we estimate the dimension of the space ${\rm ker} \,\Lambda$ and discard those of its elements that certainly are not supporting components of any field from ${\rm st} \, M_{0}$.

Thus we shall estimate the dimension of a certain subspace of the kernel of the homological operator $\Lambda$. This subspace corresponds to elements of the stabilizer of the automorphism algebra of the germ $M_{0}$. However, unlike in the classical scheme, in our case the weights of the vector polynomials $f_{\alpha},g_{\beta}$ and of the polynomial $h_{\gamma}$ are different.

The fields tangent to $M$ are divided into the following types (according to the form of the supporting component of the field):
	
1.1) $\gamma>q; \, \alpha_{\nu}=q-1$ for some $\nu$, $\, \beta_{\mu}>q-2$ for all $\mu$,
	
1.2) $\gamma>q; \, \alpha_{\nu}>q-1$ for all $\nu$, $\, \beta_{\mu}=q-2$ for some $\mu$,
	
1.3) $\gamma>q; \, \alpha_{\nu}=q-1$ for some $\nu$, $\, \beta_{\mu}=q-2$ for some $\mu$,
	
2.1) $\gamma=q; \, \alpha_{\nu}>q-1$ for all $\nu$, $\, \beta_{\mu}>q-2$ for all $\mu$,
	
2.2) $\gamma=q; \, \alpha_{\nu}=q-1$ for some $\nu$, $\, \beta_{\mu}>q-2$ for all $\mu$,
	
2.3) $\gamma=q; \, \alpha_{\nu}>q-1$ for all $\nu$, $\, \beta_{\mu}=q-2$ for some $\mu$,
	
2.4) $\gamma=q; \, \alpha_{\nu}=q-1$ for some $\nu$, $\, \beta_{\mu}=q-2$ for some $\mu$.
	
	\vspace{3ex}
	
That is, the type of a field is determined by the quantities $\gamma, \, \alpha_{\nu}, \, \beta_{\mu}$ as follows: if $\gamma=q$, then the quantity $h_{\gamma}$ enters (nontrivially) the lowest weighted component of the tangency condition (namely, the component of weight $q$); if $\alpha_{\nu}=q-1$, then the quantity $f_{\alpha_{\nu}}$ enters (nontrivially) the lowest weighted component of the tangency condition; and if $\beta_{\mu}=q-2$, then the quantity $g_{\beta_{\mu}}$ enters (nontrivially) the lowest weighted component of the tangency condition. Thus, for example, for a field of type 1.1 the component of weight $q$ is a nontrivial relation on some of the quantities $f_{\alpha_{\nu}}$. We also note that the inequality sign "$<$" cannot occur in the definitions of the types (since $q$ is the weight of the lowest component of the tangency condition).
	
From the formal point of view, type 2.4 is the most general among the types (it includes all the other types): for such fields the quantities $h_{\gamma}, \, f_{\alpha_{\nu}}, \, g_{\beta_{\mu}}$ enter the lowest weighted component of the tangency condition for some $\nu, \, \mu$. At the same time, type 1.3 includes types 1.1 and 1.2, type 2.2 includes types 1.1 and 2.1, and type 2.3 includes types 1.2 and 2.1. However, below, when estimating the number of linearly independent fields of different types, we consider the spaces generated by fields of different types and quotient spaces by the spaces of fields of the types entering them. Thus, for fields of type 1.3 we consider the quotient space by fields of types 1.1 and 1.2; for type 2.2, the quotient space by fields of types 1.1 and 2.1; for type 2.3, the quotient space by fields of types 1.2 and 2.1; and for type 2.4, the quotient space by fields of all other types.
	
The conditions on the types can be written briefly in table form as follows:
	
	\begin{table}[h]
		\begin{center}
			\begin{tabular}{|c|c|c|c|}
				\hline
				Type number & $\gamma$ & $\alpha_{\nu}$ & $\beta_{\mu}$ \\
				\hline
				1.1 & > & = & > \\
				\hline
				1.2 & > & > & = \\
				\hline
				1.3 & > & = & = \\
				\hline
				2.1 & = & > & > \\
				\hline
				2.2 & = & = & > \\
				\hline
				2.3 & = & > & = \\
				\hline
				2.4 & = & = & = \\
				\hline
			\end{tabular}
		\end{center}
	\end{table} 
	
	\vspace{3ex}
	
The sign "$>$" in the table means the absence of the corresponding polynomials in the lowest weighted component of the tangency condition, and the equality sign means their presence. We note that the table has no row with three signs "$>$" (it corresponds to the zero field).  
	
	\vspace{3ex}
	
Introduce the following notation: $X_{0}=2 \, {\rm Re} \, \frac{\partial}{\partial w}, \, X_{1}=[X_{0},X], \, X_{\nu}=[X_{0},X_{\nu-1}]$ for $\nu\geq 2$. Below, the lower index of a field will be used precisely in this sense. Notice that $X_{0}\in {\rm aut} \, M_{0}$, since the defining equation of the hypersurface does not depend on $u$. Moreover, $X_{\nu}$ also belongs to ${\rm aut} \, M_{0}$. 
	
	\vspace{3ex}
	
We shall also need to distinguish additional subtypes of fields:
	
1(0) -- fields of the first type whose supporting component does not depend on $w$,
	
1.1(0) -- fields of type 1.1 whose supporting component does not depend on $w$,
	
1.2(0) -- fields of type 1.2 whose supporting component does not depend on $w$,
	
1.3(0) -- fields of type 1.3 whose supporting component does not depend on $w$,
	
1.1(1) -- fields $X$ such that $X_{1}$ has type 1.1 
	
1.1(2) -- fields $X$ such that $X_{2}$ has type 1.1 
	
$2.b(c)$ -- fields $X$ of the second type such that the coefficient at $\frac{\partial}{\partial w}$ of the field $\hat{X}_{c}$ is a homogeneous polynomial of degree $b$ depending only on the variable $z$.
	
In addition we introduce the types $2.0(\geq 3)$ (the union of fields of types $2.0(c)$ for $c\geq 3$), $2.1(\geq 3)$ (the union of fields of types $2.1(c)$ for $c\geq 3$), $2.2(\geq 2)$ (the union of fields of types $2.2(c)$ for $c\geq 2$), and $2.\geq 3(\geq 0)$ (the union of fields of types $2.b(c)$ for $b\geq 3, \, c\geq 0$).

We note that fields $X$ of type $2.1(c)$, by definition, contain terms of the form $d_{1}z_{s}w^{c}\frac{\partial}{\partial w}$ in the supporting component $\hat{X}$, and hence they also contain terms of the form $d \, \bar{w}^{c}\frac{\partial}{\partial \bar{z}_{s}}$ in $\hat{X}$. Fields $X$ of type $2.2(c)$, by definition, contain terms of the form $d_{1}z_{s_{1}}z_{s_{2}}w^{c}\frac{\partial}{\partial w}$ in $\hat{X}$, and hence they also contain terms of the form $d \, \bar{w}^{c}\frac{\partial}{\partial \bar{y}_{s}}$ in $\hat{X}$. Indeed, if $\hat{X}$ contains terms of the form $d_{1} z_{s}w^{c}\frac{\partial}{\partial w}$ or of the form $d_{1} z_{s_{1}}z_{s_{2}}w^{c}\frac{\partial}{\partial w}$, then in the lowest component of the tangency condition there occurs, respectively, either a term of the form $d_{2}z_{s}u^{c}$ or a term of the form $d_{2}z_{s_{1}}z_{s_{2}}u^{c}$ (they enter the expressions $d_{1} z_{s}w^{c}\frac{\partial}{\partial w}(\frac{w}{2 i})$ and $d_{1} z_{s_{1}}z_{s_{2}}w^{c}\frac{\partial}{\partial w}(\frac{w}{2 i})$). Hence the terms $-d_{2}z_{s}u^{c}$ and $-d_{2}z_{s_{1}}z_{s_{2}}u^{c}$, respectively, must also occur in the lowest component. These terms can be obtained only in the following ways (apart from those already indicated): the first by means of the term $z_{s}\bar{z}_{s}$ in the defining relation and a term in $\hat{X}$ of the form $d \bar{w}^{c}\frac{\partial}{\partial \bar{z}_{s}}$, and the second by means of the term $p_{j}\bar{y}_{j}$ for some $j$ in the defining relation and a term in $\hat{X}$ of the form $d \bar{w}^{c}\frac{\partial}{\partial \bar{y}_{j}}$.

At the same time, we consider fields of type $2.b(c)$ as the quotient space by fields of types $2.b'(c')$, where $b'>b, \, b+2c=b'+2c'$.
	
	\vspace{3ex}

Before passing to the proofs, let us give, as an example, the decomposition by types of fields from the automorphism algebra of a Levi-nondegenerate (i.e. 1-nondegenerate) hyperquadric. A hyperquadric $Q$ in $\mathbb{C}^{n+1}$ with coordinates $(z_{1},...,z_{n},w=u+iv)=(z,w)$ is a hypersurface given by
	\begin{equation}\label{eq43}
	v=\mathcal{H}(z,\bar{z}),
	\end{equation}
where $\mathcal{H}(z,\bar{z})$ is a Hermitian form. Accordingly, the hyperquadric is Levi-nondegenerate if $\mathcal{H}(z,\bar{z})$ is nondegenerate. 
	
From the formal point of view, we consider only Levi-degenerate hypersurfaces; however, the decomposition by types can also be considered for Levi-nondegenerate ones. 
	
The stabilizer of the automorphism algebra of the Levi-nondegenerate hyperquadric given by equation \eqref{eq43} has the form:
	
	\begin{equation}\label{eq19}
	\mathfrak{g}_{0}= \Big\{2 \, {\rm Re} \, \Big(\Omega z \frac{\partial}{\partial z}+\tau w \frac{\partial}{\partial w}\Big)\Big\},
	\end{equation}
	
	$$\mathfrak{g}_{1}=\Big\{2 \, {\rm Re} \, \Big((aw+2i\mathcal{H}(z,\bar{a})) \frac{\partial}{\partial z}+2i\mathcal{H}(z,\bar{a}) w \frac{\partial}{\partial w}\Big)\Big\},$$
	
	$$\mathfrak{g}_{2}=\Big\{2 \, {\rm Re} \, \Big(Rwz \frac{\partial}{\partial z}+Rw^{2} \frac{\partial}{\partial w}\Big)\Big\},$$
where $R\in \mathbb{R}, \, a\in\mathbb{C}^{n}, \, \Omega$ is a matrix of size $n\times n$, and $\tau$ is a real number such that $2 \, {\rm Re} \, \mathcal{H}(\Omega z,\bar{z})=\tau\mathcal{H}(z,\bar{z})$.
	
Fields from $\mathfrak{g}_{0}$ with $\tau=0$ have type 1.1, while with $\tau\neq 0$ (the quotient space by fields of type 1.1) they have type 2.2. Fields from $\mathfrak{g}_{1}$ and $\mathfrak{g}_{2}$ have type 2.2. We note that in the case of a hyperquadric $m=n$, so the variables $y_{j}$ (and the corresponding field types 1.2, 1.3, 2.3, 2.4) are absent.
	
	\vspace{3ex}
	
At the same time, the decomposition by types holds not only for hypersurfaces but also for manifolds of higher codimension. The difference is that the variable $w$ becomes vector-valued. We also note that, as the examples in \cite{29} show, the automorphism algebras of nondegenerate quadrics of codimension $K>1$ are significantly more complicated.

We shall need the following two lemmas.
	
	\vspace{3ex}
	
\textbf{Lemma 10.} The coefficients of the field $\hat{X}$ do not depend on $y$.
	
	\textbf{Proof.} Consider condition \eqref{eq13} -- the lowest weighted component of the tangency condition, i.e. the homogeneous component of weight $q$. Let $a,b,e,e_{1}$ be multi-degrees, let $c$ be a nonnegative integer, and let $d,\tilde{d},d_{1},d_{2},d_{3}\in \mathbb{C}$.

	1) Monomials of the form $\sigma=d \cdot z^{a}y^{b}u^{c}$ with $b\neq 0$ can occur in condition \eqref{eq13} only if the supporting component $\hat{X}$ of the field $X$ contains a term $\Theta=2 i d \cdot z^{a}y^{b}w^{c}\frac{\partial}{\partial w}$. Then $\sigma$ is one of the monomials in the expression $\Theta(\frac{w}{2i})$ (which enters \eqref{eq13}, since $\frac{w}{2i}$ enters \eqref{eq28}). Hence we obtain $d=0$, i.e. $\frac{\partial}{\partial y_{\lambda}}(h_{\gamma})=0$ for all $\lambda$. Therefore only terms $\Theta$ in $\hat{X}$ of the form $\theta(z,y,w)\frac{\partial}{\partial z_{\lambda}}, \, \theta(\bar{z},\bar{y},\bar{w})\frac{\partial}{\partial \bar{z}_{\lambda}}, \, \theta(z,y,w)\frac{\partial}{\partial y_{\lambda}}, \, \theta(\bar{z},\bar{y},\bar{w})\frac{\partial}{\partial \overline{y}_{\lambda}}$ can depend on the variable $y$, where $\theta$ is a monomial.

	2) Monomials of the form $\sigma=d \cdot \bar{z}_{\lambda}z^{a}y^{b}u^{c}$ with $b\neq 0$ can occur in condition \eqref{eq13} only if $\hat{X}$ contains a term $\Theta= d \cdot z^{a}y^{b}w^{c}\frac{\partial}{\partial z_{\lambda}}$. Then $\sigma$ is one of the monomials in the expression $\Theta(|z_{\lambda}|^{2})$ (which enters \eqref{eq13}, since $|z_{\lambda}|^{2}$ enters \eqref{eq28}).
	Hence we obtain $d=0$, i.e. $\frac{\partial}{\partial y_{\lambda}}(f_{\nu,\alpha_{\nu}})=0$ for all admissible $\lambda,\nu$. Therefore only terms $\Theta$ in $\hat{X}$ of the form $\theta(z,y,w)\frac{\partial}{\partial y_{\lambda}}, \, \theta(\bar{z},\bar{y},\bar{w})\frac{\partial}{\partial \overline{y}_{\lambda}}$ can depend on the variable $y$, where $\theta$ is a monomial.

	3) Monomials of the form $\sigma=d \cdot \bar{z}^{e}z^{a}y^{b}u^{c}$ with $b\neq 0, \, |e|=2$, can occur in condition \eqref{eq13} only if $\hat{X}$ contains a term $\Theta= \tilde{d} \cdot z^{a}y^{b}w^{c}\frac{\partial}{\partial y_{\lambda}}$ for some $\lambda$. Then $\sigma$ is one of the monomials in the expression $\Theta(\bar{p}_{\lambda}y_{\lambda})$ (which enters \eqref{eq13}, since $\bar{p}_{\lambda}y_{\lambda}$ enters \eqref{eq28}). Then \eqref{eq13} also contains a monomial $\sigma_{1}=d_{1}\bar{z}^{e_{1}}z^{a}y^{b}u^{c}$ such that $z^{e_{1}}=d_{2}\cdot \tilde{p}_{\lambda}$ for some $d_{2}$ (recall that the monomials $\tilde{p}_{\lambda}$ are defined in Lemma 5).

	Hence, separating in condition \eqref{eq13} the monomials of the form $d_{3}\cdot\overline{\tilde{p}_{\lambda}}\cdot z^{a}y^{b}u^{c}$ with $b\neq 0$, we obtain $d_{3}=0$, whence $\frac{\partial}{\partial y_{\lambda}}(g_{\mu,\beta_{\mu}})=0$ for all admissible $\lambda,\mu$.

	Lemma 10 is proved.

	\vspace{3ex}

	\textbf{Lemma 11.} If $X\in {\rm st} \, M_{0}$ is a field of the first type and $\hat{X}$ depends on $w$, then $X_{1}\in {\rm st} \, M_{0}$.
	
	\textbf{Proof.} By contradiction. Suppose that $X_{1}$ does not belong to ${\rm st} \, M_{0}$. Then the following holds:
	
	a) either the coefficient of the supporting component $\hat{X}_{1}$ of the field $X_{1}$ at $\frac{\partial}{\partial z_{j}}$ is nonzero for some $j$,
	
	b) or the coefficient of the supporting component $\hat{X}_{1}$ of the field $X_{1}$ at $\frac{\partial}{\partial y_{j}}$ is nonzero for some $j$.
	
	We have
	$X_{1}=2 \, {\rm Re} \, \Big(\frac{\partial}{\partial w}(f)\cdot\frac{\partial}{\partial z}+\frac{\partial}{\partial w}(g)\cdot\frac{\partial}{\partial y}+\frac{\partial}{\partial w}(h)\cdot\frac{\partial}{\partial w}\Big)$. 
	
	In case a), we obtain that the lowest component of the tangency condition for $X_{1}$ contains a term of the form $c \, \bar{z}_{j}$ for some constant $c\neq 0$. Therefore the lowest component of the tangency condition must also contain the term $-c \, \bar{z}_{j}$. It follows that the supporting component $\hat{X}_{1}$ of the field $X_{1}$ contains the term $2 i c \bar{z}_{j} \frac{\partial}{\partial \bar{w}}$, i.e. the coefficient at $\frac{\partial}{\partial w}$ of the field $\hat{X}_{1}$ is nonzero. This is impossible, since the field $X_{1}$ is a field of the first type (because the field $X$ is a field of the first type).
	
	In case b), we obtain that the lowest component of the tangency condition for $X_{1}$ contains a term of the form $c \, \overline{\tilde{p}}_{j}$ for some constant $c\neq 0$. Therefore the lowest component of the tangency condition must also contain the term $-c \, \overline{\tilde{p}}_{j}$. It follows that the supporting component $\hat{X}_{1}$ of the field $X_{1}$ contains the term $2 i c \overline{\tilde{p}}_{j} \frac{\partial}{\partial \bar{w}}$, i.e. the coefficient at $\frac{\partial}{\partial w}$ of the field $\hat{X}_{1}$ is nonzero. This is impossible, since the field $X_{1}$ is a field of the first type (because the field $X$ is a field of the first type).
	
	The contradictions obtained mean that $X_{1}\in {\rm st} \, M_{0}$.
	
	Lemma 11 is proved.

	\vspace{3ex}
	We shall need one more lemma. Introduce the necessary notation.
	Let $X=\Xi_{1}+\Xi_{2}+\Xi_{3} \in {\rm aut} \, M_{0}$, where $\Xi_{1}=2 \, {\rm Re} \, \Big(f\frac{\partial}{\partial z}\Big)=\sum_{j=q-2}^{\infty}\Xi_{1}^{(j)}, \, \Xi_{2}=2 \, {\rm Re} \, \Big(g\frac{\partial}{\partial y}\Big)=\sum_{j=q-3}^{\infty}\Xi_{2}^{(j)}, \, \Xi_{3}=2 \, {\rm Re} \, \Big(h\frac{\partial}{\partial w}\Big)=\sum_{j=q-2}^{\infty}\Xi_{3}^{(j)},$ and where the $\Xi_{s}^{(j)}$, in turn, are homogeneous components of weight $j$. The field $X$ has order $q$, and its supporting component $\hat{X}$ is equal to $\Xi_{1}^{(q-2)}+\Xi_{2}^{(q-3)}+\Xi_{3}^{(q-2)}$.
	
	Let also $Y=\Psi_{1}+\Psi_{2}+\Psi_{3}$ be another field from ${\rm aut} \, M_{0}$, where $\Psi_{1}=2 \, {\rm Re} \, \Big(f'\frac{\partial}{\partial z}\Big)=\sum_{j=q'-2}^{\infty}\Psi_{1}^{(j)}, \, \Psi_{2}=2 \, {\rm Re} \, \Big(g'\frac{\partial}{\partial y}\Big)=\sum_{j=q'-3}^{\infty}\Psi_{2}^{(j)}, \, \Psi_{3}=2 \, {\rm Re} \, \Big(h'\frac{\partial}{\partial w}\Big)=\sum_{j=q'-2}^{\infty}\Psi_{3}^{(j)},$ and where the $\Psi_{s}^{(j)}$, in turn, are homogeneous components of weight $j$. The field $Y$ has order $q'$, and its supporting component $\hat{Y}$ is equal to $\Psi_{1}^{(q'-2)}+\Psi_{2}^{(q'-3)}+\Psi_{3}^{(q'-2)}$.

	Let us find the supporting component of the field $[X,Y]=\Phi_{1}+\Phi_{2}+\Phi_{3}$, where $\Phi_{1}=2 \, {\rm Re} \, \Big(\tilde{f}\frac{\partial}{\partial z}\Big)=\sum_{j=q+q'-5}^{\infty}\Phi_{1}^{(j)}, \, \Phi_{2}=2 \, {\rm Re} \, \Big(\tilde{g}\frac{\partial}{\partial y}\Big)=\sum_{j=q+q'-6}^{\infty}\Phi_{2}^{(j)}, \, \Phi_{3}=2 \, {\rm Re} \, \Big(\tilde{h}\frac{\partial}{\partial w}\Big)=\sum_{j=q+q'-5}^{\infty}\Phi_{3}^{(j)}$.
	
	Let
	
	\begin{equation}\label{eq52}
	\Phi_{1}=\Xi_{1}^{(q-2)}\Psi_{1}^{(q'-2)}+\Xi_{3}^{(q-2)}\Psi_{1}^{(q'-2)}+\Xi_{2}^{(q-3)}\Psi_{1}^{(q'-1)}-$$
	$$-\Psi_{1}^{(q'-2)}\Xi_{1}^{(q-2)}-\Psi_{3}^{(q'-2)}\Xi_{1}^{(q-2)}-\Psi_{2}^{(q'-3)}\Xi_{1}^{(q-1)},
	\end{equation}
	
	$$\tilde{\Phi}_{1}=\Xi_{1}^{(q-2)}\Psi_{1}^{(q'-2)}+\Xi_{3}^{(q-2)}\Psi_{1}^{(q'-2)}-\Psi_{1}^{(q'-2)}\Xi_{1}^{(q-2)}-\Psi_{3}^{(q'-2)}\Xi_{1}^{(q-2)},$$
	
	$$\Phi_{2}=\Xi_{1}^{(q-2)}\Psi_{2}^{(q'-3)}+\Xi_{3}^{(q-2)}\Psi_{2}^{(q'-3)}+\Xi_{2}^{(q-3)}\Psi_{2}^{(q'-2)}-$$	
	$$-\Psi_{1}^{(q'-2)}\Xi_{2}^{(q-3)}-\Psi_{3}^{(q'-2)}\Xi_{2}^{(q-3)}-\Psi_{2}^{(q'-3)}\Xi_{2}^{(q-2)},$$
	
	$$\Phi_{3}=\Xi_{1}^{(q-2)}\Psi_{3}^{(q'-2)}+\Xi_{3}^{(q-2)}\Psi_{3}^{(q'-2)}-\Psi_{1}^{(q'-2)}\Xi_{3}^{(q-2)}-\Psi_{3}^{(q'-2)}\Xi_{3}^{(q-2)}.$$
	
	\textbf{Lemma 12.} 1) If at least one of the expressions $\Phi_{1},\Phi_{2},\Phi_{3}$ is nonzero, then $\widehat{[X,Y]}=\Phi_{1}+\Phi_{2}+\Phi_{3}$.
	
	2) If the fields $\hat{X},\hat{Y}$ do not contain terms of the form $d \, \bar{w}^{c}\frac{\partial}{\partial \bar{z}_{s}}$ and $d \, \bar{w}^{c}\frac{\partial}{\partial \bar{y}_{s}}$, then $\Phi_{1}$ takes the simpler form $\tilde{\Phi}_{1}$.
	
	3) In particular, assertion 2) holds for fields $X,Y$ of the first type.
	
	\textbf{Proof.}	
	We have:
	
	$$\Phi_{1}^{(q+q'-5)}=\Xi_{2}^{(q-3)}\Psi_{1}^{(q'-2)}-\Psi_{2}^{(q'-3)}\Xi_{1}^{(q-2)},$$
	
	$$\Phi_{2}^{(q+q'-6)}=\Xi_{2}^{(q-3)}\Psi_{2}^{(q'-3)}-\Psi_{2}^{(q'-3)}\Xi_{2}^{(q-3)},$$
	
	$$\Phi_{3}^{(q+q'-5)}=\Xi_{2}^{(q-3)}\Psi_{3}^{(q'-2)}-\Psi_{2}^{(q'-3)}\Xi_{3}^{(q-2)}.$$
	
	By Lemma 10, the coefficients of the supporting components do not depend on $y,\bar{y}$, therefore $\Phi_{1}^{(q+q'-5)}=\Phi_{2}^{(q+q'-6)}=\Phi_{3}^{(q+q'-5)}=0$. Consider the next components by weight:
	
	$$\Phi_{1}^{(q+q'-4)}=\Xi_{1}^{(q-2)}\Psi_{1}^{(q'-2)}+\Xi_{3}^{(q-2)}\Psi_{1}^{(q'-2)}+\Xi_{2}^{(q-2)}\Psi_{1}^{(q'-2)}+\Xi_{2}^{(q-3)}\Psi_{1}^{(q'-1)}-$$
	$$-\Psi_{1}^{(q'-2)}\Xi_{1}^{(q-2)}-\Psi_{3}^{(q'-2)}\Xi_{1}^{(q-2)}-\Psi_{2}^{(q'-2)}\Xi_{1}^{(q-2)}-\Psi_{2}^{(q'-3)}\Xi_{1}^{(q-1)}=$$
	
	$$=\Xi_{1}^{(q-2)}\Psi_{1}^{(q'-2)}+\Xi_{3}^{(q-2)}\Psi_{1}^{(q'-2)}+\Xi_{2}^{(q-3)}\Psi_{1}^{(q'-1)}-$$
	$$-\Psi_{1}^{(q'-2)}\Xi_{1}^{(q-2)}-\Psi_{3}^{(q'-2)}\Xi_{1}^{(q-2)}-\Psi_{2}^{(q'-3)}\Xi_{1}^{(q-1)},$$
	
	$$\Phi_{2}^{(q+q'-5)}=\Xi_{1}^{(q-2)}\Psi_{2}^{(q'-3)}+\Xi_{3}^{(q-2)}\Psi_{2}^{(q'-3)}+\Xi_{2}^{(q-2)}\Psi_{2}^{(q'-3)}+\Xi_{2}^{(q-3)}\Psi_{2}^{(q'-2)}-$$	
	$$-\Psi_{1}^{(q'-2)}\Xi_{2}^{(q-3)}-\Psi_{3}^{(q'-2)}\Xi_{2}^{(q-3)}-\Psi_{2}^{(q'-2)}\Xi_{2}^{(q-3)}-\Psi_{2}^{(q'-3)}\Xi_{2}^{(q-2)}=$$
	
	$$=\Xi_{1}^{(q-2)}\Psi_{2}^{(q'-3)}+\Xi_{3}^{(q-2)}\Psi_{2}^{(q'-3)}+\Xi_{2}^{(q-3)}\Psi_{2}^{(q'-2)}-$$	
	$$-\Psi_{1}^{(q'-2)}\Xi_{2}^{(q-3)}-\Psi_{3}^{(q'-2)}\Xi_{2}^{(q-3)}-\Psi_{2}^{(q'-3)}\Xi_{2}^{(q-2)},$$
	
	$$\Phi_{3}^{(q+q'-4)}=\Xi_{1}^{(q-2)}\Psi_{3}^{(q'-2)}+\Xi_{3}^{(q-2)}\Psi_{3}^{(q'-2)}+\Xi_{2}^{(q-2)}\Psi_{3}^{(q'-2)}+\Xi_{2}^{(q-3)}\Psi_{3}^{(q'-1)}-$$
	$$-\Psi_{1}^{(q'-2)}\Xi_{3}^{(q-2)}-\Psi_{3}^{(q'-2)}\Xi_{3}^{(q-2)}-\Psi_{2}^{(q'-2)}\Xi_{3}^{(q-2)}-\Psi_{2}^{(q'-3)}\Xi_{3}^{(q-1)}=$$
	
	$$=\Xi_{1}^{(q-2)}\Psi_{3}^{(q'-2)}+\Xi_{3}^{(q-2)}\Psi_{3}^{(q'-2)}+\Xi_{2}^{(q-3)}\Psi_{3}^{(q'-1)}-$$
	$$-\Psi_{1}^{(q'-2)}\Xi_{3}^{(q-2)}-\Psi_{3}^{(q'-2)}\Xi_{3}^{(q-2)}-\Psi_{2}^{(q'-3)}\Xi_{3}^{(q-1)},$$
	where we have again used Lemma 10. If at least one of the fields $\Phi_{1}^{(q+q'-4)}, \, \Phi_{2}^{(q+q'-5)}, \, \Phi_{3}^{(q+q'-4)}$ is nonzero, then the supporting component $\widehat{[X,X_{1}]}$ is equal to $\Phi_{1}^{(q+q'-4)} + \Phi_{2}^{(q+q'-5)} + \Phi_{3}^{(q+q'-4)}$. These expressions can be simplified as follows.	
	First, the coefficient at $\frac{\partial}{\partial w}$ of the field $X$ does not depend on $y$. Indeed, if the field $X$ contains a monomial of the form $\Theta=d \, z^{a}y^{b}w^{c}\frac{\partial}{\partial w}$ with $b\neq 0$, then the monomial $\frac{d}{2i} \, z^{a}y^{b}u^{c}$ appears in the tangency condition (since the monomial $\chi=\frac{w}{2i}$ enters the defining relation). Other pairs of monomials $(\Theta,\chi)$ by means of which one could obtain the monomial $\frac{d}{2i} \, z^{a}y^{b}u^{c}$ are $(d_{1} \, \bar{w}^{c}\frac{\partial}{\partial \bar{z}_{t}},d_{2} \, z^{a}y^{b}\bar{z}_{t})$ and $(d_{1} \, \bar{w}^{c}\frac{\partial}{\partial \bar{y}_{t}},d_{2} \, z^{a}y^{b}\bar{y}_{t})$. But for $b\neq 0$, monomials of the form $d_{2} \, z^{a}y^{b}\bar{z}_{t}$ and $d_{2} \, z^{a}y^{b}\bar{y}_{t}$ are absent from the defining relation by condition 3 for the normal form. Hence $d=0$, as required. The same is true for the field $Y$. It follows that $\Xi_{2}^{(s)}\Psi_{3}^{(t)}=\Psi_{2}^{(s)}\Xi_{3}^{(t)}=0$ for all admissible values of $s$ and $t$. Therefore we obtain
	
	$$\Phi_{3}^{(q+q'-4)}=\Xi_{1}^{(q-2)}\Psi_{3}^{(q'-2)}+\Xi_{3}^{(q-2)}\Psi_{3}^{(q'-2)}-\Psi_{1}^{(q'-2)}\Xi_{3}^{(q-2)}-\Psi_{3}^{(q'-2)}\Xi_{3}^{(q-2)}.$$	
	Second, if the vector coefficient of the field $\Xi_{1}^{(q-1)}$ contains monomials linearly depending on $y$, i.e. monomials of the form $d \, z^{a}y_{t}w^{c}\frac{\partial}{\partial z_{s}}$, then monomials of the form $\sigma_{1}=d \, z^{a}y_{t}u^{c}\bar{z}_{s}$ appear in the tangency condition. Then the monomial $-\sigma_{1}$ must also appear in the tangency condition. By condition 3 for the normal form it can be obtained only by means of the following pairs $(\Theta,\chi)$, where $\Theta$ is a term of the field $X$ and $\chi$ is a monomial in the defining relation: $(\Theta,\chi)=(d_{1} \, \bar{w}^{c}\frac{\partial}{\partial \bar{z}_{s_{1}}},d_{2}y_{t}\bar{z}_{s_{1}}\bar{z}_{s})$ for $a=0$, and $(\Theta,\chi)=(d_{1} \, \bar{w}^{c}\frac{\partial}{\partial \bar{y}_{s_{1}}},d_{2}z^{a}\bar{z}_{s}y_{t}\bar{y}_{s_{1}})$ for $|a|=1$. Since a monomial of the form $d \, z^{a}y_{t}w^{c}\frac{\partial}{\partial z_{s}}$ has weight $q-1$, monomials of the form $d_{1} \, \bar{w}^{c}\frac{\partial}{\partial \bar{z}_{s_{1}}}$ and $d_{1} \, \bar{w}^{c}\frac{\partial}{\partial \bar{y}_{s_{1}}}$ have weights $q-2$ and $q-3$, respectively, and hence lie in the supporting component $\hat{X}$. The same is true for the field $Y$. And if the fields $\hat{X},\hat{Y}$ do not contain terms $d_{1} \, \bar{w}^{c}\frac{\partial}{\partial \bar{z}_{s_{1}}}$ and $d_{1} \, \bar{w}^{c}\frac{\partial}{\partial \bar{y}_{s_{1}}}$, then the terms of the field $\Phi_{1}^{(q+q'-4)}$ whose coefficients do not depend on $y$ take the simpler form:
	
	$$\Xi_{1}^{(q-2)}\Psi_{1}^{(q'-2)}+\Xi_{3}^{(q-2)}\Psi_{1}^{(q'-2)}-\Psi_{1}^{(q'-2)}\Xi_{1}^{(q-2)}-\Psi_{3}^{(q'-2)}\Xi_{1}^{(q-2)}.$$
	
	And since, by Lemma 10, the coefficients of the supporting component do not depend on $y$, this expression is exactly equal to the field $\Phi_{1}^{(q+q'-4)}$.
	
	For example, such a simpler form of the field $\Phi_{1}^{(q+q'-4)}$ will hold for fields of the first type, since if terms $d_{1} \, \bar{w}^{c}\frac{\partial}{\partial \bar{z}_{s_{1}}}$ and $d_{1} \, \bar{w}^{c}\frac{\partial}{\partial \bar{y}_{s_{1}}}$ are present in the supporting component $\hat{X}$, we also obtain, respectively, terms of the form $d_{2} \, z_{s_{1}}w^{c}\frac{\partial}{\partial w}$ and $d_{2} \, \tilde{p}_{s_{1}}w^{c}\frac{\partial}{\partial w}$ in the supporting component, and hence $X$ is a field of the second type.
	
	Lemma 12 is proved.
	
	\vspace{3ex}
	
	Moreover, below we shall use Lemma 12 not only to compute the field $\widehat{[X,Y]}$, but also to find its component of type $(1,0)$.

	Now we estimate the dimension ${\rm dim \, st} \, M_{0}$. We note that below, since $l=2$, the inequality $m<n$ holds.

	\vspace{3ex}

	\textbf{Lemma 13.}
	${\rm dim \, st} \, M_{0}\leq n^{2}+2n+1$. 
	
	\textbf{Proof.}
	A brief plan of the proof:

		 1) We consider fields $X$ of the first type such that the coefficients of the field $\hat{X}$ do not depend on $w$ (fields of type 1(0)). We prove that the number of linearly independent fields of this type does not exceed $n^{2}-2k$.
		 
		 2) We consider fields $X$ of the first type such that the coefficients of the field $\hat{X}$ depend linearly on $w$. We prove that among such fields only fields of type 1.1(1) can be nonzero (item 2.c.1 of the proof), and fields of the other types do not exist (items 2.a, 2.b, 2.c.2, and 2.c.3 of the proof).
		 
		 3) We consider fields $X$ of the first type such that the coefficients of the field $\hat{X}$ depend on $w^{\nu}, \, \nu\geq 2$. We prove that among such fields only fields of type 1.1(2) can be nonzero (item 3.a of the proof).
		 
		 3.b) We prove that fields of types 1.1(1) and 1.1(2) can be simultaneously diagonalized and that the dimension of the space of fields of each of these two types does not exceed $m$ (see item 3.b).

		 4.a) Fields of type $2.0(\geq 3)$ do not exist.
		 
		 4.b) Fields of type $2.\geq 3(\geq 0)$ do not exist.
		 
		 4.c) Fields of type $2.2(\geq 2)$ do not exist.
		 
		 4.d) Fields of type $2.1(\geq 3)$ do not exist.
		 
		 4.e) We estimate the number of linearly independent fields of the remaining types ($2.b(0), \, 2.b(1), \, 2.0(2), \, 2.1(2)$).
		 
		 5) We estimate the total number of fields of types 1(0), 1.1(1), 1.1(2), 2.1(1), 2.1(2), and 2.2(1). We obtain the upper estimate $n^{2}+2n-1$. Taking item 4 into account, we obtain the upper estimate $n^{2}+2n+1$ for the dimension of the stabilizer of the automorphism algebra.
	
	\vspace{3ex}

	 Consider three cases for fields $X$ of the first type:
	
	1) the coefficients of the field $\hat{X}$ do not depend on $w$,
	
	2) the coefficients of the field $\hat{X}$ depend on $w$, and their degree in the variable $w$ does not exceed one,
	
	3) some of the coefficients of the field $\hat{X}$ depend on $w^{\nu}, \, \nu\geq 2$, and the degree of all coefficients in the variable $w$ does not exceed $\nu$.
	
	The proof consists of considering a tree of cases of depth five. The subcases are numbered by means of a sequence of letters and digits.

	Next we estimate the dimension of the space of fields of types 1.1, 1.2; the dimension of the quotient space of fields of type 1.3 by the space of fields of types 1.1 and 1.2; the dimension of the space of fields of type 2.1; the dimension of the quotient space of fields of type 2.2 by the space of fields of types 1.1 and 2.1; the dimension of the quotient space of fields of type 2.3 by the space of fields of types 1.2 and 2.1; and the dimension of the quotient space of fields of type 2.4 by the space of fields of all other types. To estimate the dimensions of quotient spaces we shall consider a special basis of the space ${\rm aut} \, M_{0}$, which we shall choose in the course of the proof. Further in the proof of this lemma, $a,a_{j},b,c,e$ are multi-degrees with nonnegative integer components, and $d,d_{j},d_{a}\in \mathbb{C}$. These notations are valid only within each numbered subcase; that is, the same parameters may take different values when different subcases are considered.

	\vspace{3ex}
	
	\textit{Case 1.}
	
	\vspace{3ex}
	
	1.a) Let $X$ be a field of type 1.1(0). By Lemma 10 we obtain that $f_{\nu,\alpha_{\nu}}$ is a sum of monomials of the form $d_{a}\cdot z^{a}$ of the same degree (here $a$ is a multi-degree). Separating in condition \eqref{eq13} the monomials of the form $d \cdot \bar{z}_{j}z^{a}$, we obtain $d_{a}=0$ for $|a|>1$. And for $|a|=1$ we obtain that $f_{\nu,\alpha_{\nu}}$ may be equal to a sum of monomials of the form $d_{s}\cdot z_{s}, \, s\neq \nu,$ with complex coefficients $d_{s}$ and monomials of the form $d_{\nu}\cdot z_{\nu}$ in which $d_{\nu} \in i\mathbb{R}$. Altogether there are no more than $m^{2}$ such linearly independent fields.
	
	\vspace{3ex}
	
	1.b) Next, let $X$ be a field of type 1.2(0). By Lemma 10 we obtain that $g_{\mu,\beta_{\mu}}$ is a sum of monomials of the form $d_{a}\cdot z^{a}$ of the same degree. Separating in condition \eqref{eq13} the monomials of the form $d \cdot \bar{z}^{b}z^{a}$ with $|b|=2$, we obtain that $g_{\mu,\beta_{\mu}}$ may be equal to a sum of polynomials of the form $d_{s}\cdot p_{s}=d_{s}\cdot \tilde{p}_{s}+..., \, s\neq \mu,$ with complex coefficients $d_{s}$ and a polynomial of the form $d_{\mu}\cdot p_{\mu}=d_{\mu}\cdot \tilde{p}_{\mu}+...$, in which $d_{\mu} \in i\mathbb{R}$ (here dots denote the sum of monomials different from $\tilde{p}_{s}$ and $\tilde{p}_{\mu}$, respectively). Altogether there are no more than $k^{2}$ such linearly independent fields.
	
	\vspace{3ex}
	
	1.c) Now let $X$ be a field of type 1.3(0). By Lemma 10 we obtain that $f_{\nu,\alpha_{\nu}}$ is a sum of monomials of the form $d_{a}\cdot z^{a}$ of degree $|a|$, and $g_{\mu,\beta_{\mu}}$ is a sum of monomials of the form $d_{b}\cdot z^{b}$ of degree $|b|$. At the same time both the case $|a|=1$ and the case $|b|=2$ have already been considered above (in items 1.a and 1.b). Thus it remains to consider the case $|a|=2$ and $|b|=1$ (for other values of $|a|$ and $|b|$ we have $f_{\nu,\alpha_{\nu}}=g_{\mu,\beta_{\mu}}=0$). Therefore, separating in condition \eqref{eq13} the monomials of the form $d \cdot \bar{z}^{b}z^{a}$, we obtain that $f_{\nu,\alpha_{\nu}}$ may be equal to a sum of polynomials of the form $d_{1}p_{s_{1}}=d_{1}\cdot \tilde{p}_{s_{1}}+...$ with complex coefficients $d_{1}$, and $g_{\mu,\beta_{\mu}}$ may be equal to a sum of polynomials of the form $d_{2}\cdot z_{s_{2}}$ with complex coefficients $d_{2}$. This means that there are no more than $2 m k$ fields of type 1.3(0) (since there are $k$ monomials $\tilde{p}_{s}$, there are $m$ operators $\frac{\partial}{\partial z_{s}}$, and the coefficients $d_{2}$ are complex and are uniquely recovered from $d_{1}$). We note that in the case under consideration $q=3$.
	
	\vspace{3ex}
	
	Thus altogether in cases 1.a, 1.b, and 1.c there are no more than $m^{2}+k^{2}+2 m k=n^{2}$ fields, i.e. no more than fields of the first type for a nondegenerate quadric. However, as we shall now show, the upper estimate in case 1.c is not attained. For this it is necessary to prove that the maximum of the number of linearly independent tangent fields such that the degree of the nonzero polynomials $f_{\nu,\alpha_{\nu}}$ is two is not attained. This maximum is equal to $2 m k$.

	In case 1.c, the left-hand side of equality \eqref{eq36} also contains terms $f_{\nu,\alpha_{\nu}} \cdot \frac{\partial}{\partial z_{\nu}}(p_{j}) \cdot \overline{y}_{j}$. Denote by $F$ the sum of all such terms for all admissible values of $\nu, j$. Two cases are possible: either $F=0$ (case 1.c.1), or $F\neq 0$ (case 1.c.2).
	
	\vspace{3ex}
	
	In case 1.c.1, for the basis fields $f_{\nu,\alpha_{\nu}} \cdot \frac{\partial}{\partial z_{\nu}}$ there are two subcases, which are described by means of the defect $r$.
	
	\vspace{3ex}
	
	1.c.1.1) $1\leq \nu \leq r$. In this case the field $\hat{X}$ may contain one or several nonzero coefficients $f_{\nu,\alpha_{\nu}}$ (there are no more than $2 r k$ such fields, since there are $k$ monomials $\tilde{p}_{s}$, there are $r$ operators $\frac{\partial}{\partial z_{s}}$ for $1\leq s \leq r$, and the coefficients of the fields are complex).
	
	\vspace{3ex}
	
	1.c.1.2) $\nu>r$. In this case for each $\nu$ (there are $m-r$ such $\nu$) there exists at least one polynomial $p_{j}$ such that $f_{\nu,\alpha_{\nu}}\frac{\partial}{\partial z_{\nu}}(p_{j})\neq 0$.
	In particular, this is true for $\nu=m$. This means that if all coefficients $f_{\nu,\alpha_{\nu}}, \, r<\nu<m,$ are given, then the coefficient $f_{m,\alpha_{m}}$ is uniquely recovered (it can be expressed from the equality $\sum_{\nu=r+1}^{m}f_{\nu,\alpha_{\nu}}\frac{\partial}{\partial z_{\nu}}(p_{j})=0$).

	Thus altogether in case 1.c.1.2 there are no more than $2(m-r) k-2k$ linearly independent fields. Indeed, first, there are $k$ monomials $\tilde{p}_{s}$, there are $(m-r)$ operators $\frac{\partial}{\partial z_{s}}$ for $r< s \leq m$, and the coefficients of the fields are complex; altogether there are no more than $2 (m-r) k$ fields. Second, the coefficient $f_{m,\alpha_{m}}$ can take $2k$ linearly independent values, and since it is uniquely recovered from the remaining coefficients $f_{\nu,\alpha_{\nu}}, \, r<\nu<m,$ altogether there are no more than $2 (m-r) k-2k$ fields.

	\vspace{3ex}

	Thus altogether in case 1.c.1 we have no more than $2 r k+2(m-r)k-2k=2mk-2k<2mk$ fields. Note that this estimate does not depend on r. Let us show that case 1.c.2 is impossible.

	\vspace{3ex}
	
	In case 1.c.2 we have that the left-hand side of equality \eqref{eq36} must also contain the terms $-f_{\nu,\alpha_{\nu}} \cdot \frac{\partial}{\partial z_{\nu}}(p_{j}) \cdot \overline{y}_{j}$ (for the same values of $\nu, j$, and before collecting all similar terms in \eqref{eq36}), which do not enter the expression $F$ and contain terms of the form $d \cdot z^{A} \cdot \overline{y}_{j}, \, |A|=3$ (since the degrees of the polynomials $f_{\nu,\alpha_{\nu}}$ and $p_{j}$ are two). Among these terms, after collecting like terms, there are nonzero terms of the form $d_{0} z^{A_{0}}  \overline{y}_{j}, \, |A_{0}|=3$. Fix any of the values $d_{0},A_{0},j$ and put $\sigma_{1}=d_{0} z^{A_{0}}  \overline{y}_{j}$.
	
	At the same time $\sigma_{1}$ may occur as one of the terms in the expression $\Theta(\chi)$, where $\Theta$ is some term of the field $X$, and $\chi$ is some term in the defining relation. The terms $\Theta$ and $\chi$ must satisfy one of several conditions, which we list below.

    Before listing these conditions, let us make three remarks. We number the remarks in the same way as the logical subcases.
    
    \vspace{3ex}
    
    \textit{Remarks.}
    
    1.c.2.1) First, since $\sigma_{1}$ is a monomial, it is enough to consider only such $\chi$ and $\Theta$ for which $\chi$ is a monomial and $\Theta$ is equal to one of the following expressions:  $\theta(z,y,w)\frac{\partial}{\partial z_{\lambda}}, \, \theta(\bar{z},\bar{y},\bar{w})\frac{\partial}{\partial \bar{z}_{\lambda}}, \,\theta(z,y,w)\frac{\partial}{\partial y_{\lambda}}, \,\theta(\bar{z},\bar{y},\bar{w})\frac{\partial}{\partial \overline{y}_{\lambda}}, \, \theta(z,y,w)\frac{\partial}{\partial w}$ or $\theta(\bar{z},\bar{y},\bar{w})\frac{\partial}{\partial \bar{w}}$, where $\theta$ is a monomial.
    
    \vspace{3ex}
    
    1.c.2.2) Second, $\Theta$ cannot be equal to $\theta(z,y,w)\frac{\partial}{\partial w}$ or $\theta(\bar{z},\bar{y},\bar{w})\frac{\partial}{\partial \bar{w}}$. Indeed, suppose the contrary. By the rigidity condition for the defining equation (i.e. the defining equation does not depend on the variable $u$), we have $\chi=\frac{1}{2i}w$ or, respectively, $\chi=-\frac{1}{2i}\bar{w}$.

    Let $\chi=\frac{1}{2i}w, \, \Theta=\theta(z,y,w)\frac{\partial}{\partial w}$.
    Then, as was already shown in the proof of Lemma 12, $\theta$ does not depend on $y$, since otherwise monomials of the form $d\cdot z^{a}y^{b}u^{c}$ with $b\neq 0$ appear in \eqref{eq36} (here $a,b$ are multi-degrees, $c$ is a nonnegative integer, and $d\in\mathbb{C}$), and monomials of this form cannot be obtained in any other way except by means of the indicated $\Theta$ and $\chi$. Therefore $\Theta=d\cdot z^{a}w^{c}\frac{\partial}{\partial w}$.
    
    Next, the term $\sigma_{1}$ is one of the terms obtained if one replaces the variable $w=u+iv$ in $\Theta$ by one of the monomials entering the expression $iv$, where $v$ is given by formula \eqref{eq28}. It is clear that $c=1$, since the monomial $\sigma_{1}$ has degree one in the conjugate variable.
    
    Therefore the term $\sigma_{1}$ is one of the terms obtained if the variable $w=u+iv$ is replaced by the polynomial $ip_{j}\overline{y}_{j}$ entering the expression $iv$. Moreover, $iv$ cannot be replaced by any other expression entering it, since this expression must depend linearly on the conjugate variable (and it must be the variable $\overline{y}_{j}$), and terms with this property can be contained only in the expression $ip_{j}\overline{y}_{j}$ (by property 3 of the normal form). Hence we obtain $|a|=1$ (since the degree of the polynomial $p_{j}$ is two).
    
    But then the order of the field $\Theta$ must be equal to $q=3$ (since the variable $w$ in $\Theta$ can be replaced by $i\langle z,\bar{z}\rangle$), i.e. $\Theta$ enters the supporting component $\hat{X}$ of the field $X$. This is impossible, since we are considering a field $X$ of the first type.

    In the case $\chi=-\frac{1}{2i}\bar{w}$ we have $\Theta=d\cdot\bar{z}^{a}\bar{w}\frac{\partial}{\partial \bar{w}}$ (for analogous reasons). Next, it is clear that $a=0$ (since $\sigma_{1}$ does not depend on $\bar{z}$). Thus the field $\Theta$ has order two, which is less than $q=3$. Hence this case is also impossible.
    
    \vspace{3ex}

    1.c.2.3) Third, note that the coefficient of the field $\Theta$, where $\Theta$ is equal to one of the expressions $\theta(z,y,w)\frac{\partial}{\partial z_{\lambda}}, \,\theta(z,y,w)\frac{\partial}{\partial y_{\lambda}}$, cannot depend on the variable $w$, since otherwise, to obtain the term $\sigma_{1}$, one must replace the variable $w=u+iv$ by the polynomial $i p_{j}\bar{y}_{j}$ entering the expression $iv$. But then the order of the field $\Theta$ must be equal to $q=3$ (since the variable $w$ in $\Theta$ can be replaced by $i\langle z,\bar{z}\rangle$), i.e. $\Theta$ enters the supporting component $\hat{X}$ of the field $X$. This is impossible, since we are considering a field $X$ of type 1(0).
    
    Also the coefficient of the field $\Theta$, where $\Theta$ is equal to one of the expressions $\theta(\bar{z},\bar{y},\bar{w})\frac{\partial}{\partial \bar{z}_{\lambda}}, \,\theta(\bar{z},\bar{y},\bar{w})\frac{\partial}{\partial \overline{y}_{\lambda}}$, cannot depend on the variable $\bar{w}$, since otherwise, to obtain the term $\sigma_{1}$, one must replace the variable $\bar{w}=u-iv$ by the polynomial $-i p_{j}\bar{y}_{j}$ entering the expression $-iv$. But then the order of the field $\Theta$ must be equal to $q=3$ (since the variable $\bar{w}$ in $\Theta$ can be replaced by $-i\langle z,\bar{z}\rangle$), i.e. $\Theta$ enters the supporting component $\hat{X}$ of the field $X$. This is impossible, since we are considering a field $X$ of type 1(0).

    \vspace{3ex}
    
    Now we write out the conditions on $\Theta$ and $\chi$:
    
    \vspace{3ex}

	1.c.2.a) The field $X$ contains a term of the form $\Theta=d_{3}z^{c}\frac{\partial}{\partial y_{\lambda}}$, and the defining function contains a term of the form $\chi=d_{4}z^{e} y_{\lambda} \cdot \overline{y}_{j}$.
	
	This condition cannot be fulfilled, since the defining function does not contain terms of the indicated form by condition 3 for the normal form.
	
	\vspace{3ex}

	1.c.2.b) The field $X$ contains a term of the form $\Theta=d_{3}z^{c}\frac{\partial}{\partial z_{\lambda}}$, and the defining function contains a term of the form $\chi=d_{4}z^{e} \cdot \overline{y}_{j}$.
	
	Then $|e|= 2$ by property 3 of the normal form, whence $\chi$ is equal to one of the terms in the expression $p_{j}\overline{y}_{j}$. But this is impossible, since the expression $p_{j}\overline{y}_{j}$ was already used earlier in constructing the monomial $\sigma_{1}$ (if $\chi$ is equal to one of the monomials in the expression $p_{j}\overline{y}_{j}$, then $\Theta$ enters the supporting component of the field $X$ and was also already used earlier, i.e. $\Theta(\chi)$ enters $F$).
	Therefore this case is impossible.

    \vspace{3ex}

	1.c.2.c) The field $X$ contains a term of the form $\Theta=d_{3} \cdot \overline{y}_{j} \frac{\partial}{\partial \bar{z}_{\lambda}}$, and the defining function contains a term of the form $\chi=d_{4}  \cdot z^{c} \cdot \bar{z}_{\lambda}$ (in the case under consideration $|c|=3$).
	
	The defining function does not contain terms of the indicated form by condition 3 for the normal form.
	
	\vspace{3ex}
	
	1.c.2.d) The field $X$ contains a term of the form $\Theta=d_{3} \cdot \overline{y}_{j} \frac{\partial}{\partial \overline{y}_{\lambda}}$, and the defining function contains a term of the form $\chi=d_{4}  \cdot z^{c} \cdot \overline{y}_{\lambda}$, where $|c|=3>2$.
	
	This case is impossible, because the defining function does not contain terms of the indicated form by condition 3 for the normal form.

	\vspace{3ex}
	
	Thus altogether in case 1 we have no more than $m^{2}$ fields of type 1.1(0), no more than $k^{2}$ fields of type 1.2(0), and no more than $2mk-2k$ fields of type 1.3(0), hence altogether no more than $m^{2}+k^{2}+2mk-2k=n^{2}-2k$ linearly independent fields.

	\vspace{3ex}
	
	\textit{Case 2.}
	
	Consider the field $X_{1}=[X_{0},X]$, which also belongs to the automorphism algebra.
	The field $X_{1}$ has order $\tilde{q}=q-2$.
	To simplify notation, in what follows we omit the tilde over $q$.

	By Lemma 11 we have $X_{1}\in {\rm st} \, M_{0}$; hence the field $X_{1}$ falls into the preceding case (case 1), i.e. its supporting component has the special form described in the preceding case.
	
	\vspace{3ex}

	2.a) Let $X$ be a field of type 1.1.
	Then by Lemma 10 we have $f_{\alpha}=\tilde{f}_{1,\alpha}(z) w+\tilde{f}_{0,\alpha}(z)$. At the same time the field $X_{1}$ has type 1.1(0), i.e. falls into case 1, considered above, and therefore the vector polynomial $\tilde{f}_{1,\alpha}(z)$ has the special form indicated in item 1.a. Namely, $\tilde{f}_{1,\alpha}(z)$ is a vector whose coordinates are linear in $z$ or zero, and hence $\tilde{f}_{0,\alpha}(z)$ is a vector whose coordinates are cubic in $z$ or zero.

	Let $F$ be equal to the expression $\lfloor-2 \, {\rm Re} \, \langle f_{\alpha},\bar{z}\rangle\rfloor$ (the lowest weighted component of the tangency condition), in which we put $w=\bar{w}=0$. That is, $F=-2 \, {\rm Re} \, \langle \tilde{f}_{0,\alpha}(z),\bar{z}\rangle$.
	
	The remaining terms in the lowest weighted component of the tangency condition have the form $H=-2 \, {\rm Re} \, \langle \tilde{f}_{1,\alpha}(z) (u+i\langle z,\bar{z}\rangle),\bar{z}\rangle$.
	
	Thus we have $F+H=\lfloor-2 \, {\rm Re} \, \langle f_{\alpha},\bar{z}\rangle\rfloor=0$.
	Note that the monomials entering the expression $F$ are different from the monomials entering the expression $H$ (these monomials have different bidegrees in the variables $z,\bar{z}$). Therefore the equality $F=H=0$ holds.
	
	Next, we have $H=- \langle \tilde{f}_{1,\alpha}(z) (u+i\langle z,\bar{z}\rangle),\bar{z}\rangle- \langle \overline{\tilde{f}_{1,\alpha}(z)} (u-i\langle z,\bar{z}\rangle),z\rangle=- (\langle \tilde{f}_{1,\alpha}(z),\bar{z}\rangle+ \langle \overline{\tilde{f}_{1,\alpha}(z)},z\rangle)u - (\langle \tilde{f}_{1,\alpha}(z),\bar{z}\rangle- \langle \overline{\tilde{f}_{1,\alpha}(z)} ,z\rangle)i\langle z,\bar{z}\rangle$. Hence we obtain $\langle \tilde{f}_{1,\alpha}(z),\bar{z}\rangle+ \langle \overline{\tilde{f}_{1,\alpha}(z)},z\rangle =  \langle \tilde{f}_{1,\alpha}(z),\bar{z}\rangle- \langle \overline{\tilde{f}_{1,\alpha}(z)} ,z\rangle=0$, and therefore $\langle \overline{\tilde{f}_{1,\alpha}(z)},z\rangle = 0$. This is possible only if $\tilde{f}_{1,\alpha}=0$. That is, $X$ does not depend on $w$ -- a contradiction.

	\vspace{3ex}
	
	2.b) Let $X$ be a field of type 1.2.
	Then by Lemma 10 we have $g_{\beta}=\tilde{g}_{1,\beta}(z) w+\tilde{g}_{0,\beta}(z)$. At the same time the field $X_{1}$ has type 1.2(0), and hence falls into case 1 considered above; therefore the vector polynomial $\tilde{g}_{1,\beta}(z)$ has the special form indicated in item 1.b. Namely, $\tilde{g}_{1,\beta}(z)$ is a vector whose coordinates are quadratic in $z$ or zero, and hence $\tilde{g}_{0,\beta}(z)$ is a vector whose coordinates have degree four in $z$ or are zero.

	Let $F$ be equal to the expression $\lfloor-2 \, {\rm Re} \, (\mathcal{P}(z)\bar{g}_{\beta})\rfloor$ (the lowest weighted component of the tangency condition), in which we put $w=\bar{w}=0$. That is, $F=-2 \, {\rm Re} \, (\mathcal{P}(z)\overline{\tilde{g}_{0,\beta}(z)})$.
	
	The remaining terms in the lowest weighted component of the tangency condition have the form $H=-2 \, {\rm Re} \, \Big((\mathcal{P}(z)\overline{\tilde{g}_{1,\beta}(z)}) (u-i\langle z,\bar{z}\rangle)\Big)$.
	
	Thus we have $F+H=\lfloor-2 \, {\rm Re} \, (\mathcal{P}(z)\bar{g}_{\beta})\rfloor=0$.
	Note that the monomials entering the expression $F$ are different from the monomials entering the expression $H$ (these monomials have different bidegrees in the variables $z,\bar{z}$). Therefore the equality $F=H=0$ holds.
	
	Next, we have $H=- (\mathcal{P}(z)\overline{\tilde{g}_{1,\beta}(z)}) (u-i\langle z,\bar{z}\rangle)- (\overline{\mathcal{P}(z)}\tilde{g}_{1,\beta}(z)) (u+i\langle z,\bar{z}\rangle)=- (\mathcal{P}(z)\overline{\tilde{g}_{1,\beta}(z)}+ \overline{\mathcal{P}(z)}\tilde{g}_{1,\beta}(z))u+ (\mathcal{P}(z)\overline{\tilde{g}_{1,\beta}(z)} -\overline{\mathcal{P}(z)}\tilde{g}_{1,\beta}(z))i\langle z,\bar{z}\rangle$. Hence we obtain $\mathcal{P}(z)\overline{\tilde{g}_{1,\beta}(z)}+ \overline{\mathcal{P}(z)}\tilde{g}_{1,\beta}(z)=  \mathcal{P}(z)\overline{\tilde{g}_{1,\beta}(z)} -\overline{\mathcal{P}(z)}\tilde{g}_{1,\beta}(z)=0$, and therefore $\mathcal{P}(z)\overline{\tilde{g}_{1,\beta}(z)} = 0$. This is possible only if $\tilde{g}_{1,\beta}=0$, since the components of the vector $\mathcal{P}$ are linearly independent. That is, $X$ does not depend on $w$ -- a contradiction.
	
	\vspace{3ex}
	
	2.c) Let $X$ be a field of type 1.3.
	Then by Lemma 10 we have
	$$f_{\alpha}=\tilde{f}_{1,\alpha}(z) w+\tilde{f}_{0,\alpha}(z), \ \ \ g_{\beta}=\tilde{g}_{1,\beta}(z) w+\tilde{g}_{0,\beta}(z).$$ 
	At the same time the field $X_{1}$ falls into case 1, considered above. Therefore the field $X_{1}$ falls into one of the cases 1.a, 1.b, 1.c.
	
	\vspace{3ex}
	
	2.c.1) If $X_{1}$ falls into case 1.a, then $\tilde{g}_{1,\beta}(z)=0$, and the field $X$ has type $1.1(1)$. We shall estimate the number of such fields later.

	\vspace{3ex}

	2.c.2) If $X_{1}$ falls into case 1.b, then $\tilde{f}_{1,\alpha}(z)=0$, and the components of the vector $\tilde{g}_{1,\beta}(z)$ are quadratic in $z$ or zero.
	And therefore the components of the vector $\tilde{f}_{0,\alpha}(z)$ are polynomials of degree five in $z$ or zero, and the components of the vector $\tilde{g}_{0,\beta}(z)$ are polynomials of degree four in $z$ or zero.
	Therefore the lowest component of the tangency condition takes the form

	\begin{equation}\label{eq48}
	-2 \, {\rm Re} \, \Big(\overline{\mathcal{P}(z)}\tilde{g}_{1,\beta}(z)  \, \cdot (u+i \langle z,\bar{z}\rangle)+\langle \tilde{f}_{0,\alpha}(z) ,\bar{z}\rangle+\overline{\mathcal{P}(z)}\tilde{g}_{0,\beta}(z)\Big)=0.
	\end{equation}
	Separating in \eqref{eq48} the terms of bidegrees $(2,2)$ and $(3,3)$ in the variables $(z,\bar{z})$, we obtain
	
	$$-2 \, {\rm Re} \, \Big(\overline{\mathcal{P}(z)}\tilde{g}_{1,\beta}(z)\Big)  \, \cdot u=2 \, {\rm Im} \, \Big(\overline{\mathcal{P}(z)}\tilde{g}_{1,\beta}(z)\Big)   \, \cdot \langle z,\bar{z}\rangle=0,$$
	whence $\overline{\mathcal{P}(z)}\tilde{g}_{1,\beta}(z)=0$, and therefore $\tilde{g}_{1,\beta}(z)=0$. That is, case 2.c.2 is impossible.

	\vspace{3ex}

	2.c.3) If $X_{1}$ falls into case 1.c, then
	  the components of the vector polynomials $\tilde{f}_{1,\alpha}$ and $\tilde{g}_{1,\beta}$ have the special form described in the analysis of case 1 in item 1.c. Namely, first, the components of the vector $\tilde{f}_{1,\alpha}$ are quadratic in $z$ or zero, and the components of the vector $\tilde{g}_{1,\beta}$ are linear in $z$ or zero. Hence the components of the vectors $\tilde{f}_{0,\alpha}$ and $\tilde{g}_{0,\beta}$ have degrees four and three, respectively, or are zero. Second, the terms of the field $\hat{X}_{1}$ fall either into case 1.c.1.1 or into case 1.c.1.2 (see above).
	
	It is clear that the vector $\tilde{f}_{0,\alpha}$ is equal to zero (since all terms of bidegrees $(4,1)$ and $(1,4)$ in the variables $z,\bar{z}$ in the lowest weighted component of the tangency condition can be obtained in a unique way, and only by means of this vector).
	
	Then the lowest weighted component of the tangency condition takes the form
	
	\begin{equation}\label{eq41}
	-2 \, {\rm Re} \, \Big( (\langle \tilde{f}_{1,\alpha}(z) ,\bar{z}\rangle +\overline{\mathcal{P}(z)}\tilde{g}_{1,\beta}(z)) \, \cdot (u+i \langle z,\bar{z}\rangle)+\overline{\mathcal{P}(z)}\tilde{g}_{0,\beta}(z)\Big)=0.
	\end{equation}
	
	Writing out the terms linear in $u$, we obtain
	$$-2 \, {\rm Re} \, \Big( (\langle \tilde{f}_{1,\alpha}(z) ,\bar{z}\rangle +\overline{\mathcal{P}(z)}\tilde{g}_{1,\beta}(z)) \, \cdot u\Big)=0,$$
	whence
	\begin{equation}\label{eq44}
    \langle \tilde{f}_{1,\alpha}(z) ,\bar{z}\rangle =-\mathcal{P}(z)\overline{\tilde{g}_{1,\beta}(z)}.
	\end{equation}

	Then \eqref{eq41} takes the form
	$$- \Big( (\langle \tilde{f}_{1,\alpha}(z) ,\bar{z}\rangle +\overline{\mathcal{P}(z)}\tilde{g}_{1,\beta}(z)) \, \cdot (i \langle z,\bar{z}\rangle)+\overline{\mathcal{P}(z)}\tilde{g}_{0,\beta}(z)\Big)-$$
	$$-\Big(\overline{\Big( \langle \tilde{f}_{1,\alpha}(z) ,\bar{z}\rangle +\overline{\mathcal{P}(z)}\tilde{g}_{1,\beta}(z)\Big)} \, \cdot (-i \langle z,\bar{z}\rangle)+\overline{\Big(\overline{\mathcal{P}(z)}\tilde{g}_{0,\beta}(z)\Big)}\Big)=$$
	$$=- \Big( 2 \,(\langle \tilde{f}_{1,\alpha}(z) ,\bar{z}\rangle +\overline{\mathcal{P}(z)}\tilde{g}_{1,\beta}(z)) \, \cdot (i \langle z,\bar{z}\rangle)+2 \, {\rm Re} \, (\overline{\mathcal{P}(z)}\tilde{g}_{0,\beta}(z))\Big)=0.$$
	
	Hence, taking into account the degrees of the coordinates of the vectors $\tilde{f}_{1,\alpha},\tilde{f}_{0,\alpha},\tilde{g}_{0,\beta}$, we obtain
	\begin{equation}\label{eq42}
	2 \,\langle \tilde{f}_{1,\alpha}(z) ,\bar{z}\rangle \, \cdot i \langle z,\bar{z}\rangle=-\overline{\mathcal{P}(z)}\tilde{g}_{0,\beta}(z).
	\end{equation}
	Next, let $\tilde{f}_{1,\alpha}(z)\frac{\partial}{\partial z}=\sum_{j=1}^{m}\tilde{f}_{1,\alpha,j}(z)\frac{\partial}{\partial z_{j}}$. Denote $\tilde{f}_{1,\alpha,j}$ by $\varphi_{j}$. Let $S$ be the set of those $s$ for which $\varphi_{s}\neq 0$ (and $S$ is nonempty, since we are considering a field of type 1.3). If $s\in S$, then the left-hand side of equality \eqref{eq42} contains the expression $\bar{z}_{s}^{2}$ with a nonzero factor. And if $s'\notin S$ and $s\in S$, then the expression $\varphi_{s}\cdot\bar{z}_{s} \, \cdot i \langle z,\bar{z}\rangle$ contains the monomial $\bar{z}_{s}\bar{z}_{s'}$ with a nonzero factor (moreover these monomials do not enter the remaining expressions $\varphi_{j}\cdot\bar{z}_{j} \, \cdot i \langle z,\bar{z}\rangle, \, j\neq s,$ and also note that no collecting of like terms is needed here). Therefore from \eqref{eq42} we obtain that some of the polynomials $p_{j}(z)$ contain the indicated monomials $\bar{z}_{s}^{2}, \bar{z}_{s}\bar{z}_{s'}$ with nonzero coefficients, and each monomial enters at least one polynomial. Hence $p_{j}(z)$ depend on all variables $z_{j}$, i.e. the defect $r$ is zero. This means that the terms of the field $\hat{X}_{1}$ fall into case 1.c.1.2.

	Since some of the polynomials $p_{j}$ contain nonzero monomials of the form $d_{s}z_{s}^{2}$ for $s\in S, \, d_{s}\in\mathbb{C}$, in Lemma 5 one can choose the monomials in the polynomials $p_{j}$ so that at least one of these monomials $d_{s_{0}}z_{s_{0}}^{2}$ is the distinguished monomial $\tilde{p}_{j}$ for some $p_{j}$ and enters exactly one polynomial $p_{j}$ (this choice does not affect the rest of the proof).
	Then from equality \eqref{eq44} we obtain that a nonzero monomial $d_{s_{0}}'z_{s_{0}}^{2}$ enters the polynomial $\varphi_{t}$ for some $t$. At the same time, depending on $t$, two cases are possible: either $t=s_{0}$ (i.e. $d_{s_{0}}'z_{s_{0}}^{2}$ enters $\varphi_{s_{0}}$, case 2.c.3.1), or $t\neq s_{0}$ (i.e. $d_{s_{0}}'z_{s_{0}}^{2}$ enters $\varphi_{t}, \,  t\neq s_{0}$, case 2.c.3.2).
	
	2.c.3.1) Since there exists $j$ such that $\tilde{p}_{j}=d_{s_{0}}z_{s_{0}}^{2}$, we have $\varphi_{s_{0}}\frac{\partial}{\partial z_{s_{0}}}(p_{j})=2d_{s_{0}}d_{s_{0}}'z_{s_{0}}^{3}+...$, where dots do not contain monomials divisible by $z_{s_{0}}^{3}$. Therefore, in order for the condition of case 1.c.1.2 to be fulfilled, there must exist $t\neq s_{0}$ such that $\varphi_{t}$ contains a nonzero monomial $d_{s_{0}}''z_{s_{0}}^{2}$ (since the equality $\tilde{f}_{1,\alpha}\frac{\partial}{\partial z}(p_{j})=0$ must hold, and hence the left-hand side of this equality must also contain the monomial $-2d_{s_{0}}d_{s_{0}}'z_{s_{0}}^{3}$).
	
	Thus in both cases 2.c.3.1 and 2.c.3.2 there exist $s$, $t$, and a nonzero number $d_{s}$ such that $d_{s}z_{s}^{2}$ enters $\varphi_{t}$ with $t\neq s$. At the same time we obtain that $\varphi_{t}\neq 0$, and hence there exists $j'$ such that a nonzero monomial of the form $d_{t}z_{t}^{2}$ enters $p_{j'}$.

	\vspace{3ex}

	Next, consider the expression
	$$F=-2 \, {\rm Re} \, \Big( (\langle \tilde{f}_{1,\alpha}(z) ,\bar{z}\rangle +\overline{\mathcal{P}(z)}\tilde{g}_{1,\beta}(z)) \, \cdot 2 \,i \, {\rm Re} \, (\mathcal{P}(z)\bar{y})\Big).$$
	
	It enters the expression
	$-2 \, {\rm Re} \, \Big( \langle f_{\alpha},\bar{z}\rangle+\mathcal{P}(z)\bar{g}_{\beta}\Big)$ (under condition \eqref{eq28}), which is contained in the tangency condition.
	
	Taking equality \eqref{eq44} into account, we obtain
	
	$$F=-2 \, {\rm Re} \, \Big( (\langle \tilde{f}_{1,\alpha}(z) ,\bar{z}\rangle - \langle \overline{\tilde{f}_{1,\alpha}(z)} ,z\rangle) \, \cdot 2 \,i \, {\rm Re} \, (\mathcal{P}(z)\bar{y})\Big)=$$
	$$=-2 \, (\langle \tilde{f}_{1,\alpha}(z) ,\bar{z}\rangle - \langle \overline{\tilde{f}_{1,\alpha}(z)} ,z\rangle) \, \cdot 2 \,i \, {\rm Re} \, (\mathcal{P}(z)\bar{y}).$$
	
	The expression $F$ contains the terms $2 \,i \, \langle \overline{\tilde{f}_{1,\alpha}(z)} ,z\rangle \, \cdot   (\mathcal{P}(z)\bar{y})$, which, in turn, contain the terms $2 \,i \, \overline{\varphi_{t}}z_{t}  p_{j}\bar{y}_{j}$. By what has been proved, $\varphi_{t}$ contains a nonzero monomial $d_{s} z_{s}^{2}$ with $s\neq t$ for some $s$ and $t$, and there exists $j$ such that $p_{j}$ contains a nonzero monomial $d_{t}z_{t}^{2}$. Then $F$ contains (after collecting like terms) the nonzero monomial $\sigma_{1}=2 \,i \, \bar{d}_{s} \bar{z}_{s}^{2}z_{t}  d_{t}z_{t}^{2}\bar{y}_{j}=d \bar{z}_{s}^{2}  z_{t}^{3}\bar{y}_{j}$ for some nonzero constant $d$. We note that in the case 2.c.3 under consideration the order $q$ of the field $X$ is equal to five.

	At the same time, the term $-\sigma_{1}$, not entering $F$, also appears in \eqref{eq36}. We write out the conditions under which this term can occur; as in case 1, we formulate them for terms $\Theta$ of the field $X$ and terms $\chi$ in the defining relation such that the expression $\Theta(\chi)$ contains the term $-\sigma_{1}$.

	Before listing these conditions, let us make three remarks. We number the remarks in the same way as the logical subcases.

	\vspace{3ex}
	
	\textit{Remarks.}
	
	2.c.3.3) First, it is enough to consider only such $\chi$ and $\Theta$ that $\chi$ is a monomial and $\Theta$ is equal to one of the following expressions:
	$\theta(z,y,w)\frac{\partial}{\partial z_{\lambda}}, \, \theta(\bar{z},\bar{y},\bar{w})\frac{\partial}{\partial \bar{z}_{\lambda}}, \,\theta(z,y,w)\frac{\partial}{\partial y_{\lambda}}, \,\theta(\bar{z},\bar{y},\bar{w})\frac{\partial}{\partial \bar{y}_{\lambda}}, \, \theta(z,y,w)\frac{\partial}{\partial w}$ or $\theta(\bar{z},\bar{y},\bar{w})\frac{\partial}{\partial \bar{w}}$, where $\theta$ is a monomial.

	\vspace{3ex}
	
	2.c.3.4) Second, $\Theta$ cannot be equal to
	$\theta(z,y,w)\frac{\partial}{\partial w}$ or $\theta(\bar{z},\bar{y},\bar{w})\frac{\partial}{\partial \bar{w}}$.

	Indeed, suppose the contrary. By the rigidity condition for the defining equation (i.e. the defining equation does not depend on the variable $u$), we have $\chi=\frac{1}{2i}w$ or, respectively, $\chi=-\frac{1}{2i}\bar{w}$.

	Let $\chi=\frac{1}{2i}w, \, \Theta=\theta(z,y,w)\frac{\partial}{\partial w}$. As was already shown in the proof of Lemma 12, the coefficient $\theta$ cannot depend on the variable $y$, since in this case monomials of the form $d\cdot z^{a_{1}}y^{a_{2}}u^{\delta}$ with $a_{2}\neq 0$ appear in the tangency condition, and these monomials cannot occur in any other way except by means of $\Theta$ and $\chi$ of the indicated form.

	Then $\Theta=d\cdot z^{a_{1}}w^{\delta}\frac{\partial}{\partial w}$ for some $\delta$. The term $-\sigma_{1}$ is one of the terms obtained if one of the factors $w=u+iv$ in the expression $w^{\delta}$ is replaced by the polynomial $ip_{j}\bar{y}_{j}$ entering the expression $iv$, and the remaining factors are replaced by $i\langle z,\bar{z}\rangle$.
	
	At the same time it is clear that $w^{\delta}$ cannot be replaced by any other expression entering it, since this expression must depend linearly on the variable $\bar{y}_{j}$ and not depend on the other variables $\bar{y}_{s}$ or the variable $y$, and terms with this property can be contained only in the expression $i^{\delta}p_{j}\bar{y}_{j}(\langle z,\bar{z}\rangle)^{\delta-1}$.
	
	But then the order of the field $\Theta$ must be equal to $q=5$ (since the weight of the monomial $\sigma_{1}$ is six, and all variables $w$ in the expression $w^{\delta}$ in $\Theta$ can be replaced by $i \langle z,\bar{z}\rangle$). Hence $\Theta$ enters the supporting component $\hat{X}$ of the field $X$. This is impossible, since we are considering a field $X$ of the first type.

	The case $\chi=-\frac{1}{2i}\bar{w}$, in which $\Theta=d\cdot \bar{z}^{a_{1}}\bar{w}^{\delta}\frac{\partial}{\partial \bar{w}}$, is treated completely analogously.

	\vspace{3ex}

	2.c.3.5) Third, if the coefficients of the field $\Theta=\theta(z,y,w)\frac{\partial}{\partial z_{\lambda}}$ (where $\theta=d\cdot z^{a}y^{b}w^{c}$ is a monomial) depend on $w=u+i\Big(|z_{1}|^{2}+|z_{2}|^{2}+...+|z_{m}|^{2}+ 2 \, {\rm Re} \,\Big(p_{1}\bar{y}_{1}+...+p_{k}\bar{y}_{k})\Big)+o(3)\Big)$, then, in order to obtain the term $-\sigma_{1}$, one must replace each of the factors $w$ in the expression $w^{c}$ in $\theta$ by one of the
	expressions $i |z_{\lambda}|^{2}, 2\,i \, {\rm Re} \, (p_{\lambda}\bar{y}_{\lambda})$ or $o(3)$. Moreover, none of the factors $w$ can be replaced by $o(3)$, since otherwise the order of the field $\Theta$ would be less than $q=5$, which is impossible by the definition of the quantity $q$ (and it is clear that, to obtain an upper estimate for the order of the field, one can put $c=1$ in the expression for $\theta$; in this case, to estimate the order of the field one must subtract the weight of a term of the form $o(3)$ from the weight of the term $-\sigma_{1}$ and add the weight of the term $i |z_{\lambda}|^{2}$, i.e. one obtains a quantity not exceeding $6-4+2=4$). It also cannot be replaced by $2 \,i\, {\rm Re} \, (p_{\lambda}\bar{y}_{\lambda})$, since then the order of the field $\Theta$ does not exceed $q=5$ (it is also clear that, to obtain an upper estimate for the order of the field, one can put $c=1$ in the expression for $\theta$; in this case, to estimate the order of the field one must subtract the weight of a term of the form $2 \,i\, {\rm Re} \, (p_{\lambda}\bar{y}_{\lambda})$ from the weight of the term $-\sigma_{1}$ and add the weight of the term $i \,|z_{\lambda}|^{2}$, i.e. $6-3+2=5$). In this case the order of the field is equal to $q$ (by the definition of the quantity $q$), whence $c=1$ (since we are within case 2). But such terms $\Theta$ have already been taken into account, i.e. $\Theta$ and $\chi$ are such that $-\sigma_{1}$ enters $F$.

	Hence the variable $w$ can be replaced only by the factor $i \,|z_{\lambda}|^{2}$.
	
	The cases $\Theta=\theta(z,y,w)\frac{\partial}{\partial y_{\lambda}}, \, \Theta=\theta(\bar{z},\bar{y},\bar{w})\frac{\partial}{\partial \bar{z}_{\lambda}}, \, \Theta=\theta(\bar{z},\bar{y},\bar{w})\frac{\partial}{\partial \bar{y}_{\lambda}}$ (where $\theta$ is a monomial) are treated analogously.

	\vspace{3ex}
	
	Now let us list the conditions on $\Theta$ and $\chi$ (below, to obtain the term $-\sigma_{1}$, the variable $w$ is replaced by one of the terms $i \,|z_{\lambda}|^{2}$ in the expression $u+iv$, where $v$ is given by formula \eqref{eq28}):

	\vspace{3ex}
	
	2.c.3.a) The field $X$ contains a term of the form $\Theta=d_{3}z^{a_{1}}w^{a_{4}}\frac{\partial}{\partial y_{\lambda}}$, and the defining function contains a term of the form $\chi=d_{4}z^{a_{2}} \bar{z}^{a_{3}} y_{\lambda} \cdot \bar{y}_{j}$.
	
	By property 5 of the normal form we have $|a_{2}|=|a_{3}|=1$. This means that $a_{4}=1$. But then $w$ must be replaced by one of the terms $i \,|z_{\lambda'}|^{2}$. Since $\sigma_{1}$ is not divisible by the monomial $|z_{\lambda'}|^{2}$ for any $\lambda'$, condition 2.c.3.a cannot be fulfilled.

	\vspace{3ex}

	2.c.3.b.1) The field $X$ contains terms of the form $\Theta=d_{3}z^{a_{1}}w^{a_{4}}\frac{\partial}{\partial z_{\lambda}}$, and the defining function contains terms of the form $\chi=d_{4}z^{a_{2}}\bar{z}^{a_{3}} \cdot \bar{y}_{j}$, with $-\sigma_{1}=d_{4}\Theta(z^{a_{2}}) \bar{z}^{a_{3}} \cdot \bar{y}_{j}$.
	
	Since the defining function depends quadratically on $z,\bar{z}$, the pair $(|a_{2}|,|a_{3}|)$ can take only the values $(0,2),(1,1),(2,0)$. By conditions 3 and 4 for the normal form, $|a_{2}|=2, \, |a_{3}|=0$, whence $a_{4}=2$. But then $w^{2}$ must be replaced by one of the terms $i \,|z_{\lambda'}|^{2}\cdot i \,|z_{\lambda''}|^{2}$. Since $\sigma_{1}$ is not divisible by the monomial $|z_{\lambda'}|^{2}$ for any $\lambda'$, condition 2.c.3.b.1 cannot be fulfilled.

	\vspace{3ex}

	2.c.3.b.2) The field $X$ contains terms of the form $\Theta=d_{3}\bar{z}^{a_{1}}\bar{w}^{a_{4}}\frac{\partial}{\partial \bar{z}_{\lambda}}$, and the defining function contains terms of the form $\chi=d_{4}z^{a_{2}}\bar{z}^{a_{3}} \cdot \bar{y}_{j}$, with $-\sigma_{1}=d_{4}z^{a_{2}} \Theta(\bar{z}^{a_{3}}) \cdot \bar{y}_{j}$ (and hence $|a_{3}|>0$).
	
	By condition 3 for the normal form, such $\chi$ are absent from the defining relation. Hence case 2.c.3.b.2 is impossible.

	\vspace{3ex}

	2.c.3.c) The field $X$ contains terms of the form $\Theta=d_{3}\bar{z}^{a_{1}} \bar{w}^{a_{4}}  \bar{y}_{j} \frac{\partial}{\partial \bar{z}_{\lambda}}$, and the defining function contains terms of the form $\chi=d_{4}z^{a_{2}}\bar{z}^{a_{3}}$.
	
	Since the defining function depends quadratically on $z,\bar{z}$, the pair $(|a_{2}|,|a_{3}|)$ can take only the values $(0,2),(1,1),(2,0)$.	
	By conditions 3 and 4 for the normal form, $|a_{2}|=|a_{3}|=1$. Therefore $a_{4}=2$ (since $\sigma_{1}$ is divisible by $z_{t}^{3}$). But then $w^{2}$ must be replaced by one of the terms $i \,|z_{\lambda'}|^{2}\cdot i \,|z_{\lambda''}|^{2}$. Since $\sigma_{1}$ is not divisible by the monomial $|z_{\lambda'}|^{2}$ for any $\lambda'$, condition 2.c.3.c cannot be fulfilled.

	\vspace{3ex}

	2.c.3.d) 
	The field $X$ contains a term of the form $\Theta=d_{3}\bar{z}^{a_{1}}\bar{w}^{a_{4}}\frac{\partial}{\partial \bar{y}_{\lambda}}$, and the defining function contains a term of the form $\chi=d_{4}z^{a_{2}} \bar{z}^{a_{3}} \bar{y}_{j}\bar{y}_{\lambda}$.
	
	By property 6 of the normal form, the defining equation does not contain such terms. Hence case 2.c.3.d is impossible.

	\vspace{3ex}

	2.c.3.e) The field $X$ contains terms of the form $\Theta=d_{3}\bar{z}^{a_{1}} \bar{w}^{a_{4}} \bar{y}_{j} \frac{\partial}{\partial \bar{y}_{\lambda}}$, and the defining function contains terms of the form $\chi=d_{4}z^{a_{2}} \bar{z}^{a_{3}}  \cdot \bar{y}_{\lambda}$.
	
	By properties 3 and 4 of the normal form we have $|a_{2}|=2, \, |a_{3}|=0$. This means that $a_{4}=1$. But then $w$ must be replaced by one of the terms $i \,|z_{\lambda'}|^{2}$. Since $\sigma_{1}$ is not divisible by the monomial $|z_{\lambda'}|^{2}$ for any $\lambda'$, condition 2.c.3.e cannot be fulfilled.

	\vspace{3ex}

	Thus nonzero fields for case 2.c.3 do not exist. Nonzero fields in case 2 can occur only in item 2.c.1.

	\vspace{3ex}

	\textit{Case 3.} Recall that the fields $X_{\lambda}$ are defined by the formula $X_{1}=[X_{0},X], \, X_{\lambda}=[X_{0},X_{\lambda-1}], \, \lambda\geq 2$. At the same time $X_{\lambda}$ has order $\tilde{q}=q-2\lambda$, and the field $\hat{X}_{\lambda}$ has the form
	
	$$\hat{X}_{\lambda}=2 \, {\rm Re} \, \Big(\sum_{\nu}\Big[\frac{\partial^{\lambda}}{\partial w^{\lambda}}(f_{\nu,\alpha_{\nu}})\Big]_{\tilde{q}-1}\frac{\partial}{\partial z_{\nu}}+\sum_{\mu}\Big[\frac{\partial^{\lambda}}{\partial w^{\lambda}}(g_{\mu,\beta_{\mu}})\Big]_{\tilde{q}-2}\frac{\partial}{\partial y_{\mu}}+\Big[\frac{\partial^{\lambda}}{\partial w^{\lambda}}(h_{\gamma})\Big]_{\tilde{q}}\frac{\partial}{\partial w}\Big).$$

	If the highest degree in $w$ of the coefficients of the field $\hat{X}$ is equal to $\nu>1$, then the highest degree in $w$ of the coefficients of the field $\hat{X}_{\nu-1}$ is equal to one, and $\hat{X}_{\nu-1}$ belongs to ${\rm st} \, M_{0}$ (the proof is analogous to the proof of Lemma 11; one only needs to replace $\hat{X}_{1}$ by $\hat{X}_{\nu-1}$). Thus $\hat{X}_{\nu-1}$ falls into case 2. Therefore case 3 is impossible if $\hat{X}_{\nu-1}$ does not fall into case 2.c.1.
	
	\vspace{3ex}
	
	3.a) Let $X$ be such that $\hat{X}_{\nu-1}$ falls into case 2.c.1. We shall prove that $\hat{X}=0$ for $\nu\geq 3$ (i.e. nonzero fields of the first type whose supporting component depends on $w$ can only be fields of type 1.1(1) and 1.1(2)). It is enough to consider the case $\nu=3$ (if $\nu>3$, then one can consider the field $X_{\nu-3}$).

	Thus let $\hat{X}=2 \, {\rm Re} \, \Big( (\tilde{f}_{3,\alpha}w^{3}+\tilde{f}_{2,\alpha}w^{2}+\tilde{f}_{1,\alpha}w+\tilde{f}_{0,\alpha})\frac{\partial}{\partial z}+(\tilde{g}_{3,\beta}w^{3}+\tilde{g}_{2,\beta}w^{2}+\tilde{g}_{1,\beta}w+\tilde{g}_{0,\beta})\frac{\partial}{\partial y}\Big)$. At the same time, by our choice of conditions, $\tilde{g}_{3,\beta}=0$; the components of the vector polynomial $\tilde{f}_{3,\alpha}$ are linear in $z$ or zero; the components of the vector polynomial $\tilde{f}_{2,\alpha}$ have degree three in $z$ or are zero; the components of the vector polynomial $\tilde{f}_{1,\alpha}$ have degree five in $z$ or are zero; the components of the vector polynomial $\tilde{f}_{0,\alpha}$ have degree seven in $z$ or are zero; the components of the vector polynomial $\tilde{g}_{2,\beta}$ have degree two in $z$ or are zero; the components of the vector polynomial $\tilde{g}_{1,\beta}$ have degree four in $z$ or are zero; and the components of the vector polynomial $\tilde{g}_{0,\beta}$ have degree six in $z$ or are zero. The lowest component of the tangency condition has the form

	\begin{equation}\label{eq49}
	-2 \, {\rm Re} \, \Big(\langle \tilde{f}_{3,\alpha}(z) ,\bar{z}\rangle  \, \cdot (u+i \langle z,\bar{z}\rangle)^{3}+(\langle \tilde{f}_{2,\alpha}(z) ,\bar{z}\rangle +\overline{\mathcal{P}(z)}\tilde{g}_{2,\beta}(z))  \, \cdot (u+i \langle z,\bar{z}\rangle)^{2}+
	\end{equation}
	$$+(\langle \tilde{f}_{1,\alpha}(z) ,\bar{z}\rangle+\overline{\mathcal{P}(z)}\tilde{g}_{1,\beta}(z))  \, \cdot (u+i \langle z,\bar{z}\rangle)+\langle \tilde{f}_{0,\alpha}(z) ,\bar{z}\rangle+\overline{\mathcal{P}(z)}\tilde{g}_{0,\beta}(z)\Big)=0.$$
	
	Separating in \eqref{eq49} the terms of bidegrees $(7,1),(5,1),(3,1),(6,2),(4,2)$ in the variables $(z,\bar{z})$, we obtain $\tilde{f}_{0,\alpha}(z)=\tilde{f}_{1,\alpha}(z)=\tilde{f}_{2,\alpha}(z)=\tilde{g}_{0,\beta}(z)=\tilde{g}_{1,\beta}(z)=0$. Equality \eqref{eq49} takes the form
	\begin{equation}\label{eq50}
	-2 \, {\rm Re} \, \Big(\langle \tilde{f}_{3,\alpha}(z) ,\bar{z}\rangle  \, \cdot (u+i \langle z,\bar{z}\rangle)^{3}+(\overline{\mathcal{P}(z)}\tilde{g}_{2,\beta}(z))  \, \cdot (u+i \langle z,\bar{z}\rangle)^{2}\Big)=0.
	\end{equation}
	Separating in \eqref{eq50} the terms of degrees $3,2,1,0$ in the variable $u$, we obtain a full-rank system with respect to the variables $\langle\tilde{f}_{3,\alpha}(z) ,\bar{z}\rangle, \langle \overline{\tilde{f}_{3,\alpha}(z)} ,z\rangle, \overline{\mathcal{P}(z)}\tilde{g}_{2,\beta}(z), \mathcal{P}(z)\overline{\tilde{g}_{2,\beta}(z)}$:
	
	\begin{equation}\label{eq51}
	6\langle \tilde{f}_{3,\alpha}(z) ,\bar{z}\rangle + 6\langle \overline{\tilde{f}_{3,\alpha}(z)} ,z\rangle =0,
	\end{equation}
	
	$$(6i \langle \tilde{f}_{3,\alpha}(z) ,\bar{z}\rangle-6i \langle \overline{\tilde{f}_{3,\alpha}(z)} ,z\rangle)  \, \cdot \langle z,\bar{z}\rangle+2 \overline{\mathcal{P}(z)}\tilde{g}_{2,\beta}(z)+2 \mathcal{P}(z)\overline{\tilde{g}_{2,\beta}(z)} =0,$$
	
	$$(-3\langle \tilde{f}_{3,\alpha}(z) ,\bar{z}\rangle -3\langle \overline{\tilde{f}_{3,\alpha}(z)} ,z\rangle) \, \cdot (\langle z,\bar{z}\rangle)^{2}+(2i \overline{\mathcal{P}(z)}\tilde{g}_{2,\beta}(z)-2i \mathcal{P}(z)\overline{\tilde{g}_{2,\beta}(z)})  \, \cdot  \langle z,\bar{z}\rangle=0,$$
	
	$$(-i \langle \tilde{f}_{3,\alpha}(z) ,\bar{z}\rangle +i \langle \overline{\tilde{f}_{3,\alpha}(z)} ,z\rangle) \, \cdot (\langle z,\bar{z}\rangle)^{3}+(- \overline{\mathcal{P}(z)}\tilde{g}_{2,\beta}(z)- \mathcal{P}(z)\overline{\tilde{g}_{2,\beta}(z)})  \, \cdot  (\langle z,\bar{z}\rangle)^{2}=0.$$
	
	In more detail: after simplification this system has the form
	
	$$\langle \overline{\tilde{f}_{3,\alpha}(z)} ,z\rangle=-\langle \tilde{f}_{3,\alpha}(z) ,\bar{z}\rangle, \overline{\mathcal{P}(z)}\tilde{g}_{2,\beta}(z)= \mathcal{P}(z)\overline{\tilde{g}_{2,\beta}(z)},$$
	
	$$\mathcal{P}(z)\overline{\tilde{g}_{2,\beta}(z)}=-3i \langle \tilde{f}_{3,\alpha}(z) ,\bar{z}\rangle \langle z,\bar{z}\rangle, \langle \tilde{f}_{3,\alpha}(z) ,\bar{z}\rangle =0,$$
	i.e. all its solutions are zero. From the equality $\langle \tilde{f}_{3,\alpha}(z) ,\bar{z}\rangle =0$ we obtain $\tilde{f}_{3,\alpha}(z)=0$ -- a contradiction. Thus such fields do not exist.	
	
	\vspace{3ex}

	3.b) Next, let $\nu=2$, i.e. the field has type 1.1(2). We have $\hat{X}=2 \, {\rm Re} \, \Big((\tilde{f}_{2,\alpha}w^{2}+\tilde{f}_{1,\alpha}w+\tilde{f}_{0,\alpha})\frac{\partial}{\partial z}+(\tilde{g}_{2,\beta}w^{2}+\tilde{g}_{1,\beta}w+\tilde{g}_{0,\beta})\frac{\partial}{\partial y}\Big)$. At the same time, by our choice of conditions, $\tilde{g}_{2,\beta}=0$; the components of the vector polynomial $\tilde{f}_{2,\alpha}$ are linear in $z$ or zero; the components of the vector polynomial $\tilde{f}_{1,\alpha}$ have degree three in $z$ or are zero; the components of the vector polynomial $\tilde{f}_{0,\alpha}$ have degree five in $z$ or are zero; the components of the vector polynomial $\tilde{g}_{1,\beta}$ have degree two in $z$ or are zero; and the components of the vector polynomial $\tilde{g}_{0,\beta}$ have degree four in $z$ or are zero. As in the case $\nu=3$, we obtain $\tilde{f}_{0,\alpha}(z)=\tilde{f}_{1,\alpha}(z)=\tilde{g}_{0,\beta}(z)=0$ (the proof is analogous).
	
	Let $X'$ be another field of type 1.1(2). We have $\hat{X'}=2 \, {\rm Re} \, \Big(\tilde{f}_{2,\alpha}'w^{2}\frac{\partial}{\partial z}+\tilde{g}_{1,\beta}'w\frac{\partial}{\partial y}\Big)$.
	
	The fields $\hat{X}$ and $\hat{X}'$ do not contain terms $d \, \bar{w}^{c}\frac{\partial}{\partial \bar{z}_{s}}$ and $d \, \bar{w}^{c}\frac{\partial}{\partial \bar{y}_{s}}$. Therefore, by Lemma 12, the coefficient at $\frac{\partial}{\partial z}$ of the field $\widehat{[X,X']}$ is equal to
	
	\begin{equation}\label{eq55}
	w^{4}\Big(\Big(\tilde{f}_{2,\alpha}\frac{\partial}{\partial z}(\tilde{f}_{2,\alpha}')\Big)\frac{\partial}{\partial z}-\Big(\tilde{f}_{2,\alpha}'\frac{\partial}{\partial z}(\tilde{f}_{2,\alpha})\Big)\frac{\partial}{\partial z}\Big),
	\end{equation}
	if this expression is nonzero. Suppose that expression \eqref{eq55} is nonzero. Since $X$ and $X'$ are fields of the first type, the coefficients at $\frac{\partial}{\partial w}$ of the fields $\hat{X},\hat{X}'$ are zero. Then by Lemma 12 the field $[X,X']$ is also a field of the first type. At the same time $\widehat{[X,X']}$ depends on $w^{4}$.
	But by what was proved above in subitems 3 and 3.a, the field $\widehat{[X,X']}$ depends on $w$ at most quadratically -- a contradiction; hence expression \eqref{eq55} must be equal to zero. Therefore the fields $\tilde{f}_{2,\alpha}\frac{\partial}{\partial z}$ and $\tilde{f}_{2,\alpha}'\frac{\partial}{\partial z}$ commute. At the same time $2 \, {\rm Re} \,(\tilde{f}_{2,\alpha}\frac{\partial}{\partial z})=\frac{1}{2}\hat{X}_{2}, \, 2 \, {\rm Re} \,(\tilde{f}_{2,\alpha}'\frac{\partial}{\partial z})=\frac{1}{2}\hat{X}_{2}'$. Since the fields $\hat{X}_{2},\hat{X}_{2}'$ belong to the automorphism algebra, the fields $2 \, {\rm Re} \,(\tilde{f}_{2,\alpha}\frac{\partial}{\partial z}), \, 2 \, {\rm Re} \,(\tilde{f}_{2,\alpha}'\frac{\partial}{\partial z})$ are linear combinations with real coefficients of the fields $2 \, {\rm Re} \,(iz_{j}\frac{\partial}{\partial z_{j}}), \, 2 \, {\rm Re} \,(z_{j}\frac{\partial}{\partial z_{s}}-z_{s}\frac{\partial}{\partial z_{j}}), \, 2 \, {\rm Re} \,(iz_{j}\frac{\partial}{\partial z_{s}}+iz_{s}\frac{\partial}{\partial z_{j}})$. Therefore
	the matrices, which correspond to the fields $2 \, {\rm Re} \,(\tilde{f}_{2,\alpha}\frac{\partial}{\partial z})$ and $2 \, {\rm Re} \,(\tilde{f}_{2,\alpha}'\frac{\partial}{\partial z})$  are skew-Hermitian. Since the fields $2 \, {\rm Re} \,(\tilde{f}_{2,\alpha}\frac{\partial}{\partial z})$ and $2 \, {\rm Re} \,(\tilde{f}_{2,\alpha}'\frac{\partial}{\partial z})$ also commute, there exists a unitary transformation, which diagonalizes the matrices of these fields. Hence in the same basis the fields $\hat{X}_{2},\hat{X}_{2}'$ are equal to a linear combination with real coefficients of the fields $2 \, {\rm Re} \, (iz_{1}\frac{\partial}{\partial z_{1}}),...,2 \, {\rm Re} \, (iz_{m}\frac{\partial}{\partial z_{m}})$, which we shall call \textit{diagonal} (if the field $\hat{X_{\nu}}$ is equal to such a linear combination, we shall say that the field $X$ \textit{is diagonalized}). Hence we obtain that there are no more than $m$ such linearly independent fields $X$.

	\vspace{3ex}
	
	Next, if nonzero fields $X$ of type 1.1(2) exist, consider a field $Y$ of type 1.1(1). We have $\hat{Y}=2 \, {\rm Re} \, \Big(\tilde{f}_{1,\alpha}w\frac{\partial}{\partial z}+\tilde{g}_{0,\beta}\frac{\partial}{\partial y}\Big)$ for some $\tilde{f}_{1,\alpha},\tilde{g}_{0,\beta}$.
	
	By Lemma 12, the coefficient at $\frac{\partial}{\partial z}$ of the field $\widehat{[X,Y]}$ is equal to
	
	$$w^{3}\Big(\Big(\tilde{f}_{2,\alpha}\frac{\partial}{\partial z}(\tilde{f}_{1,\alpha})\Big)\frac{\partial}{\partial z}-\Big(\tilde{f}_{1,\alpha}\frac{\partial}{\partial z}(\tilde{f}_{2,\alpha})\Big)\frac{\partial}{\partial z}\Big),$$
	if this expression is nonzero.
	But if this coefficient is nonzero, then the field $[X,Y]$ is a field of the first type (the reasoning is analogous to that given above for the fields $X,X'$), whence by what was proved in subitems 3 and 3.a this coefficient must be equal to zero (since the degree of this coefficient in the variable $w$ is greater than two) -- a contradiction. Hence the fields $\tilde{f}_{2,\alpha}\frac{\partial}{\partial z}$ and $\tilde{f}_{1,\alpha}\frac{\partial}{\partial z}$ commute. Therefore in some basis all fields $X$ of type 1.1(2) and $Y$ of type 1.1(1) are diagonalized, i.e. $\hat{X_{2}}$ and $\hat{Y_{1}}$ are equal to a linear combination with real coefficients of the fields $2 \, {\rm Re} \, (iz_{1}\frac{\partial}{\partial z_{1}}),...,2 \, {\rm Re} \, (iz_{m}\frac{\partial}{\partial z_{m}})$.
	
	Next, for any other field $Y'=2 \, {\rm Re} \, \Big((\tilde{f}_{1,\alpha}'w)\frac{\partial}{\partial z}+(\tilde{g}_{0,\beta}')\frac{\partial}{\partial y}\Big)$ of the same type 1.1(1), we have that $\hat{Y}_{1}'$ commutes with all fields $\hat{X}_{2}$ for fields $X$ of type 1.1(2). At the same time, by Lemma 12, the coefficient at $\frac{\partial}{\partial z}$ of the field $\widehat{[Y,Y']}$ is equal to
	
	$$F=w^{2}\Big(\Big(\tilde{f}_{1,\alpha}\frac{\partial}{\partial z}(\tilde{f}_{1,\alpha}')\Big)\frac{\partial}{\partial z}-\Big(\tilde{f}_{1,\alpha}'\frac{\partial}{\partial z}(\tilde{f}_{1,\alpha})\Big)\frac{\partial}{\partial z}\Big),$$
	if $F\neq 0$. But if $F\neq 0$, then the field $[Y,Y']$ is a field of the first type (the proof is analogous to that given above for the fields $X,X'$). Moreover, the weight of the vector coefficient $F$ is five, and hence the weight of the vector coefficient $G$ at $\frac{\partial}{\partial y}$ of the field $\widehat{[Y,Y']}$ is four (if $G\neq 0$). $G$ cannot depend on $w$ quadratically, since otherwise $G$ contains terms of the form $d \, w^{2} \frac{\partial}{\partial y_{s}}$, but in this case the field $\widehat{[Y,Y']}$ also contains nonzero terms of the form $\tilde{d} \, w^{2} \tilde{p}_{s} \frac{\partial}{\partial w}$ in the supporting component, and then $\widehat{[Y,Y']}$ cannot be a field of the first type. Therefore $G$ depends on $w$ at most linearly.
	That is, if $F\neq 0$, then the field $\widehat{[Y,Y']}$ has type 1.1(2), and therefore it is diagonalizable simultaneously with all fields of type 1.1(2) and the field $Y$. If the field $\hat{Y}_{1}'$ contains nondiagonal terms in the basis in which all fields $X$ of type 1.1(2) and the field $Y$ are diagonalized, then the field $2 \, {\rm Re} \, (F\frac{\partial}{\partial z})$ is also a sum of nondiagonal terms, whence $F=0$ -- a contradiction. Therefore either $Y'$ can also be diagonalized simultaneously with all $X$ and the field $Y$, or $F=0$, from which the same follows. That is, all fields $Y$ can also be diagonalized simultaneously with the fields $X$ in some basis, and therefore there are also no more than $m$ linearly independent fields $Y$.
	
	If, however, nonzero fields $X$ of type 1.1(2) do not exist, then $F=0$, since the degree of $F$ in $w$ is two. This means that all fields $Y$ of type 1.1(1) are diagonalizable in some basis (the proof is analogous), and therefore there are no more than $m$ linearly independent fields $Y$.

    \vspace{3ex}
	
	4) Let $X$ be a field of the second type. By Lemma 10 and the weighted homogeneity of the polynomials $f_{\nu,\alpha_{\nu}},g_{\mu,\beta_{\mu}},h_{\gamma}$, we have $f_{\alpha}=\tilde{f}_{1,\alpha}(z)\cdot w^{a_{1}}+\tilde{f}_{2,\alpha}(z)\cdot w^{a_{1}-1}+...+\tilde{f}_{a_{1}+1,\alpha}(z), \, g_{\beta}=\tilde{g}_{1,\beta}(z)\cdot w^{a_{2}}+\tilde{g}_{2,\beta}(z)\cdot w^{a_{2}-1}+...+\tilde{g}_{a_{2}+1,\beta}(z),  \,
	h_{\gamma}=\tilde{h}_{1,\gamma}(z)\cdot w^{a_{3}}+\tilde{h}_{2,\gamma}(z)\cdot w^{a_{3}-1}+...+\tilde{h}_{a_{3}+1,\gamma}(z)$, where dots denote monomials of lower degree in the variable $w$; the coordinates of the vectors $\tilde{f}_{1,\alpha},\tilde{f}_{2,\alpha}$ are homogeneous polynomials of degrees $b_{1},b_{1}+2$, respectively; and the functions $\tilde{h}_{1,\gamma}, \tilde{h}_{2,\gamma}$ are homogeneous polynomials of degrees $b_{3}$ and $b_{3}+2$, respectively, with $b_{1}+1+2a_{1}=b_{3}+2a_{3}$. Everywhere below we assume that $\tilde{h}_{1,\gamma}\neq 0$ (otherwise $h_{\gamma}=0$, and we obtain a field of the first type, contradicting our assumption). We also assume that $\tilde{f}_{1,\alpha}\neq 0$ if $f_{\alpha}\neq 0$, and that $\tilde{g}_{1,\beta}\neq 0$ if $g_{\beta}\neq 0$.

	\vspace{3ex}

	We split the consideration of fields of the second type into five steps.
	
	4.a) Fields of type $2.0(\geq 3)$ do not exist.
	
	4.b) Fields of type $2.\geq 3(\geq 0)$ do not exist.
	
	4.c) Fields of type $2.2(\geq 2)$ do not exist.
	
	4.d) Fields of type $2.1(\geq 3)$ do not exist.
	
	4.e) We estimate the number of linearly independent fields of the remaining types ($2.b(0), \, 2.b(1), \, 2.0(2), \, 2.1(2)$).

	\vspace{3ex}

	4.a)
	Let $X$ be a field of type $2.0(\geq 3)$ such that its order $q$ takes the maximal value. We note that the order $q$ of the field $X$ is equal to $2 a_{3}+b_{3}$ (the weight of the polynomial $h_{\gamma}$). We shall prove that if $X$ is not equal to zero, then the quantity $q$ can be increased (i.e. there exists a field of the same type $2.0(\geq 3)$ whose order is greater than $2a_{3}+b_{3}$).
	
	Since $b_{1}+1+2a_{1}=b_{3}+2a_{3}$, we have $2a_{1}=b_{3}-b_{1}-1+2a_{3}=-b_{1}-1+2a_{3}$, i.e. $a_{1}<a_{3}$. At the same time it is clear that if $\tilde{f}_{1,\alpha}\neq 0$, then $b_{1}$ is odd, whence $b_{1}\geq 1$. Therefore, by Lemma 12, the coefficient at $\frac{\partial}{\partial w}$ of the field $\widehat{[X,X_{1}]}$ is equal to

	$$H=-w^{2a_{3}-2}a_{3}\tilde{h}_{1,\gamma}^{2}
	+...,$$
	where dots denote terms of lower degree in the variable $w$, if $H\neq 0$. Therefore, if $\tilde{h}_{1,\gamma}\neq 0$, then the order of the field $\widehat{[X,X_{1}]}$ is equal to $2(2a_{3}-2)+2b_{3}=4a_{3}-4=2a_{3}+(2a_{3}-4)\geq 2a_{3}+2>2a_{3}+b_{3}$, i.e. the order of the field has increased. At the same time the field $\widehat{[X,X_{1}]}$ has type $2.0(2a_{3}-2)$, and hence falls into the same type $2.0(\geq 3)$, since $2a_{3}-2\geq 4$. The contradiction obtained (with the maximality of the order $q$) means that $\tilde{h}_{1,\gamma}= 0$. Therefore nonzero fields of type $2.0(\geq 3)$ do not exist.
		
		\vspace{3ex}
	
	4.b)
	Let $X$ be a field of type $2.\geq 3(\geq 0)$.
	
	In this case $b_{3}\geq 3$, and then the only term of bidegree $(b_{3},0)$ in the variables $(z,\bar{z})$ and of degree $a$ in the variable $u$ in the lowest component of the tangency condition is $\tilde{h}_{1,\gamma}u^{a}$. This means that $\tilde{h}_{1,\gamma}=0$ -- a contradiction.

	Hence fields of type $2.\geq 3(\geq 0)$ do not exist.

	\vspace{3ex}

	4.c) Now we prove that fields of type $2.2(\geq 2)$ do not exist.		
	Such fields $X$ contain terms of the form $d_{1} \, \bar{w}^{c}\frac{\partial}{\partial \bar{y}_{s_{1}}}$ with $c\geq 2$ and terms of the form $d_{3} \,\tilde{p}_{s_{1}} w^{c}\frac{\partial}{\partial w}$. We shall prove that for $c>1$ such fields are zero. It is enough to consider the case $c=2$ (otherwise, instead of $X$, one can consider the field $X_{c-2}$). At the same time $a_{1}\leq 2,a_{2}=2,a_{3}=2,b_{2}=0,b_{3}=2$, and the degrees of the components of the vector polynomials $\tilde{h}_{1,\gamma}, \, \tilde{f}_{2,\alpha}, \, \tilde{g}_{1,\beta}, \, \tilde{g}_{3,\beta}$ are equal to two, three, zero, and four, respectively (if these components do not vanish). Without loss of generality we shall assume that $a_{2}=2$, while allowing $\tilde{f}_{1,\alpha}$ to be zero.
	
	We also have $\tilde{h}_{j,\gamma}=0$ for $j\neq 1$, since otherwise the lowest component of the tangency condition contains terms of the form $d z^{A}u^{C}$ with $|A|\geq 4$, and they can be obtained in a unique way (by means of $\tilde{h}_{j,\gamma}$). Next, $\tilde{f}_{j,\alpha}=0$ for $j\neq 1, j\neq 2$, since otherwise the lowest component of the tangency condition contains terms of the form $d z^{A}\bar{z}_{t}u^{C}$ with $|A|\geq 5$, and they can be obtained in a unique way (by means of $\tilde{f}_{j,\alpha}$).
	   
	   We write out the terms of bidegrees $(s+2,s)$ and $(s,s+2)$ in the variables $z,\bar{z}$ in the tangency condition:

	   $$2 \, {\rm Re} \, \Big(\tilde{h}_{1,\gamma}(u+i\langle z,\bar{z}\rangle)^{2}-\langle\tilde{f}_{2,\alpha},\bar{z}\rangle(u+i\langle z,\bar{z}\rangle)-$$
	   
	   $$-\overline{\mathcal{P}}\tilde{g}_{1,\beta}(u+i\langle z,\bar{z}\rangle)^{2}-\overline{\mathcal{P}}\tilde{g}_{3,\beta}\Big)=0.$$
	   
	   Separating in this equality the terms of degrees 2,1,0 in $u$, we obtain a system of equalities:
	   
	   $$- i\tilde{h}_{1,\gamma}+ i \overline{\tilde{h}_{1,\gamma}}-2 \overline{\mathcal{P}}\tilde{g}_{1,\beta} -2 \mathcal{P}\overline{\tilde{g}_{1,\beta}}=0,$$
	   
	   $$\Big(\tilde{h}_{1,\gamma} + \overline{\tilde{h}_{1,\gamma}} -2 i  \overline{\mathcal{P}}\tilde{g}_{1,\beta}+2i \mathcal{P}\overline{\tilde{g}_{1,\beta}}\Big) \langle z,\bar{z}\rangle-$$
	   
	   $$ -  \langle \tilde{f}_{2,\alpha},\bar{z}\rangle +  \langle \overline{\tilde{f}_{2,\alpha}},z\rangle=0,$$
	   
	   $$\Big(\frac{1}{2} i \tilde{h}_{1,\gamma} -\frac{1}{2} i \overline{\tilde{h}_{1,\gamma}} +  \overline{\mathcal{P}}\tilde{g}_{1,\beta}+ \mathcal{P}\overline{\tilde{g}_{1,\beta}}\Big) (\langle z,\bar{z}\rangle)^2+$$
	   
	   $$+\Big(- i \langle \tilde{f}_{2,\alpha},\bar{z}\rangle + i \langle \overline{\tilde{f}_{2,\alpha}},z\rangle\Big) \langle z,\bar{z}\rangle-\overline{\mathcal{P}}\tilde{g}_{3,\beta}-\mathcal{P}\overline{\tilde{g}_{3,\beta}}=0.$$
	   
	   Hence we obtain the equalities
	   
	   \begin{equation}\label{eq57}
	   \tilde{h}_{1,\gamma} = 2i \mathcal{P}\overline{\tilde{g}_{1,\beta}}, \langle \tilde{f}_{2,\alpha},\bar{z}\rangle = 4i \mathcal{P}\overline{\tilde{g}_{1,\beta}} \langle z,\bar{z}\rangle.
	   \end{equation}
	   
	   From the second equality it follows that $\tilde{f}_{2,\alpha} = (4i \mathcal{P}\overline{\tilde{g}_{1,\beta}}) z$.

	   Now we write out the coefficient of the field $\widehat{[X,X_{1}]}$ at $\frac{\partial}{\partial w}$. By Lemma 12 it is equal to
	   
	   \begin{equation}\label{eq56}
	   w^{2}\Big(-2\tilde{h}_{1,\gamma}^{2}
	   +\tilde{f}_{2,\alpha}\frac{\partial}{\partial z}(\tilde{h}_{1,\gamma})\Big),
	   \end{equation}
	   if this expression is nonzero.	   
	   Moreover, if expression \eqref{eq56} is nonzero, then the field $[X,X_{1}]$ has type 2.4(2), whence by what was proved in item 4.b expression \eqref{eq56} vanishes -- a contradiction. Hence, taking the equalities \eqref{eq57} into account, expression \eqref{eq56} takes the form
	   
	   $$w^{2}\Big(-2(2i \mathcal{P}\overline{\tilde{g}_{1,\beta}})^{2}
	   +(4i \mathcal{P}\overline{\tilde{g}_{1,\beta}}) z\frac{\partial}{\partial z}(2i \mathcal{P}\overline{\tilde{g}_{1,\beta}})\Big)=$$
	   
	   $$=w^{2}\Big(8( \mathcal{P}\overline{\tilde{g}_{1,\beta}})^{2}
	   -16(\mathcal{P}\overline{\tilde{g}_{1,\beta}})^{2})\Big)=-8w^{2}( \mathcal{P}\overline{\tilde{g}_{1,\beta}})^{2}
	   =0.$$
	   
	   We used the fact that $z\frac{\partial}{\partial z}(\mathcal{P})=2\mathcal{P}$, since $\mathcal{P}$ is a homogeneous vector polynomial of degree two. Thus $ \mathcal{P}\overline{\tilde{g}_{1,\beta}}=0$, and therefore from \eqref{eq57} we obtain $\tilde{h}_{1,\gamma}=0$ -- a contradiction. Hence such fields $X$ do not exist.
	   
	   \vspace{3ex}

	4.d) We prove that fields of type $2.1(\geq 3)$ do not exist. Such fields contain terms of the form $d \,z_{s_{1}} w^{c}\frac{\partial}{\partial w}$ and terms of the form $d_{1} \, w^{c}\frac{\partial}{\partial z_{s_{1}}}$.
	We prove that such fields do not exist for $c\geq 3$. One may assume that $c=3$ (otherwise, instead of $X$, one can consider the field $X_{c-3}$). In this case $a_{1}=a_{3}=3, \tilde{f}_{1,\alpha}$ is a constant vector, and $\tilde{h}_{1,\gamma}$ is linear in $z$, i.e. $b_{1}=0,b_{3}=1$.
	
	We have $\tilde{h}_{j,\gamma}=0$ for $j\neq 1$, since otherwise the lowest component of the tangency condition contains terms of the form $d z^{A}u^{C}$ with $|A|\geq 3$, and they can be obtained in a unique way (by means of $\tilde{h}_{j,\gamma}$). Next, $\tilde{f}_{j,\alpha}=0$ for $j\neq 1, j\neq 2$, since otherwise the lowest component of the tangency condition contains terms of the form $d z^{A}\bar{z}_{t}u^{C}$ with $|A|\geq 4$, and they can be obtained in a unique way (by means of $\tilde{f}_{j,\alpha}$). We also have that the degrees of the vector polynomial $\tilde{g}_{j,\beta}$ are odd if $g_{\beta}\neq 0$. Moreover, $\tilde{g}_{j,\beta}=0$ if the degree of the polynomial $\tilde{g}_{j,\beta}$ is different from one and three, since otherwise the lowest component of the tangency condition contains terms of the form $d z^{A}\bar{z}^{B}u^{C}$ with $|A|\geq 5, |B|=2$, and they can be obtained in a unique way (by means of $\tilde{g}_{j,\beta}$).
	
	At the same time $a_{2}\leq 2$, and without loss of generality one can assume that $g_{\beta}=\tilde{g}_{1,\beta}w^{2}+\tilde{g}_{2,\beta}w$, while allowing $\tilde{g}_{1,\beta}$ or $\tilde{g}_{2,\beta}$ to be zero. Then the lowest component of the tangency condition takes the form

	$$2 \, {\rm Re} \, \Big(\tilde{h}_{1,\gamma}(u+i\langle z,\bar{z}\rangle)^{3}-\langle\tilde{f}_{1,\alpha},\bar{z}\rangle(u+i\langle z,\bar{z}\rangle)^{3}-\langle\tilde{f}_{2,\alpha},\bar{z}\rangle(u+i\langle z,\bar{z}\rangle)^{2}-$$
	
	$$-\overline{\mathcal{P}}\tilde{g}_{1,\beta}(u+i\langle z,\bar{z}\rangle)^{2}-\overline{\mathcal{P}}\tilde{g}_{2,\beta}(u+i\langle z,\bar{z}\rangle)\Big).$$
	
	All terms in this equality have bidegrees $(s+1,s)$ or $(s,s+1)$ in the variables $z,\bar{z}$. Separating in this equality the terms of degrees 3,2,1,0 in $u$, we obtain a system of equalities:
	
	$$-3 i\tilde{h}_{1,\gamma}+3 i \overline{\tilde{h}_{1,\gamma}}-6 \langle \tilde{f}_{1,\alpha},\bar{z}\rangle -6 \langle \overline{\tilde{f}_{1,\alpha}},z\rangle=0,$$
	
	$$\Big(3 \tilde{h}_{1,\gamma} +3 \overline{\tilde{h}_{1,\gamma}} - 6 i \langle \tilde{f}_{1,\alpha},\bar{z}\rangle + 6 i \langle \overline{\tilde{f}_{1,\alpha}},z\rangle\Big) \langle z,\bar{z}\rangle-$$
	
	$$-2 \langle \tilde{f}_{2,\alpha},\bar{z}\rangle-2 \langle \overline{\tilde{f}_{2,\alpha}},z\rangle-2  \overline{\mathcal{P}}\tilde{g}_{1,\beta}-2 \mathcal{P}\overline{\tilde{g}_{1,\beta}}=0,$$
	
	$$\Big(\frac{3}{2} i \tilde{h}_{1,\gamma} -\frac{3}{2} i \overline{\tilde{h}_{1,\gamma}} +3 \langle \tilde{f}_{1,\alpha},\bar{z}\rangle +3 \langle \overline{\tilde{f}_{1,\alpha}},z\rangle\Big) (\langle z,\bar{z}\rangle)^2+$$
	
	$$+\Big(-2 i \langle \tilde{f}_{2,\alpha},\bar{z}\rangle +2 i \langle \overline{\tilde{f}_{2,\alpha}},z\rangle - 2 i \overline{\mathcal{P}}\tilde{g}_{1,\beta} +2 i \mathcal{P}\overline{\tilde{g}_{1,\beta}}\Big) \langle z,\bar{z}\rangle-\overline{\mathcal{P}}\tilde{g}_{2,\beta}-\mathcal{P}\overline{\tilde{g}_{2,\beta}}=0,$$
	
	$$\Big(-\frac{1}{2} \tilde{h}_{1,\gamma}-\frac{1}{2} \overline{\tilde{h}_{1,\gamma}} +i \langle \tilde{f}_{1,\alpha},\bar{z}\rangle -i \langle \overline{\tilde{f}_{1,\alpha}},z\rangle\Big) (\langle z,\bar{z}\rangle)^3+\Big(\langle \tilde{f}_{2,\alpha},\bar{z}\rangle +$$
	
	$$+\langle \overline{\tilde{f}_{2,\alpha}},z\rangle +\overline{\mathcal{P}}\tilde{g}_{1,\beta} +\mathcal{P}\overline{\tilde{g}_{1,\beta}}\Big) (\langle z,\bar{z}\rangle)^2+\Big(-i \overline{\mathcal{P}}\tilde{g}_{2,\beta} +i \mathcal{P}\overline{\tilde{g}_{2,\beta}}\Big) (\langle z,\bar{z}\rangle)=0.$$
	
	Simplifying this system, we obtain the equalities
	
	\begin{equation}\label{eq53}
	\tilde{h}_{1,\gamma} = 2i \langle \overline{\tilde{f}_{1,\alpha}},z\rangle, \mathcal{P}\overline{\tilde{g}_{1,\beta}} = 6i \langle \overline{\tilde{f}_{1,\alpha}},z\rangle \langle z,\bar{z}\rangle-\langle \tilde{f}_{2,\alpha},\bar{z}\rangle,
	\end{equation}
	
	$$\overline{\mathcal{P}}\tilde{g}_{2,\beta} = -12 \langle \overline{\tilde{f}_{1,\alpha}},z\rangle \langle z,\bar{z}\rangle^2-4i \langle \tilde{f}_{2,\alpha},\bar{z}\rangle \langle z,\bar{z}\rangle, 
	 \langle \tilde{f}_{2,\alpha},\bar{z}\rangle = 4i \langle \overline{\tilde{f}_{1,\alpha}},z\rangle \langle z,\bar{z}\rangle$$
	 and their conjugates. From the fourth equality we obtain $\tilde{f}_{2,\alpha} = 4i \langle \overline{\tilde{f}_{1,\alpha}},z\rangle z$
	 
	 Now we write out the coefficient of the field $\widehat{[X,X_{1}]}$ at $\frac{\partial}{\partial w}$. By Lemma 12 it is equal to
	 
	 \begin{equation}\label{eq54}
	 w^{4}\Big(-3\tilde{h}_{1,\gamma}^{2} +\tilde{f}_{2,\alpha} \frac{\partial}{\partial z}(\tilde{h}_{1,\gamma})\Big),
	 \end{equation}
	 if this expression is nonzero.
	 Moreover, if expression \eqref{eq54} is nonzero, then $[X,X_{1}]$ has type 2.2(4), whence by what was proved in item 4.c expression \eqref{eq54} vanishes -- a contradiction. Taking the equalities \eqref{eq53} into account, expression \eqref{eq54} takes the form
	 
	 $$w^{4}\Big(-3(2i \langle \overline{\tilde{f}_{1,\alpha}},z\rangle)^{2} +\tilde{f}_{2,\alpha} \frac{\partial}{\partial z}(2i \langle \overline{\tilde{f}_{1,\alpha}},z\rangle)\Big)=$$
	 
	 $$=w^{4}\Big(12(\langle \overline{\tilde{f}_{1,\alpha}},z\rangle)^{2} + 2i \langle \overline{\tilde{f}_{1,\alpha}},\tilde{f}_{2,\alpha}\rangle\Big)=$$
	
	$$=w^{4}\Big(12(\langle \overline{\tilde{f}_{1,\alpha}},z\rangle)^{2} + 2i \langle \overline{\tilde{f}_{1,\alpha}},4i \langle \overline{\tilde{f}_{1,\alpha}},z\rangle z\rangle\Big)=$$
	
	$$=w^{4}\Big(12(\langle \overline{\tilde{f}_{1,\alpha}},z\rangle)^{2} -8 \langle \overline{\tilde{f}_{1,\alpha}},z\rangle \langle \overline{\tilde{f}_{1,\alpha}}, z\rangle\Big)=4 w^{4}(\langle \overline{\tilde{f}_{1,\alpha}},z\rangle)^{2}=0.$$
	
	Hence $\langle \overline{\tilde{f}_{1,\alpha}},z\rangle= 0$, and therefore from \eqref{eq53} we find $\tilde{h}_{1,\gamma}= 0$. The contradiction obtained means that such fields do not exist.

	 \vspace{3ex}

	 4.e) Now we estimate the number of remaining fields.
	 
	 Consider four subcases.
	 
	 \vspace{3ex}
	 
	 4.e.1) Let the field have type $2.b(0)$. Such fields do not lie in the stabilizer of the automorphism algebra, since in this case either $h_{1,\gamma}$ (for $b=0$), or $f_{1,\alpha}$ (for $b=1$), or $g_{1,\beta}$ (for $b=2$) is equal to a constant vector; hence the field does not vanish at the origin. Moreover, for $b\geq 3$ fields of this type do not exist, as shown in item 4.b.
	 
	 \vspace{3ex}
	 
	 4.e.2) Let the field have type $2.b(1)$.
	 
	 First let $b=0$, i.e. $\tilde{h}_{1,\gamma}=const$. Separating in the lowest component of the tangency condition the monomials of the form $d \cdot u$, we obtain that $\tilde{h}_{1,\gamma}\in \mathbb{R}$. There is no more than one such field.
	 
	 Now let $b>1$, i.e. $\tilde{h}_{1,\gamma}\neq const$. Separating in the lowest component of the tangency condition the monomials of the form $d \cdot z^{a} u$, we obtain that $\tilde{h}_{1,\gamma}$ contains a linear combination of monomials $\tilde{p}_{j}w, \, 1\leq j\leq k,$ and monomials $z_{j}w, \, 1\leq j\leq m,$ with complex coefficients. Moreover, in this case the components $\varphi_{j}$ of the vector $\tilde{f}_{1,\alpha}$ contain terms $A_{j}w, \, A_{j}\in \mathbb{C}$, and the components $\psi_{j}$ of the vector $\tilde{g}_{1,\beta}$ contain terms $B_{j}w, \, B_{j}\in \mathbb{C}$. We also note that $\tilde{h}_{1,\gamma}\neq 0$ if and only if at least one of the numbers $A_{1},...,A_{m},B_{1},...,B_{k}$ is nonzero. That is, there are no more than $2n$ fields of type $2.b(1)$: no more than $2m$ fields of type 2.1(1), no more than $2k$ fields of type 2.2(1), and fields of type $2.b(1)$ with $b\geq 3$ do not exist.
	 
	 \vspace{3ex}
	 
	 4.e.3) Let the field have type $2.0(2)$. In this case $\tilde{h}_{1,\gamma}=const$. Separating in the lowest component of the tangency condition the monomials of the form $d \cdot u^{2}$, we obtain that $\tilde{h}_{1,\gamma}\in \mathbb{R}$. There is no more than one such field.
	 
	 \vspace{3ex}
	 
	 4.e.4) Let the field have type $2.1(2)$. Then $\tilde{h}_{1,\gamma}$ contains a linear combination of monomials $z_{j}w^{2}, \, 1\leq j\leq m$ with complex coefficients. There are no more than $2m$ such fields.

	 \vspace{3ex}

	 5) We prove that the total number of linearly independent fields of types 1(0), 1.1(1), 1.1(2), 2.1(1), 2.1(2), and 2.2(1) does not exceed $n^{2}+2n-1$. Denote by $Y$ a field of type 1.1(1), by $W$ a field of type 1.1(2), by $Z^{1},Z^{2}$ fields of types 2.1(1) and 2.1(2), respectively, by $T$ a field of type 1.3(0), and by $Z^{3}$ a field of type 2.2(1). Let $S_{1},S_{2}$ be the number of linearly independent fields $Y$ and $W$, respectively.

	 The coefficients $\zeta_{1}$ and $\zeta_{2}$ at $\frac{\partial}{\partial w}$ of the fields $\hat{Z}^{1},\hat{Z}^{2}$ are sums of expressions $d_{j} z_{j} w+...$ and $d_{j} z_{j} w^{2}+...$, respectively, where $d_{j} \in\mathbb{C}$ and dots here and below denote terms of lower degree in the variable $w$. The coefficients $\chi_{s}$ and $\upsilon_{s}$ at $\frac{\partial}{\partial z_{s}}$ of the fields $\hat{W},\hat{Y}$ are sums of expressions $i c_{j} z_{j} w^{2}+...$ and $i c_{j} z_{j} w+..., \, c_{j}\in \mathbb{R}$, respectively (we assume that these fields have been diagonalized;
	 diagonalization is a unitary transformation of the variable $z$ which does not change the Hermitian form $\langle z,\bar{z}\rangle$ and can change the polynomials $p_{j}$, but the equation still remains written in normal form, and the type of the field is preserved under diagonalization; therefore considering diagonalized fields does not restrict generality). By Lemma 12, the coefficient at $\frac{\partial}{\partial w}$ of the field $\widehat{[W,Z^{1}]}$ is equal to $\sum_{s}\chi_{s}\frac{\partial}{\partial z_{s}}(\zeta_{1})$ if this expression is nonzero. But fields of type 2.1(3) have such a coefficient at $\frac{\partial}{\partial w}$, and they are zero by what was proved in item 4.d. Therefore the coefficient at the highest degree of the variable $w$ (i.e. at $w^{3}$) in the expression $\sum_{s}\chi_{s}\frac{\partial}{\partial z_{s}}(\zeta_{1})$ is equal to zero. Also, the coefficient at $\frac{\partial}{\partial w}$ of the field $\widehat{[Y,Z^{2}]}$ is equal to $\sum_{s}\upsilon_{s}\frac{\partial}{\partial z_{s}}(\zeta_{2})$ if this expression is nonzero. But fields of type 2.1(3) have such a coefficient at $\frac{\partial}{\partial w}$, and they are zero by what was proved in item 4.d. Therefore the coefficient at the highest degree of the variable $w$ (i.e. at $w^{3}$) in the expression $\sum_{s}\upsilon_{s}\frac{\partial}{\partial z_{s}}(\zeta_{2})$ is equal to zero. Hence, if a term $i c_{j} z_{j}w^{2}\frac{\partial}{\partial z_{j}}$ is present in the field $W$, the corresponding term $d_{j} z_{j}w\frac{\partial}{\partial w}$ in the field $Z^{1}$ is equal to zero; and if a term $i c_{j} z_{j}w\frac{\partial}{\partial z_{j}}$ is present in the field $Y$, the corresponding term $d_{j} z_{j}w^{2}\frac{\partial}{\partial w}$ in the field $Z^{2}$ is equal to zero. That is, the coefficients at $w^{2}$ of the expressions $\chi_{s}$ and at $w$ of the expressions $\zeta_{1}$ depend on disjoint sets of variables $z_{j}$, and also the coefficients at $w$ of the expressions $\upsilon_{s}$ and at $w^{2}$ of the expressions $\zeta_{2}$ depend on disjoint sets of variables $z_{j}$. Thus the number of linearly independent fields $Z^{1}$ and $Z^{2}$ does not exceed $2m-2S_{2}$ and $2m-2S_{1}$, respectively. Hence, in particular, we obtain that the total number of linearly independent fields $Z^{1}$ and $W$ does not exceed $2m-2S_{2}+S_{2}=2m-S_{2}\leq 2m$.
	 
	 Let the fields $\hat{Z}^{2}_{2}$ depend on the variables $z_{j}$ with $j\in J_{1}$. Then the fields $\hat{Y}_{1}$ can depend only on the variables $z_{j}, \, j\in J_{3}=\{1,...,m\}\setminus J_{1}$. That is, for the field $Y$ we have $\upsilon_{j}\frac{\partial}{\partial z_{j}}=(i c_{j} z_{j} w+...)\frac{\partial}{\partial z_{j}},$ with $c_{j}=0$ for $j\in J_{1}$. Next, let $T$ be a field of type 1.3(0). Its coefficient $\varphi_{s}$ at $\frac{\partial}{\partial z_{s}}$ is equal to a sum of expressions of the form $e_{j} p_{j}, \, e_{j} \in \mathbb{C}$. The coefficient at $\frac{\partial}{\partial z_{j}}$ of the field $\widehat{[T,Y]}$ is equal to $(ic_{j}w\varphi_{j}-\sum_{s\in J_{3}}\upsilon_{s}\frac{\partial}{\partial z_{s}}(\varphi_{j}))+...=iw(c_{j}\varphi_{j}-\sum_{s\in J_{3}}c_{s}z_{s}\frac{\partial}{\partial z_{s}}(\varphi_{j}))+...$ by Lemma 12, if this expression is nonzero. Suppose that this expression is nonzero. It has the form $d p w+...$, where $p$ is a homogeneous polynomial of degree two, but fields of type $1.3$, considered in item 2.c.3, have such a coefficient at $\frac{\partial}{\partial z_{s}}$, and they are zero by what was proved in the same item. Moreover, the coefficient at $\frac{\partial}{\partial w}$ of the field $\widehat{[T,Y]}$ is equal to zero by Lemma 12, since the coefficients at $\frac{\partial}{\partial w}$ of the fields $\hat{T},\hat{Y}$ are zero, i.e. $[T,Y]$ is a field of the first type. And the coefficient at $\frac{\partial}{\partial y}$ of the field $\widehat{[T,Y]}$ has weight three, since the coefficient at $\frac{\partial}{\partial z}$ of this field has weight four; hence $G$ depends on $w$ at most linearly. That is, $[T,Y]$ is a field of type 1.3. Therefore our assumption is false, and $d p w=0$. Hence $iw(c_{j}\varphi_{j}-\sum_{s\in J_{3}}c_{s}z_{s}\frac{\partial}{\partial z_{s}}(\varphi_{j}))=0$, whence $c_{j}\varphi_{j}-\sum_{s\in J_{3}}c_{s}z_{s}\frac{\partial}{\partial z_{s}}(\varphi_{j})=0$.
	 
	 Next, the coefficient at $\frac{\partial}{\partial w}$ of the field $\widehat{[T,Z^{2}]}$ is equal to $\sum_{s}\varphi_{s}\frac{\partial}{\partial z_{s}}(\zeta_{2})$ by Lemma 12, if this expression is nonzero. $\sum_{s}\varphi_{s}\frac{\partial}{\partial z_{s}}(\zeta_{2})$ is an expression of the form $d p w^{2}+...$, where $p$ is a homogeneous polynomial of degree two, but fields of type $2.2(2)$ have such a coefficient at $\frac{\partial}{\partial w}$, and they are zero by what was proved in item 4.c. Hence $d p w^{2}=0$, i.e. the coefficient at $w^{2}$ in the expression $\sum_{s}\varphi_{s}\frac{\partial}{\partial z_{s}}(\zeta_{2})$ is equal to zero, i.e. $\frac{\partial^{2}}{\partial w^{2}}\Big(\sum_{s}\varphi_{s}\frac{\partial}{\partial z_{s}}(\zeta_{2})\Big)=0$.

	 Let $|Z^{2}|$ be the number of linearly independent fields $Z^{2}$. Choosing suitable linear combinations of the fields $Z^{2}$, one can assume that each of the fields $Z^{2}$ contains a term $\sigma$ of the form $\delta_{j}z_{j}w^{2}\frac{\partial}{\partial w}$ or $i \delta_{j}z_{j}w^{2}\frac{\partial}{\partial w}$ with $\delta_{j}\in \mathbb{R}$, which does not enter the other fields $Z^{2}$. Moreover, we can successively choose the terms $\sigma$ as follows. Let the first term $\sigma _{1}$ be equal to $\delta_{j_{1}}z_{j_{1}}w^{2}\frac{\partial}{\partial w}$ or $i\delta_{j_{1}}z_{j_{1}}w^{2}\frac{\partial}{\partial w}$. If $\sigma_{1}$ is present in other fields $Z^{2}$, we pass to linear combinations of the fields $Z^{2}$ such that the other fields do not contain the term $\sigma_{1}$. Then, if the term paired to $\sigma_{1}$, of the form $i \delta_{j_{1}}'z_{j_{1}}w^{2}\frac{\partial}{\partial w}$ or $\delta_{j_{1}}'z_{j_{1}}w^{2}\frac{\partial}{\partial w}$, respectively, enters some other field $Z^{2}$, we choose $\sigma_{2}$ as the second term. If there are no paired terms in the other fields $Z^{2}$, then we choose any other term as $\sigma_{2}$. And so on: each next term in the remaining fields is chosen either as paired to the previous one or arbitrarily, after first removing the previous term from the remaining fields (passing to a linear combination of the fields $Z^{2}$). Let the terms $\sigma$ depend on $z_{j}$ for $j\in J_{2} \subset J_{1}$. The terms $\sigma$ depend on no fewer than $\frac{|Z^{2}|}{2}$ variables $z_{j}$ if $|Z^{2}|$ is even, and on no fewer than $\frac{|Z^{2}|+1}{2}$ variables $z_{j}$ if $|Z^{2}|$ is odd. That is, the cardinality of the set $J_{2}$ is at least $\frac{|Z^{2}|}{2}$ or $\frac{|Z^{2}|+1}{2}$, respectively. Hence there are at least $\frac{|Z^{2}|}{2}$ complex linearly independent equalities of the form $\frac{\partial^{2}}{\partial w^{2}}\Big(\sum_{s}\varphi_{s}\frac{\partial}{\partial z_{s}}(\zeta_{2})\Big)=0$. Therefore the number of free coefficients $\varphi_{s}, \, s\in J_{1}\setminus J_{2},$ through which the remaining coefficients $\varphi_{s}$ are expressed does not exceed $m-|J_{2}|\leq m-\frac{|Z^{2}|}{2}$, where $|\cdot|$ denotes the cardinality of a set. The remaining $\varphi_{s}, \, s\in J_{2}$ are uniquely recovered from $\varphi_{s}, \, s\in J_{1}\setminus J_{2}$. Moreover, the equalities $\frac{\partial^{2}}{\partial w^{2}}\Big(\sum_{s}\varphi_{s}\frac{\partial}{\partial z_{s}}(\zeta_{2})\Big)=0$ impose no restrictions on $\varphi_{s}$ for $s\in J_{3}$, and hence each such $\varphi_{s}$ can take $2k$ real linearly independent values. Therefore altogether there are no more than $2k(|J_{1}|+|J_{3}|-\frac{|Z^{2}|}{2})=2k(m-\frac{|Z^{2}|}{2})=2km-|Z^{2}|k$ linearly independent fields $T$ (since each coefficient $\varphi_{s}$ is a sum of no more than $2k$ real linearly independent polynomials).

	 At the same time four cases are possible.
	 
	 \vspace{3ex}
	 
	 5.1) There are no fields $Y$ (fields of type 1.1(1)).
	 
	 Recall that the total number of fields $W$ of type 1.1(2) and fields $Z^{1}$ of type 2.1(1) does not exceed $2m$. Then the total number of fields of types 1(0), 1.1(1), 1.1(2), 2.1(1), 2.1(2), and 2.2(1) does not exceed $(m^{2}+k^{2}+2km-|Z^{2}|k)+0+2m+|Z^{2}|+2k\leq n^{2}+(1-k)|Z^{2}|+2n$.
	 
	 5.1.1) If $|Z^{2}|>0$ and $k>1$, then this quantity does not exceed $n^{2}+2n-1$.
	 
	 5.1.2) If $|Z^{2}|=0$ and $k\geq 1$, then we use the estimate $2km-2k$ for the number of fields of type 1.3(0) (see item 1.c). Therefore the total number of fields of types 1(0), 1.1(1), 1.1(2), 2.1(1), 2.1(2), and 2.2(1) does not exceed $(m^{2}+k^{2}+2km-2k)+0+2m+0+2k\leq n^{2}+2n-2$.
	 
	 5.1.3) If $m>1$, $|Z^{2}|>0$, and $k= 1$, then we prove that there is not a maximal number of fields of type 1.1(0) (i.e. their number does not exceed $m^{2}-1$).	 
	 Suppose that the number of fields of type 1.1(0) (from case 1.a) is equal to $m^{2}$ (i.e. maximal). Then the fields $2 \, {\rm Re} \, (iz_{j}\frac{\partial}{\partial z_{j}}), \, 2 \, {\rm Re} \, (z_{j}\frac{\partial}{\partial z_{s}}-z_{s}\frac{\partial}{\partial z_{j}})$ and $2 \, {\rm Re} \, (iz_{j}\frac{\partial}{\partial z_{s}}+iz_{s}\frac{\partial}{\partial z_{j}})$ enter the supporting components of different fields from the automorphism algebra.	 
	 If the field $2 \, {\rm Re} \, (f\frac{\partial}{\partial z})$ enters the supporting component of some field of type 1.1(0) from the automorphism algebra, then the expression $f\frac{\partial}{\partial z}(p_{1})$ must either be equal to zero or be proportional to the polynomial $p_{1}$ with some nonzero complex coefficient. (Indeed, otherwise a term $f\frac{\partial}{\partial z}(p_{1})\bar{y}_{1}$ appears in the tangency condition, and a term proportional to it cannot be obtained in any other way except by means of the field $f\frac{\partial}{\partial z}$ and the term $p_{1}\bar{y}_{1}$ in the defining relation; hence the tangency condition cannot be fulfilled.) And then the field contains a term of the form $d \,y_{1} \frac{\partial}{\partial y_{1}}$ (with its help one can obtain the term $-f\frac{\partial}{\partial z}(p_{1})\bar{y}_{1}$).	 
	 Then $(z_{j}\frac{\partial}{\partial z_{s}}-z_{s}\frac{\partial}{\partial z_{j}})(p_{1})=d_{1}p_{1}, \, (iz_{j}\frac{\partial}{\partial z_{s}}+iz_{s}\frac{\partial}{\partial z_{j}})(p_{1})=d_{2}p_{1}$, whence $(z_{j}\frac{\partial}{\partial z_{s}})(p_{1})=d_{3}p_{1}$. Also $(iz_{s}\frac{\partial}{\partial z_{s}})(p_{1})=d_{4}p_{1}$ holds. It follows that $\frac{\partial}{\partial z_{s}}(p_{1})=0$, i.e. $p_{1}$ does not depend on $z_{s}$. By the arbitrariness of $s$, the polynomial $p_{1}$ does not depend on any variable -- a contradiction. Hence the number of fields from case 1.a is not maximal, as required. Therefore the total number of fields of types 1(0), 1.1(1), 1.1(2), 2.1(1), 2.1(2), and 2.2(1) does not exceed $(m^{2}-1+k^{2}+2km-|Z^{2}|k)+0+2m+|Z^{2}|+2k\leq n^{2}-1+2n$.
	 
	 5.1.4) If, however, $m=1$, $|Z^{2}|>0$, and $k= 1$, then we use the estimate $2km-2k=0$ for the number of fields of type 1.3(0) (see item 1.c). Thus, the fields of type 1.3(0) do not exist (the condition $k= 1$ is automatically satisfied for $m=1$).

	 We prove that fields $Z^{3}$ of type 2.2(1) do not exist. The coefficient $\zeta_{3}$ at $\frac{\partial}{\partial w}$ of the field $\hat{Z}^{3}$ of type 2.2(1) is equal to an expression of the form $d z^{2} w+..., \, d\in \mathbb{C}$. And the coefficient $\zeta$ at $\frac{\partial}{\partial z}$ of the field $\hat{Z}^{2}$ has the form $d' w^{2}+...$, where $d' \in\mathbb{C}$. Then the coefficient at $\frac{\partial}{\partial w}$ of the field $\widehat{[Z^{2},Z^{3}]}$ is equal to $2 d d' w^{3} z+...$ by Lemma 12, if this expression is nonzero. But fields of type 2.1(3) have such a coefficient at $\frac{\partial}{\partial w}$, and they are zero by what was proved in item 4.d. Therefore $2 d d' w^{3} z=0$, whence $d=0$ (since $d'\neq 0$ by the condition $|Z^{2}|>0$). Hence fields of type 2.2(1) do not exist. Therefore the total number of fields of types 1(0), 1.1(1), 1.1(2), 2.1(1), 2.1(2), and 2.2(1) does not exceed $(m^{2}+k^{2})+0+2m+|Z^{2}|+0\leq n^{2}-2km+2n-2k+2\leq n^{2}+2n-2$.
	 
	 \vspace{3ex}
	 
	 5.2) There exists exactly one field $Y$ (up to proportionality) such that $\hat{Y}_{1}=2 \, {\rm Re} \, \Big(\sum_{j=1}^{m}iz_{j}\frac{\partial}{\partial z_{j}}\Big)$. In this case $|Z^{2}|=0$, since the field $\hat{Y}_{1}$ depends on all variables $z_{j}$.
	 
	 Since the nonzero coefficients $\varphi_{j}$ of the field $T$ are homogeneous quadratic polynomials, considering the field $[T,Y]$ we obtain $c_{j}\varphi_{j}-\sum_{s\in J_{3}}c_{s}z_{s}\frac{\partial}{\partial z_{s}}(\varphi_{j})=\varphi_{j}-2\varphi_{j}=-\varphi_{j}=0$ (the coefficients $c_{j}$ are equal to one). That is, fields $T$ of type 1.3(0) do not exist.
	 
	 Next, the coefficient $\zeta_{3}$ at $\frac{\partial}{\partial w}$ of the field $\hat{Z}^{3}$ of type 2.2(1) is equal to a sum of expressions of the form $d_{j} p_{j} w+..., \, d_{j}\in \mathbb{C}$. And the coefficient at $\frac{\partial}{\partial w}$ of the field $\widehat{[Y,Z^{3}]}$ is equal to $\hat{Y}(\zeta_{3})=2iw \zeta_{3}+...$ by Lemma 12, if this expression is nonzero. The expression $2iw \zeta_{3}+...$ has the form $d p w^{2}+...$, where $p$ is a homogeneous polynomial of degree two, but fields of type $2.2(2)$ have such a coefficient at $\frac{\partial}{\partial w}$, and they are zero by what was proved in item 4.c. Hence the coefficient at the highest degree of the variable $w$ (i.e. at $w$) in the expression $\zeta_{3}$ is equal to zero; consequently, fields of type 2.2(1) do not exist. Therefore the total number of fields of types 1(0), 1.1(1), 1.1(2), 2.1(1), 2.1(2), and 2.2(1) does not exceed $(m^{2}+k^{2})+1+2m+0+0\leq n^{2}-2km+2n-2k+1\leq n^{2}+2n-3$.
	 
	 \vspace{3ex}

	 5.3) There exists exactly one field $Y$ (up to proportionality) such that $\hat{Y}_{1}\neq 2 \, {\rm Re} \, \Big(\sum_{j=1}^{m}iz_{j}\frac{\partial}{\partial z_{j}}\Big)$.
	 
	 This means that there exist $j,j'$ such that $c_{j}\neq c_{j'}$. We prove that the number of linearly independent fields $Z^{4}$ of type 1.1(0) does not exceed $m^{2}-2$. Suppose first that there exists a field $Z^{4}$ such that the field $\hat{Z}^{4}$ is equal to $2 \, {\rm Re} \, (z_{j}\frac{\partial}{\partial z_{j'}}-z_{j'}\frac{\partial}{\partial z_{j}})$. The coefficient at $\frac{\partial}{\partial z}$ of the field $\widehat{[Y,Z^{4}]}$ is equal to $(\sum_{s=1}^{m}\upsilon_{s}\frac{\partial}{\partial z_{s}})(\hat{Z}^{4})-\hat{Z}^{4}(\sum_{s=1}^{m}\upsilon_{s}\frac{\partial}{\partial z_{s}})$ by Lemma 12, if this expression is nonzero. The coefficient at $\frac{\partial}{\partial z_{j}}$ of the field $\widehat{[Y,Z^{4}]}$ is equal to $iwz_{j'}(c_{j}-c_{j'})+...$. And the coefficient at $\frac{\partial}{\partial z_{j'}}$ of the field $\widehat{[Y,Z^{4}]}$ is equal to $iwz_{j}(c_{j}-c_{j'})+...$. That is, the fields $\widehat{[Y,Z^{4}]}$ and $\hat{Y}$ are linearly independent -- a contradiction. Hence the field $Z^{4}$ does not lie in the automorphism algebra. Next, suppose that there exists a field $Z^{5}$ such that the field $\hat{Z}^{5}$ is equal to $2 \, {\rm Re} \, (iz_{j}\frac{\partial}{\partial z_{j'}}+iz_{j'}\frac{\partial}{\partial z_{j}})$. The coefficient at $\frac{\partial}{\partial z}$ of the field $\widehat{[Y,Z^{5}]}$ is equal to $(\sum_{s=1}^{m}\upsilon_{s}\frac{\partial}{\partial z_{s}})(\hat{Z}^{5})-\hat{Z}^{5}(\sum_{s=1}^{m}\upsilon_{s}\frac{\partial}{\partial z_{s}})$ by Lemma 12, if this expression is nonzero. The coefficient at $\frac{\partial}{\partial z_{j}}$ of the field $\widehat{[Y,Z^{5}]}$ is equal to $wz_{j'}(c_{j}-c_{j'})+...$. And the coefficient at $\frac{\partial}{\partial z_{j'}}$ of the field $\widehat{[Y,Z^{5}]}$ is equal to $-wz_{j}(c_{j}-c_{j'})+...$. That is, the fields $\widehat{[Y,Z^{4}]}, \, \widehat{[Y,Z^{5}]}$, and $\hat{Y}$ are linearly independent. Hence the fields $\hat{Z}^{4},\hat{Z}^{5}$ can enter the supporting component of a field from the automorphism algebra only together with some other fields. Therefore the number of linearly independent fields of type 1.1(0) is at least two less than the maximal one, i.e. does not exceed $m^{2}-2$.
	 
	 In this case the total number of fields of types 1(0), 1.1(1), 1.1(2), 2.1(1), 2.1(2), and 2.2(1) does not exceed $(m^{2}-2+k^{2}+2km-|Z^{2}|k)+1+2m+|Z^{2}|+2k\leq n^{2}-1+(1-k)|Z^{2}|+2n\leq n^{2}-1+2n$.
	 
	 \vspace{3ex}
	 
	 5.4) The number $|Y|$ of linearly independent fields $Y$ is at least two.
	 
	 Hence, passing if necessary to a linear combination of the fields $Y$, one may assume that there exist $|Y|$ terms of the form $i c_{j}wz_{j}\frac{\partial}{\partial z_{j}}, \, j\in J_{4}\subset J_{3}, \, |J_{4}|=|Y|$, entering exactly one field $Y$. Then for all fields $Y$ except one, the coefficient $c_{j}, \, j\in J_{4}$, at the term $iwz_{j}\frac{\partial}{\partial z_{j}}$ is equal to zero.
	 
	 Consider the field $[T,Y]$. We have the equality $c_{j}\varphi_{j}-\sum_{s\in J_{3}}c_{s}z_{s}\frac{\partial}{\partial z_{s}}(\varphi_{j})=0$ for the coefficients $\varphi_{j}$ of the field $T$. The polynomial $\varphi_{j}$ can be a solution of this equation either only for zero $c_{j}$ or only for nonzero $c_{j}$. Fix a number $t, \, 1\leq t\leq k,$ and check whether the equality $\varphi_{j}=d_{t}p_{t}, \, d_{t}\in \mathbb{C}$ can hold for some $j\in J_{4}$. Two cases are possible.
	 
	 \vspace{3ex}
	 
	 5.4.1) The equality $c_{j}(d_{t}p_{t})-\sum_{s\in J_{3}}c_{s}z_{s}\frac{\partial}{\partial z_{s}}(d_{t}p_{t})=0$ is not fulfilled for a nonzero value of $c_{j}$ for any fixed $j\in J_{4}$.
	 
	 For each field $Y$ consider the set $J(Y)$ of those $j\in J_{4}$ for which $c_{j}=0$. The intersection of the sets $J(Y)$ over all $Y$ is empty. Therefore $\varphi_{j}\neq d_{t}p_{t}$ for $j\in J_{4}$.	 
	 Hence the coefficient $\varphi_{j}, \, j\in J_{4}$, can take no more than $(2k-2)$ linearly independent values (since $2k$ is the maximal number of linearly independent values of the form $\sum_{s=1}^{k}d_{s}p_{s}$ that $\varphi_{j}$ can take). Therefore the number of linearly independent terms $\varphi_{s}\frac{\partial}{\partial z_{s}}, s\in J_{3}$, in fields $T$ of type 1.3(0) does not exceed $2k(|J_{3}|-|J_{4}|)+(2k-2)|J_{4}|=2k|J_{3}|-2|J_{4}|$. Thus we have improved the estimate $2k(|J_{1}|+|J_{3}|)-|Z^{2}|k)$ for the number of fields of type 1.3(0), obtained in item 5 above, as follows: the number of fields of type 1.3(0) does not exceed $2k(|J_{1}|+|J_{3}|)-|Z^{2}|k-2|J_{4}|$.
	 
	 In this case the total number of fields of types 1(0), 1.1(1), 1.1(2), 2.1(1), 2.1(2), and 2.2(1) does not exceed $(m^{2}+k^{2}+2k(|J_{1}|+|J_{3}|)-2|J_{4}|-|Z^{2}|k)+|J_{4}|+2m+|Z^{2}|+2k\leq n^{2}-|J_{4}|+(1-k)|Z^{2}|+2n\leq n^{2}-2+2n$.
	 
	 \vspace{3ex}
	 
	 5.4.2) The equality $c_{j}(d_{t}p_{t})-\sum_{s\in J_{3}}c_{s}z_{s}\frac{\partial}{\partial z_{s}}(d_{t}p_{t})=0$ is fulfilled for some nonzero value $c_{j}=\tilde{c}_{j}$ for some fixed $j\in J_{4}$.
	 
	 Fix a field $Y$ for which the coefficient $c_{j}$ is equal to $\tilde{c}_{j}$. In this case $\sum_{s\in J_{3}}c_{s}z_{s}\frac{\partial}{\partial z_{s}}(\varphi_{j})\neq 0$ for the coefficients $c_{s}$ of the fixed field $Y$. Next, suppose that the coefficient $\zeta_{3}$ at $\frac{\partial}{\partial w}$ of the field $\hat{Z}^{3}$ of type 2.2(1) is equal to an expression of the form $d_{t} p_{t} w+...$ . Then the coefficient at $\frac{\partial}{\partial w}$ of the field $\widehat{[Y,Z^{3}]}$ is equal to $\hat{Y}(\zeta_{3})=iw^{2}\sum_{s\in J_{3}}c_{s}z_{s}\frac{\partial}{\partial z_{s}}( d_{t} p_{t})+...$ by Lemma 12, if this expression is nonzero. The expression $iw^{2}\sum_{s\in J_{3}}c_{s}z_{s}\frac{\partial}{\partial z_{s}}( d_{t} p_{t})+...$ has the form $d p w^{2}+...$, where $p$ is a homogeneous polynomial of degree two, but fields of type $2.2(2)$ have such a coefficient at $\frac{\partial}{\partial w}$, and they are zero by what was proved in item 4.c. Hence the coefficient at the highest degree of the variable $w$ (i.e. at $w^{2}$) in the expression $iw^{2}\sum_{s\in J_{3}}c_{s}z_{s}\frac{\partial}{\partial z_{s}}( d_{t} p_{t})+...$ is equal to zero, whence $\sum_{s\in J_{3}}c_{s}z_{s}\frac{\partial}{\partial z_{s}}( d_{t} p_{t})=0$, which is impossible by our assumption. Consequently, $\zeta_{3}\neq d_{t} p_{t}w+...$, and therefore there are no more than $2k-2$ fields of type 2.2(1) (since $2k$ is the maximal number of linearly independent values of the form $\sum_{s=1}^{k}d_{s}p_{s}w+...$ that $\zeta_{3}$ can take).
	 
	 Next, for the fixed field $Y$, the equality $c_{j}(d_{t}p_{t})-\sum_{s\in J_{3}}c_{s}z_{s}\frac{\partial}{\partial z_{s}}(d_{t}p_{t})=0$ holds for some $j\in J_{4}$. Then for $s\neq j, \, s\in J_{4}$ this equality is not fulfilled, since $c_{s}=0$ for such values of $s$. Therefore $\varphi_{s}\neq d_{t}p_{t}$ for $s\in J_{4}, \, s\neq j$, and hence the coefficient $\varphi_{s}$ for such values of $s$ can take no more than $(2k-2)$ linearly independent values (since $2k$ is the maximal number of linearly independent values of the form $\sum_{s=1}^{k}d_{s}p_{s}$ that $\varphi_{s}$ can take). Hence the number of linearly independent terms $\varphi_{s}\frac{\partial}{\partial z_{s}}, s\in J_{3}$, in fields $T$ of type 1.3(0) does not exceed $2k(|J_{3}|-|J_{4}|)+(2k-2)(|J_{4}|-1)+2k=2k|J_{3}|-2|J_{4}|+2$. Thus we have improved the estimate $2k(|J_{1}|+|J_{3}|)-|Z^{2}|k)$ for the number of fields of type 1.3(0), obtained in item 5 above, as follows: the number of fields of type 1.3(0) does not exceed $2k(|J_{1}|+|J_{3}|)-|Z^{2}|k-2|J_{4}|+2$.
	 
	 In this case the total number of fields of types 1(0), 1.1(1), 1.1(2), 2.1(1), 2.1(2), and 2.2(1) does not exceed $(m^{2}+k^{2}+2k(|J_{1}|+|J_{3}|)-2|J_{4}|+2-|Z^{2}|k)+|J_{4}|+2m+|Z^{2}|+2k-2\leq n^{2}-|J_{4}|+(1-k)|Z^{2}|+2n\leq n^{2}-2+2n$.

	 \vspace{3ex}
	 
	 Taking item 4 of the proof into account, and adding the estimate two for the number of fields of types 2.0(1) and 2.0(2) from items 4.e.2 and 4.e.3, we obtain the estimate $n^{2}+2n+1$.

	Lemma 13 is proved.

	\vspace{3ex}

	Since at most $2n+1$ fields can lie in the complement to the stabilizer, the general estimate for the dimension of the full automorphism algebra is $n^{2}+4n+2$, which is strictly less than the corresponding dimension $n^{2}+4n+3$ for Levi-nondegenerate quadrics in $\mathbb{C}^{n+1}$.

	Thus we have proved the following.
	
	\vspace{3ex}
	
	\textbf{Theorem 14.} ${\rm dim \, aut}\,M_{0}<{\rm dim \, aut}\,Q$.

	\vspace{3ex}

	As a consequence, we also obtain the following assertion.

	\vspace{3ex}
	
	\textbf{Assertion 15.} The automorphism algebra of a hypersurface in $\mathbb{C}^{n+1}$ that is $2$-nondegenerate and has a sign-definite Levi form at a generic point is uniquely determined by its $3$-jet.

	\textbf{Proof.} From the proof of Lemma 13, in which the form of the supporting component of fields from ${\rm st} \, M_{0}$ is described, it is clear that the automorphism algebra of the germ of a hypersurface in $\mathbb{C}^{n+1}$ that is $2$-nondegenerate with sign-definite Levi form at a generic point is uniquely determined at a generic point on this hypersurface:
	
	1) by the derivatives $\frac{\partial^{|\delta|}}{\partial z^{\delta}}(f_{\alpha}), \, \frac{\partial^{|\delta|}}{\partial z^{\delta}}(g_{\beta})$ up to order $2$ inclusive (fields of types 1.1(0), 1.2(0), 1.3(0)),
	
	2) by the derivatives $\frac{\partial}{\partial w}(f_{\alpha}), \, \frac{\partial}{\partial w}(g_{\beta}), \, \frac{\partial}{\partial w}(h_{\gamma}), \, \frac{\partial^{2}}{\partial w^{2}}(h_{\gamma})$ (fields of types 2.1(1), 2.2(1), 2.0(1), 2.0(2)),
	
	3) by the derivatives $\frac{\partial^{2}}{\partial w \partial z_{j}}(f_{\alpha}), \, \frac{\partial^{2}}{\partial w^{2}}(f_{\alpha}), \, \frac{\partial^{3}}{\partial w^{2} \partial z_{j}}(f_{\alpha})$ (fields of types 1.1(1), 2.1(2), 1.1(2)).
	
	And the fields from the automorphism algebra lying in the complement to the stabilizer are determined by their value at the point.
	
	Assertion 15 is proved.

	\vspace{3ex}
	
	In other words, if at a generic point there are two automorphisms with the same sets of the indicated derivatives, then these two automorphisms coincide. Example 16 (see below) shows that the estimate for the jet order is not lower than two (the automorphism algebra of the hypersurface from this example contains the maximal number, i.e. $k^{2}$, of linearly independent fields of type 1.2(0)). At the same time the derivatives listed in items 1 and 2 of Assertion 15 are necessary for specifying an automorphism, whereas we do not know examples of hypersurfaces for whose automorphisms the derivatives from item 3 of the proof of Assertion 15 are required.
	
	\vspace{3ex}

	For $j=1$ the formulated result is known (and more precise: the estimate for the jet order is attained and is equal to two) and is contained in \cite{43}. We note, however, that we consider an arbitrary 2-nondegenerate hypersurface with a sign-definite Levi form at a generic point, whereas \cite{43} concerns the germ of a 1-nondegenerate hypersurface at an arbitrary point.
	And in \cite{30} an analogous result is obtained for hypersurfaces of finite multitype at an arbitrary point.
	
	\vspace{3ex}
	Together with the estimates for fields of each of the types, the following questions arise: are they sharp, and do there exist hypersurfaces with nontrivial tangent fields of the types considered? The answers to these questions are contained in the following series of examples, which shows, first, that there exist $2$-model surfaces realizing the maximum number of linearly independent fields of types 1.1(0), 1.2(0), and 1.3(0) for case 1.c.1.1, and second, that there exists a $2$-model surface whose automorphism algebra contains a nonzero field of type 1.3(0) for case 1.c.1.2.

	\vspace{3ex}
	
	\textbf{Example 16.} Let $k$ (the dimension of the vector variable $y$) be equal to the dimension of the space of homogeneous quadratic polynomials in $m$ variables with real coefficients (i.e. $k=\frac{m(m+1)}{2}$), and let $p_{j}$ be the basis quadratic monomials in $m$ variables. Then the $2$-model surface given by the equation $v=\langle z,\bar{z}\rangle+2 \, {\rm Re} \, (\mathcal{P}(z)\bar{y})$ in the notation adopted above (for any dimension $m$ of the variable $z$) realizes the maximum number of fields of types 1.1 and 1.2.
	
	Indeed, the automorphism algebra contains all fields $2 \, {\rm Re} \, (iz_{j}\frac{\partial}{\partial z_{j}})+L_{1}(y), \, 2 \, {\rm Re} \, (z_{\mu}\frac{\partial}{\partial z_{\nu}}-z_{\nu}\frac{\partial}{\partial z_{\mu}})+L_{2}(y), \, 2 \, {\rm Re} \, (iz_{\mu}\frac{\partial}{\partial z_{\nu}}+iz_{\nu}\frac{\partial}{\partial z_{\mu}})+L_{3}(y)$ (type 1.1(0), $m^{2}$ fields; here $L_{1}(y), \, L_{2}(y), \, L_{3}(y)$ are vector fields of the form $2 \, {\rm Re} \,\Big(\sum_{s,t}c_{s,t}y_{s}\frac{\partial}{\partial y_{t}}\Big), \, c_{s,t} \in \mathbb{C}$, which are uniquely recovered from the explicitly written terms of the tangent fields) and all fields $2 \, {\rm Re} \, (ip_{j}(z)\frac{\partial}{\partial y_{j}}), \, 2 \, {\rm Re} \, (p_{\mu}\frac{\partial}{\partial y_{\nu}}-p_{\nu}\frac{\partial}{\partial y_{\mu}}), \, 2 \, {\rm Re} \, (ip_{\mu}\frac{\partial}{\partial y_{\nu}}+ip_{\nu}\frac{\partial}{\partial y_{\mu}})$ (type 1.2(0), $k^{2}$ fields).
	
	At the same time it is clear that as $n$ grows, the dimension of the automorphism algebra of this hypersurface grows quadratically, i.e. as $O(n^{2})$. This happens because the dimension of the space of fields of the first type, which potentially make the largest contribution to the dimension of the algebra, is indeed large.
	
	\vspace{3ex}
	
	\textbf{Example 17.} Let $k$ be equal to the dimension of the space of homogeneous quadratic polynomials in $(m-r)$ variables with real coefficients (i.e. $k=\frac{(m-r)(m-r+1)}{2}$), and let $p_{j}$ be the basis quadratic monomials in the $(m-r)$ variables $z_{r+1},...,z_{m}$. Then the $2$-model surface given by the equation $v=\langle z,\bar{z}\rangle+2 \, {\rm Re} \, (\mathcal{P}(z)\bar{y})$ in the notation adopted above realizes the maximum number of fields of type 1.3(0) for case 1.c.1.1.
	
	Indeed, the automorphism algebra contains all fields $2 \, {\rm Re} \, (p_{\mu}\frac{\partial}{\partial z_{\nu}}-z_{\nu}\frac{\partial}{\partial y_{\mu}}), \, 2 \, {\rm Re} \, (ip_{\mu}\frac{\partial}{\partial z_{\nu}}+iz_{\nu}\frac{\partial}{\partial y_{\mu}})$ for $1\leq\nu \leq r$.

	\vspace{3ex}
	
	\textbf{Example 18.} The automorphism algebra of the $2$-model surface in $\mathbb{C}^{8}$ given by the equation $v=|z_{1}|^{2}+|z_{2}|^{2}+|z_{3}|^{2}+|z_{4}|^{2}+2 \, {\rm Re} \, ((z_{1}z_{3}+z_{2}z_{4})\bar{y}_{1}+z_{1}^{2}\bar{y}_{2}+z_{1}z_{2}\bar{y}_{3})$ contains a nontrivial field of type 1.3(0) for case 1.c.1.2. This field has the following form: $2 \, {\rm Re} \, \Big(z_{1}z_{2}\frac{\partial}{\partial z_{3}}-z_{1}^{2}\frac{\partial}{\partial z_{4}}-z_{3}\frac{\partial}{\partial y_{3}}+z_{4}\frac{\partial}{\partial y_{2}}\Big)$.

	\vspace{3ex}
	
	\textbf{Remark 19.} As can be seen from the proof of the dimension conjecture, the supporting component $\hat{X}$ of a field $X$ from the stabilizer of the automorphism algebra is uniquely determined by a linear combination of the terms entering $\hat{X}$ of the following form (we also indicate the type of the terms):
	
	1) type 1.1(0): $2\,{\rm Re}\Big(iz_{\mu}\frac{\partial}{\partial z_{\mu}}\Big), \, 2\,{\rm Re}\Big(z_{\mu}\frac{\partial}{\partial z_{\nu}}-z_{\nu}\frac{\partial}{\partial z_{\mu}}\Big), \, 2\,{\rm Re}\Big(iz_{\mu}\frac{\partial}{\partial z_{\nu}}+iz_{\nu}\frac{\partial}{\partial z_{\mu}}\Big)$;
	
	2) type 1.1(1): $2\,{\rm Re}\Big(iz_{\mu}\frac{\partial}{\partial z_{\mu}}\Big)w, \, 2\,{\rm Re}\Big(z_{\mu}\frac{\partial}{\partial z_{\nu}}-z_{\nu}\frac{\partial}{\partial z_{\mu}}\Big)w, \, 2\,{\rm Re}\Big(iz_{\mu}\frac{\partial}{\partial z_{\nu}}+iz_{\nu}\frac{\partial}{\partial z_{\mu}}\Big)w$;
	
	3) type 1.1(2): $2\,{\rm Re}\Big(iz_{\mu}\frac{\partial}{\partial z_{\mu}}\Big)w^{2}, \, 2\,{\rm Re}\Big(z_{\mu}\frac{\partial}{\partial z_{\nu}}-z_{\nu}\frac{\partial}{\partial z_{\mu}}\Big)w^{2}, \, 2\,{\rm Re}\Big(iz_{\mu}\frac{\partial}{\partial z_{\nu}}+iz_{\nu}\frac{\partial}{\partial z_{\mu}}\Big)w^{2}$;
	
	4) type 1.2(0): $2\,{\rm Re}\Big(p_{\nu}\frac{\partial}{\partial y_{\mu}}-p_{\mu}\frac{\partial}{\partial y_{\nu}}\Big), \, 2\,{\rm Re}\Big(ip_{\nu}\frac{\partial}{\partial y_{\mu}}+ip_{\mu}\frac{\partial}{\partial y_{\nu}}\Big)$;
	
	5) type 1.3(0): $2\,{\rm Re}\Big(p_{\nu}\frac{\partial}{\partial z_{\mu}}- z_{\mu}\frac{\partial}{\partial y_{\nu}}\Big), \, 2\,{\rm Re}\Big(ip_{\nu}\frac{\partial}{\partial z_{\mu}}+iz_{\mu}\frac{\partial}{\partial y_{\nu}}\Big)$;
	
	6) type 2.0(1): $2\,{\rm Re}\Big(w\frac{\partial}{\partial w}+...\Big)$, where here and below in items 7--10 dots denote fields with zero coefficient at $\frac{\partial}{\partial w}$;
	
	7) type 2.1(1): $2\,{\rm Re}\Big(z_{\mu}w\frac{\partial}{\partial w}+...\Big), \,  2\,{\rm Re}\Big(iz_{\mu}w\frac{\partial}{\partial w}+...\Big)$;
	
	8) type 2.1(2): $2\,{\rm Re}\Big(z_{\mu}w^{2}\frac{\partial}{\partial w}+...\Big), \,  2\,{\rm Re}\Big(iz_{\mu}w^{2}\frac{\partial}{\partial w}+...\Big)$;
	
	9) type 2.2(1): $2\,{\rm Re}\Big(p_{\mu}w\frac{\partial}{\partial w}+...\Big), \, 2\,{\rm Re}\Big(ip_{\mu}w\frac{\partial}{\partial w}+...\Big)$;
	
	10) type 2.0(2): $2\,{\rm Re}\Big(w^{2}\frac{\partial}{\partial w}+...\Big)$.
	
	At the same time, if the field contains the terms written in items 2, 3, and 8, then the number of linearly independent fields of other types decreases.
	
	Let
	$Q$ be the hyperquadric in $\mathbb{C}^{n+1}$ given by the equation $v=\langle z,\bar{z}\rangle+|y_{1}|^{2}+...+|y_{k}|^{2}$.
	
	Suppose that the terms from items 2, 3, and 8 are absent. Define a linear operator $\mathcal{L}: {\rm st} \, M_{0}\longrightarrow {\rm st} \,Q$ for each of the remaining items (by Assertion 15 it is enough to define it on the supporting component $\hat{X}$ of fields $X$ from the automorphism algebra) as follows.
	
	1), 4), 5): one must replace the polynomial $p_{\nu}$ everywhere by $y_{\nu}$. For example, for type 1.1(0) the action is the identity, and for one of the fields of type 1.2(0) we have $\mathcal{L}\Big(2\,{\rm Re}\Big(p_{\nu}\frac{\partial}{\partial y_{\mu}}-p_{\mu}\frac{\partial}{\partial y_{\nu}}\Big)\Big)=2\,{\rm Re}\Big(y_{\nu}\frac{\partial}{\partial y_{\mu}}-y_{\mu}\frac{\partial}{\partial y_{\nu}}\Big)$.

	6) type 2.0(1): $\mathcal{L}\Big(2\,{\rm Re}\Big(w\frac{\partial}{\partial w}+...\Big)\Big)=2\,{\rm Re}\Big(w\frac{\partial}{\partial w}+\frac{1}{2}\Big(z\frac{\partial}{\partial z}+y\frac{\partial}{\partial y}\Big)\Big)$;
	
	7) type 2.1(1): $\mathcal{L}\Big(2\,{\rm Re}\Big(z_{\mu}w\frac{\partial}{\partial w}+...\Big)\Big)=2\,{\rm Re}\Big(z_{\mu}w\frac{\partial}{\partial w}+z_{\mu}(z\frac{\partial}{\partial z}+y\frac{\partial}{\partial y})-\frac{1}{2i}w\frac{\partial}{\partial z_{\mu}}\Big), \, \mathcal{L}\Big(2\,{\rm Re}\Big(iz_{\mu}w\frac{\partial}{\partial w}+...\Big)\Big)=2\,{\rm Re}\Big(iz_{\mu}w\frac{\partial}{\partial w}+iz_{\mu}(z\frac{\partial}{\partial z}+y\frac{\partial}{\partial y})+\frac{1}{2}w\frac{\partial}{\partial z_{\mu}}\Big)$;
	
	9) type 2.2(1): $\mathcal{L}\Big(2\,{\rm Re}\Big(p_{\mu}w\frac{\partial}{\partial w}+...\Big)\Big)=2\,{\rm Re}\Big(y_{\mu}w\frac{\partial}{\partial w}+y_{\mu}(z\frac{\partial}{\partial z}+y\frac{\partial}{\partial y})-\frac{1}{2i}w\frac{\partial}{\partial y_{\mu}}\Big), \, \mathcal{L}\Big(2\,{\rm Re}\Big(ip_{\mu}w\frac{\partial}{\partial w}+...\Big)\Big)=2\,{\rm Re}\Big(iy_{\mu}w\frac{\partial}{\partial w}+iy_{\mu}(z\frac{\partial}{\partial z}+y\frac{\partial}{\partial y})+\frac{1}{2}w\frac{\partial}{\partial y_{\mu}}\Big)$;
	
	10) type 2.0(2): $\mathcal{L}\Big(2\,{\rm Re}\Big(w^{2}\frac{\partial}{\partial w}+...\Big)\Big)=2\,{\rm Re}\Big(w^{2}\frac{\partial}{\partial w}+w(z\frac{\partial}{\partial z}+y\frac{\partial}{\partial y})\Big)$.

	It is clear that the operator $\mathcal{L}$ so defined has no kernel. Thus we obtain a parametrization of the stabilizer of the algebra ${\rm aut} \, M_{0}$ by the stabilizer of the algebra ${\rm aut} \,Q$ written in \eqref{eq19}). That is, the parameter space of the stabilizer of the algebra ${\rm aut} \,M_{0}$ is a subspace in the parameter space of the algebra ${\rm aut} \,Q$. But since the hyperquadric is holomorphically homogeneous, the fields from the complement to the stabilizer are also parametrized in an obvious way by the corresponding fields from the algebra of the quadric. That is, we obtain a parametrization of the algebra ${\rm aut} \, M_{0}$ by the algebra ${\rm aut} \,Q$. It is also clear that instead of the indicated quadric $Q$, which is most natural in our situation, one can also take any other nondegenerate hyperquadric of the same CR dimension (the mapping will be arranged in an analogous way).

	A similar result on parametrization was obtained in \cite{42}, where the automorphism algebra of a uniformly 2-nondegenerate surface in $\mathbb{C}^{3}$ is parametrized by the automorphism algebra of a hyperquadric in $\mathbb{C}^{2}$.
	
	If the terms from items 2, 3, and 8 are present, then the operator $\mathcal{L}$ can be defined for each hypersurface so that linear combinations of terms from items 2, 3, and 8 replace linear combinations of terms from other items. However, it is worth noting that we do not know examples of hypersurfaces for which the terms from items 2, 3, and 8 are present in the automorphism algebra.
	
	We also emphasize that, since the stabilizer of the algebra ${\rm aut} \,Q$ is a linear space whose dimension exceeds the dimension of the algebra ${\rm aut} \,M_{0}$, the parametrization certainly exists even if it is not given by a concrete formula.
	
	\vspace{3ex}

	From the geometric point of view, the conjecture on the maximal symmetry of hyperquadrics in the class of holomorphically nondegenerate hypersurfaces looks quite natural. Namely, all holomorphically nondegenerate surfaces are divided into two groups:
	
	1) Levi-nondegenerate at a generic point,
	
	2) Levi-degenerate at a generic point.
	
	The first group consists exactly of hypersurfaces that are $1$-nondegenerate at a generic point, and the second group contains all surfaces that are $l$-nondegenerate at a generic point for $l\geq 2$. At the same time Levi-degenerate hypersurfaces have a nontrivial kernel of the Levi form, which imposes additional restrictions on automorphisms of such surfaces, since the kernel is biholomorphically invariant.
	
	Therefore one may expect the results obtained to have a geometric justification, along with the analytic proof obtained by us, which uses a certain special form of the defining equation -- the reduced form.

	Let us formulate the corresponding question:
	
	\vspace{3ex}
	
	\textbf{Question 20.} Give a geometric proof of the dimension conjecture for a 2-nondegenerate hypersurface with a sign-definite Levi form.

	\vspace{3ex}
	
	However, it is worth noting that a Levi-nondegenerate surface need not be more symmetric than a Levi-degenerate one: it is not difficult to give examples of Levi-nondegenerate hypersurfaces with trivial automorphism algebra and Levi-degenerate hypersurfaces with large algebras (see Example 16).

	\vspace{3ex}

	A question about the exact estimate for the dimension of the automorphism algebra is also appropriate:
	
	\textbf{Question 21.} Is the found estimate $n^{2}+4n+2$ for the dimension of the automorphism algebra of 2-nondegenerate hypersurfaces with a sign-definite Levi form attained? If so, on which hypersurfaces?

	\vspace{3ex}

	In conclusion, we formulate the remaining open question for 2-nondegenerate germs of arbitrary signature of the Levi form and for arbitrary $l$-nondegenerate germs.
	
	\vspace{3ex}
	
	\textbf{Question 22.} Is the dimension conjecture true for 2-nondegenerate germs of arbitrary signature of the Levi form and for $l$-nondegenerate hypersurfaces with $l>2$?

\end{document}